\documentclass[11pt]{amsart}
\usepackage{amsmath,latexsym,amsfonts,amssymb,amsthm,mathrsfs,longtable,lscape,dsfont,fancybox,yfonts}
\usepackage{amsthm}
\usepackage{amssymb}
\usepackage{amscd}
\usepackage{amsfonts}
\usepackage{amsbsy}
\usepackage{fancyhdr}
\usepackage{graphicx,color}
\usepackage[all]{xy}
\usepackage{epsfig}
\usepackage{subfigure}
\usepackage{comment}
\usepackage{float}

\usepackage{placeins}

\usepackage{xcolor, colortbl}
\definecolor{skyblue6}{rgb}{.9, .5, .4}
\definecolor{firebrick}{rgb}{.95, .9, .9}
\definecolor{lightcopper}{rgb}{.93, .76, .58}
\definecolor{blueice}{rgb}{.85, .96, .94}

\usepackage{tikz-cd}
\newtheorem {teo} {Theorem} [section]
\newtheorem{lem}[teo]{Lemma}
\newtheorem {prop} [teo]{Proposition}
\newtheorem {cor} [teo] {Corollary}

\newtheorem{rem}[teo]{Remark}
\newcommand{\eps}{\epsilon}

\begin{document}

\title[Arnold diffusion in the EH4BP]{Arnold diffusion in the elliptic Hill  four-body problem: geometric method  and numerical verification}
\author{Jaime Burgos--Garc\'\i a}
\address{Facultad de Ciencias F\'isico Matematicas, Universidad Aut\'onoma de Coahuila,
Unidad Campo Redondo, Edificio A, 25020 Saltillo, Coahuila, M\'exico}
\email{jburgos@uadec.edu.mx}
\author{Marian Gidea$^1$}\footnote{Corresponding author}
\address{Yeshiva University, Department of Mathematical Sciences, New York, NY 10016, USA }
\email{Marian.Gidea@yu.edu}

\author{Claudio Sierpe}
\address{Department of Mathematics, Sciences Faculty, University of B\'io-B\'io, Casilla 5-C, Concepci\'on, VIII-Regi\'on, Chile}
\curraddr{ Departamento de Ciencias B\'asicas, Facultad de Ciencias, Universidad Santo Tom\'as, Avenida H\'eroes de la Concepci\'on 2885, I-Regi\'on, Chile} \email{csierpe@santotomas.cl}

\thanks{The work of  M.G. was partially supported by NSF grant  DMS-2307718.
The work  of C.S. was supported by ANID-Chile through Beca de Doctorado Nacional, folio 21191713.
Part of the work   was carried out while C.S. was visiting Yeshiva University during 2022-2023.
}

\subjclass{Primary: 70F10,  70H08, 37J40; Secondary: 70F15, 37E40.}

\keywords{Elliptic four-body problem, Hill's problem, normally hyperbolic invariant manifolds, Arnold diffusion}

\begin{abstract}
We present a mechanism for  Arnold diffusion in energy in a model of the elliptic Hill four-body problem.
Our model is expressed as a small perturbation of the  circular  Hill four-body problem,
with the small parameter being the eccentricity of the orbits of the primaries.
The mechanism relies on the existence of two normally hyperbolic invariant manifolds (NHIMs), and on the corresponding homoclinic and heteroclinic connections.
The dynamics along homoclinic/heteroclinic orbits is encoded via scattering maps, which we compute numerically.
We provide two arguments for diffusion. In the first argument, we successively apply a single scattering map and use Birkhoff's Ergodic Theorem to obtain pseudo-orbits that, on average, gain energy.
In the second argument, we use two scattering maps and, at each step, select one that increases energy.
Either argument yields pseudo-orbits of scattering maps along which the energy grows by an amount independent of the small parameter. A shadowing lemma concludes  the existence of diffusing orbits.
\end{abstract}

\maketitle

\section{Introduction}

\subsection{Overview of results} We consider the planar elliptic restricted four-body problem (ER4BP), describing the dynamics of a massless body under the gravitational influence of three massive bodies (primaries) of masses $m_1>m_2>m_3$ forming an equilateral central configuration, where each  primary moves on  an elliptic orbit  about the  common  center of mass. This  represents a homographic solution of the general three-body problem.
We derive the elliptic  Hill four-body problem (EH4BP), which is an  approximation of the ER4BP  describing the  dynamics of the infinitesimal body in a
neighborhood of the smaller body $m_3$, in the limit  case when $m_3\to 0$ and $m_1,m_2$ are sent to infinity.
The EH4BP can be written as a perturbation of the circular  Hill four-body problem (CH4BP), with the eccentricity $\eps$ of the elliptic orbits being  the small parameter.
While for some orbits the effect of the
perturbation may average out, for other orbits it may accumulate in the long run.
We show that the EH4BP  exhibits Arnold diffusion, in the sense that there exist orbits of the infinitesimal body that undergo significant changes over time.
In particular, the energy along these orbits increases by $O(1)$ with respect to the perturbation parameter $\eps$.
We note that the CH4BP is given by a two-degree of freedom Hamiltonian, and the EH4BP by a two-and-a-half degree of freedom Hamiltonian (as a time-periodic perturbation of the  CH4BP), which  is the minimum possible number of degrees of freedom for which Arnold diffusion can occur.

A motivation for this work is the system consisting of  Sun, Jupiter, the Trojan asteroid (624) Hektor, and its  moonlet Skamandrios,
which can be modeled by the EH4BP.
Hektor  is located close to the  Lagrangian point  $L_4$ of the Sun-Jupiter system.
The masses of Sun, Jupiter and Hektor are $m_1= 1.989\times10^{30}$ kg, $m_2=1.898\times10^{27}$ kg, and $m_3=7.91\times10^{18}$ kg, respectively.
Some models based on astronomical observations suggest that a small perturbation could potentially lead to the ejection of  Skamandrios, or to  its collision  with the asteroid \cite{Marchis}.
Our work indicates  the theoretical possibility that Arnold diffusion may lead to such scenarios.

\subsection{Methodology} Our arguments for Arnold diffusion rely on geometric methods and numerical verifications.

The geometric method is concerned with identifying the objects that organize the dynamics.
For the unperturbed problem, corresponding to the CH4BP, we find $4$ equilibrium points: $L_1$, $L_2$, of center-saddle linear stability  type, and $L_3$, $L_4$, of center-center type. We focus on the dynamics near $L_1$, $L_2$. For a range of energies slightly above that of $L_1$, $L_2$, we find families of Lyapunov periodic orbits around these points.
Each family forms a normally hyperbolic invariant manifold (NHIM), which is an annulus that can be parameterized via symplectic action-angle coordinates.
The NHIMs posses stable and unstable manifolds.
For these energies, the motion of the massless body is confined to an inner region around $m_3$ connected to an outer region opening towards $m_1$ and $m_2$. The two regions are separated by two `bottlenecks' where the Lyapunov orbits near $L_1$, $L_2$ lie.
There exist transverse homoclinic and heteroclinic connections between the two NHIMs within both the inner region and the outer region.
Each homoclinic orbit is asymptotic in both forward and backwards time to the same Lyapunov orbit, and each  heteroclinic orbit is asymptotic in  forward and backwards time to two Lyapunov orbits -- one near $L_1$ and the other near $L_2$ -- of the same energy. We can encode the asymptotic behavior of homoclinic/heteroclinic orbits via the scattering map -- a geometrically defined map acting on the NHIMs. It turns out that, when expressed in action-angle coordinates, the unperturbed scattering map preserves the action and shifts the angle coordinate.

We then consider the effect of the perturbations on the NHIMs and on the transverse homoclinic/heteroclinic orbits. For small perturbations, the NHIMs survive while being slightly deformed, due to the standard normally hyperbolicity theory. The  transverse homoclinic/heteroclinic orbits also survive, due to the robustness of transversality. However the perturbed scattering map can increase or decrease the action coordinate. Using deformation theory, one can  compute the perturbed scattering map as an expansion in terms of the perturbation parameter: the zero-th order term in the expansion is the  unperturbed scattering map, and the first order term is given by a Hamiltonian function in the action-angle variables. This Hamiltonian can be explicitly computed via Melnikov theory in terms of improper, convergent integrals  of the perturbation along homoclinic/heteroclinic orbits of the unperturbed system. Moreover, these integrals converge exponentially fast, which allows for their efficient numerical computation.
Through the computation of the perturbed scattering map, we can identify regions where a scattering map increases the action coordinate by $O(\eps)$, as well as regions where a scattering map decreases the action coordinate by $O(\eps)$.
 In our system, there are several  scattering maps. We describe two mechanisms of diffusion. In the first mechanism, we use a single scattering map  (which is globally defined except for a measure-zero set) and apply Birkhoff's Ergodic Theorem to show the
 existence of pseudo-orbits generated by this map  along which the time average of the  actions grows by~$O(\eps)$. Hence, over a time $O(1/\eps)$ the net growth
 of the action is $O(1)$. In the second mechanism, we use two scattering maps and show that, at each step, we can always choose one of them that increases the action by ~$O(\eps)$, leading to pseudo-orbits that \emph{increase} the action by $O(1)$ in $O(1/\eps)$ steps.
Finally,  a dichotomy argument and a shadowing lemma  yield   true orbits that \emph{change}  the action by $O(1)$.
With additional information on the dynamics, we can obtain true orbits along which the action \emph{grows} by $O(1)$.

When the action increases, the corresponding Lyapunov periodic orbits get bigger, and their stable/unstable manifolds in the inner region pass closer to $m_3$. At the same time,  the `bottlenecks' around $L_1$, $L_2$ get wider, and orbits can escape more easily from the  inner region  to the  outer region.
In the context of the Sun-Jupiter-Hektor system, it is possible that, in the long run,   Arnold diffusion can push the moonlet Skamandrios towards collision with Hektor, or, on the contrary, to escape from Hektor's capture.

While our numerical calculations are not validated via  computer assisted proofs, they are sufficiently accurate by numerical analysis standards. A rigorous numerical verification  can   be done using  the CAPD library \cite{kapela2021capd}.

The novel aspects of this work can be summarized as follows:
\begin{itemize}
\item We derive the Hamiltonian for the EH4BP and write it is as a perturbation of the CH4BP.
\item We show the existence of Arnold diffusion for the EH4BP, in a model with realistic parameters for the masses of the bodies and for the energy levels.
\item We describe a new mechanism of diffusion based on iterating a single scattering map,  which, on average, leads to energy growth along pseudo-orbits.
\item We compare the growth of energy via homoclinic orbits with that via heteroclinic orbits.
\item We compute the level sets of the Hamiltonian functions that generate the scattering maps, which provide a complete chart of the scattering map dynamics.
\end{itemize}

\subsection{Related works}
The Arnold Diffusion problem \cite{Arnold64} asserts that integrable Hamiltonian systems subject to
small, generic perturbations, have diffusing
orbits along which the action coordinate changes by an
amount independent of the size of the perturbation.
While Arnold illustrated this phenomenon on a simple model consisting of a rotator and a pendulum plus a small time-dependent perturbation of a  special type,  he conjectured \emph{``I~believe that the mechanism of transition
chain that guarantees that nonstability in our example is
applicable to the general case (for example, to the problem of
three bodies)”}.

Despite its long history, there are few results on proving Arnold diffusion in celestial mechanics,  particularly with realistic parameters (e.g., for the masses of the bodies, or for the energy levels). We mention some works that prove Arnold diffusion in realistic models of the  restricted three-body problem by combining analytical methods with numerical computations.  The paper \cite{fejoz2016kirkwood} shows that  diffusion along mean motion resonances  in  the Sun-Jupiter system can be used to explain
the Kirkwood gaps in the main asteroid belt.  The paper \cite{CapinskiGL17} shows the existence of orbits that exhibit  Arnold diffusion in energy for   the Sun-Jupiter system at energy levels close to that of the comet Oterma. A computer assisted proof for Arnold diffusion in  energy for the Neptune-Triton system is given in \cite{capinski2023arnold}; moreover, this paper shows that the energy of diffusing orbits evolves according to a Brownian motion with drift.
The paper \cite{CG2025} extends this result to the full three-body problem.
Numerical evidence for Arnold diffusion in the dynamics of Jupiter’s Trojan asteroids is provided in \cite{robutel2006resonant}, by the means of  the Frequency
Map Analysis \cite{LASKAR1990}.

We mention some works that provide purely analytical proofs of Arnold diffusion in celestial mechanics models. The paper \cite{DelshamsKRS19}
considers the planar elliptic restricted three-body problem, for mass ratio and eccentricities of the primaries sufficiently
small, and constructs orbits with large drift in angular momentum. This result is extended for arbitrary masses in  \cite{guardia2023degenerate}.
The paper \cite{clarke2022inner} proves Arnold’s conjecture in a planetary spatial
4-body problem as well as in the corresponding hierarchical problem (where the bodies are increasingly
separated); in particular, they show that one of  the planet may change inclination in a chaotic fashion, and, moreover, it can flip from a prograde, nearly horizontal orbit to a retrograde one.

In the present paper, we use analytical methods and numerical verifications to show Arnold diffusion in  an elliptic Hill  four-body problem modeling the dynamics of the moonlet Skamandrios relative to the  Sun-Jupiter-Hektor system, for realistic values of the masses and energies.
In this model, the diffusion mechanism can affect the orbit of the moonlet dramatically, so that it can pass closer and closer to the asteroid, or it can escape from the asteroid's capture.
A model for this system  based on the circular  Hill  four-body problem is  derived in \cite{Burgos_Gidea_15}, and further studied in \cite{burgos2016families,burgos2022spatial}. The paper \cite{burgos2020hill} develops a model that takes into account the highly oblate shape of Hektor, and derives a circular  Hill  four-body problem in this setting.  In a follow-up paper \cite{belbruno2023regularization}, collision orbits are studied via McGehee regularization.
The geometric approach  that we use in the present paper has the advantage that it can be adapted to incorporate  additional realistic effects, such as the oblateness of the asteroid. We can also include the effect  of inclination of the orbits of the asteroid and of the moonlet, as in \cite{delshams2016arnold}, or  the effect of dissipation, as in \cite{akingbade2023arnold}.

The methodology to prove diffusion in our model is similar to that in \cite{CapinskiGL17}. The main tool that we use is the scattering map.
This was used in \cite{DelshamsLS06a} to prove Arnold diffusion in a priori unstable Hamiltonian systems; see also the related papers \cite{DelshamsLS00, DelshamsLS06b}.
The geometric properties of the scattering map have been studied in \cite{DelshamsLS08a}.  A general mechanism of diffusion that relies mostly on the scattering map, and uses only the Poincar\'e recurrence  of the  inner dynamics (which is automatically satisfied in Hamiltonian systems over regions of bounded measure) was developed in \cite{GideaLlaveSeara20-CPAM} (see also \cite{GideaLlaveSeara20-DCDS}).  This mechanism applies only in the case when the unperturbed scattering map is the identity.
Unlike in \cite{GideaLlaveSeara20-CPAM}, the unperturbed scattering map in the EH4BP is a shift in the angle variable.
In the present paper, we introduce a new mechanism of diffusion based on iterating a single scattering map and using Birkhoff's Ergodic Theorem to show that there exist pseudo-orbits whose energy increases on average.
The idea of estimating the average growth of energy along orbits also appears in \cite{GideaL17}.
Alternatively, we can use the mechanism from
\cite{CapinskiGL17}, where for each point in the NHIM we can always select a scattering map that increases the energy.

An important quality of the scattering map that we explore in this paper is that it can be computed explicitly, either analytically, via perturbation theory, or numerically \cite{DelshamsLS06a,delshams2016arnold}. In the perturbative case, the perturbed scattering map can be expanded in powers of the perturbation parameter, such  that the
zero-th order term in the expansion is the unperturbed scattering map, and the  first order term is given by a Hamiltonian function on the NHIM. These level sets can be used to detect the fastest pathways for diffusing orbits \cite{Delshams_Schaefer_2018,delshams2017arnold,DelshamsGR2025}. In the present work, we perform a direct  computation of the level sets of the Hamiltonians generating the scattering maps associated to homoclinic as well as to heteroclinic connections, and discuss their role in diffusion.

\subsection{Structure of the paper}
In Section \ref{sec:ER4BP} we describe the ER4BP.
In Section \ref{sec:EH4BP} we derive the Hill approximation EH4BP for the ER4BP, and we express  the EH4BP as a perturbation of the CH4BP.
In Section \ref{sec:equilibria} we describe the geometric objects that organize the dynamics in the CH4BP:  equilibrium points,  Lyapunov periodic orbits, and their stable and unstable manifolds. We also provide
numerical evidence for the existence of homoclinic and heteroclinic connections.
In Section \ref{sec:diffusion_EH4BP}, we present two mechanisms  to show the existence of Arnold diffusion.
In Section \ref{sec:numerical_verification}, we implement these mechanisms numerically  in the case of the EH4BP.

\section{The elliptic restricted four-body problem}
\label{sec:ER4BP}
Consider a massless particle that moves on the same plane of three massive bodies, called primaries, without affecting their motion. The equations for the massless particle in a $(X,Y)$ inertial coordinate system with origin at the center of mass of the system  are given by

\begin{equation}\begin{split}\label{eqn:fixedequations}
\frac{d^2X}{dt^2}&=-G\sum_{i=1}^{3}\frac{m_{i}(X-X_{i})}{\rho_{i}^3},\quad
\frac{d^2Y}{dt^2}=-G\sum_{i=1}^{3}\frac{m_{i}(Y-Y_{i})}{\rho_{i}^3},
\end{split}
\end{equation}

where $(X_{i},Y_{i})$ are the positions of the three primaries,
$\rho_{i}$ is the distance between the massless particle and each primary for $i=1,2,3$, and $G$ denotes the gravitational constant. Using the complex coordinate $Z=X+\textit{i}Y$ we can write the equations \eqref{eqn:fixedequations} in the following form
\begin{equation}\label{eqn:complexfixedequation}
\frac{d^2Z}{dt^2}=-G\sum_{i=1}^{3}\frac{m_{i}(Z-Z_{i})}{\rho_{i}^3},
\end{equation}
with $\rho_{i}=\vert Z-Z_{i}\vert$. The positions of the primaries can be considered as solutions of the three-body problem such as the Euler solution \cite{maranhao1999ejection}, the eight-shaped solution \cite{lara2022restricted}, or the Lagrange solution to produce different kinds of  \textit{restricted four-body problems}. Some authors have considered elliptical motions for the primaries which are not actual solutions of the three-body to generate other kinds of  \textit{elliptic restricted four-body problems} \cite{liu2018hill,carletta2022characterization}.

The restricted four-body problem, with the three primaries in a Lagrange configuration and moving in circular orbits, has been recently studied by several authors. See, for instance, \cite{burgos2013periodic,burgos2013blue,burgos2019spatial} and references therein. However, there are only a few studies of the elliptic case; see  \cite{chakraborty2019new,ansari2020generalized}, where the equations of motion were derived for special values of the masses. Since the Lagrange solution exists for all  physical values of the masses,  and we are interested in the case when  one primary is much smaller than the other two, we provide a comprehensive discussion of the problem and derive the equations of motion for the elliptic restricted four-body problem (ER4BP).

\begin{figure}
	\centering
	\includegraphics[width=0.75\textwidth]{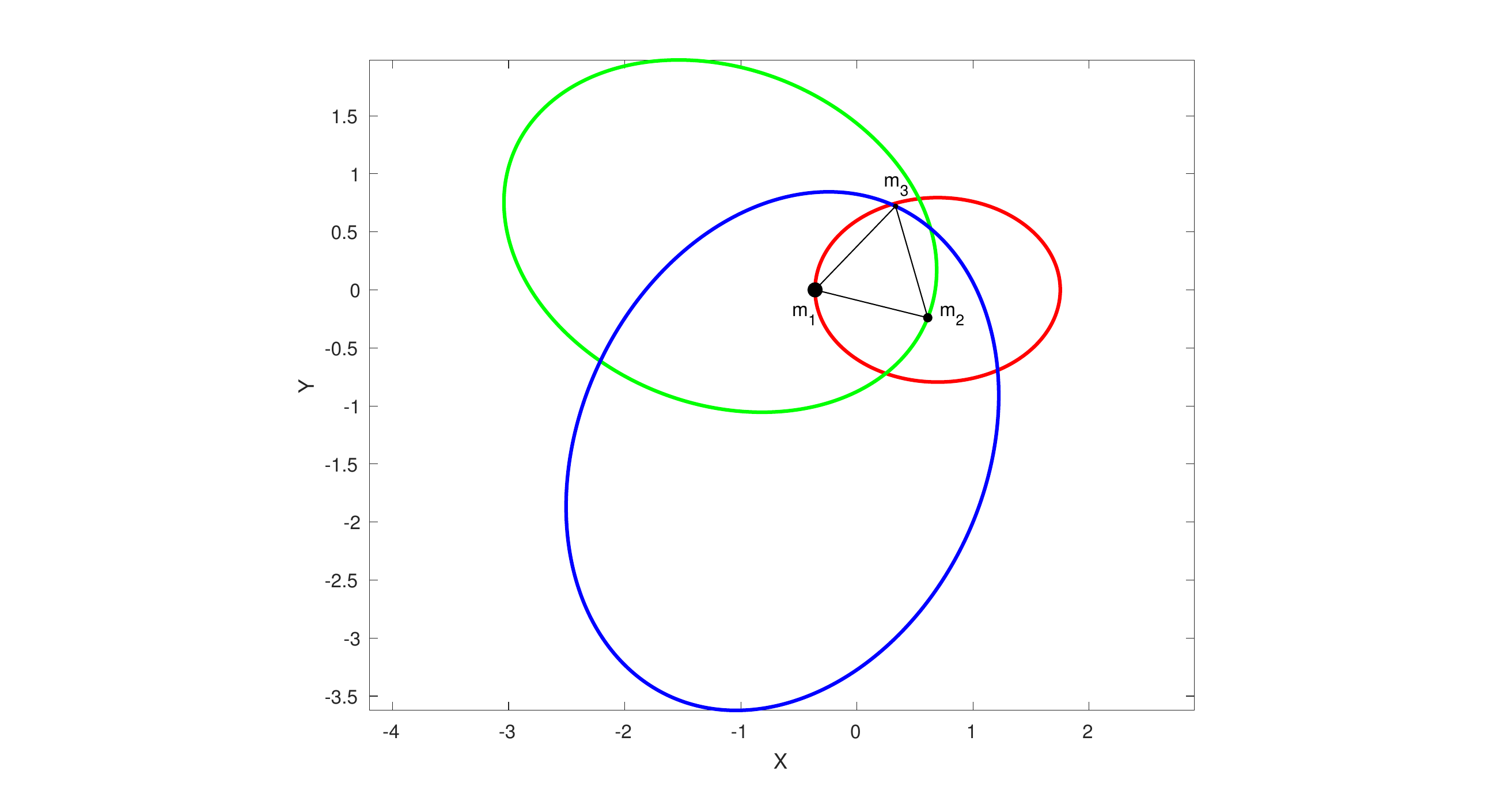}
	\caption{A Lagrange central configuration. 
}\label{fig:lagrangiansolution}
\end{figure}

We begin by expressing the solutions of the planar three-body problem  (in complex notation) in the form $q_{i}(t)=\phi(t)a_i$, where each $a_i$ is a  complex number and $\phi(t)$ is a time-dependent complex-valued function. Substituting this expression into the equations of the three-body problem gives rise to the following equations \cite{meyer2009introduction}:
\begin{align}
\frac{d^2\phi}{dt^2} &=-\frac{\lambda\phi}{\vert\phi\vert^3},\label{eqn:cc_kepler_part}\\
\frac{\partial U}{\partial q}(a)+\lambda\frac{\partial I}{\partial q}(a) &= 0,\label{eqn:cc_equation}
\end{align}
where $q=(q_1, q_2, q_3)$, $a=(a_1, a_2, a_3)$, $U$ is the potential of the three-body problem and $I=\frac{1}{2}\left(m_1\vert q_1\vert^2+m_2\vert q_2\vert^2+m_3\vert q_3^2\vert\right)$ is the moment of inertia of the system. The solutions of \eqref{eqn:cc_equation} give a central configurations of the three-body problem, and the equation \eqref{eqn:cc_kepler_part}  is a two-dimensional Kepler problem.  The solutions of \eqref{eqn:cc_equation}  are the celebrated Euler (collinear) and Lagrange (equilateral) configurations. See Fig.~\ref{fig:lagrangiansolution}. We note that the Lagrange  configuration has also been studied in the case where the primaries are assumed to be oblate bodies \cite{burgos2020hill}. The main conclusion is that the equilateral shape of the configuration deforms into either an  isosceles or  a scalene triangle, depending on the oblateness of the bodies.

For the equilateral triangle configuration, the distances between the three primaries $\rho_{ij}=\| q_{i}-q_{j}\|$ satisfy
\begin{equation*}
\rho_{12}=\rho_{23}=\rho_{31}=\left(\frac{GM}{\lambda}\right)^{-1/3}
\end{equation*}
where $\lambda$ is a scale factor that fixes the size of the system and $M=m_{1}+m_{2}+m_{3}$.  For the elliptic solutions of (\ref{eqn:cc_kepler_part}) the function $\phi(t)$ can be written in polar coordinates as $\phi(t)=re^{\textit{i}f}$, where $r$ and $f$ satisfy the two dimensional Kepler problem
\begin{eqnarray}
\frac{dr}{dt}-r\left(\frac{df}{dt}\right)^2&=&-\frac{\lambda}{r^2},\label{eqn:polarkepler}\\
r^2\frac{df}{dt}&=&c,\label{eqn:angularmomentum}
\end{eqnarray}
and they are related by
\begin{equation}\label{eqn:keplerellipse}
r=\frac{c^2/\lambda}{1+\eps\cos(f)}.
\end{equation}
The constant $c$ is the magnitude of the angular momentum, $\eps$ denotes the eccentricity, and $f$ is the true anomaly.

We will consider  $f$ as the new  time variable (which is $2\pi$-periodic), and  derive the equations of motion with respect to the new time~$f$.

To fix the positions of the primaries, we introduce the so-called pulsating coordinates
\begin{equation}
Z=re^{\textit{i}f}z \label{eqn:pulsatingchange}.
\end{equation}
We apply the chain rule to obtain the expressions
\begin{equation*}
\frac{dz}{dt}=\frac{dz}{df}\frac{df}{dt},\quad \frac{d^2z}{dt^2}=\frac{dz}{df}\frac{d^2f}{dt^2}+\frac{d^2z}{df^2}\left(\frac{df}{dt}\right)^2,
\end{equation*}
which are necessary to compute the second derivative in the expression (\ref{eqn:pulsatingchange}). A straightforward computation shows that the second derivative of $Z$ is given by
\begin{equation*}
\begin{split}
\frac{d^2Z}{dt^2}e^{-\textit{i}f}&=r\left(\frac{df}{dt}\right)^2\left(\frac{d^2z}{df^2}+2\textit{i}\frac{dz}{df}\right)+z\left(\frac{dr}{dt}-r\left(\frac{df}{dt}\right)^2\right)\\
&+\left(\frac{dz}{df}+\textit{i}z\right)\left(r^2\frac{d^2f}{dt^2}+2\frac{dr}{dt}\frac{df}{dt}\right).
\end{split}
\end{equation*}

The second term on the right side of the above equation corresponds to the equation (\ref{eqn:polarkepler}), and the third term are equal to zero because the term
$$r^2\frac{d^2f}{dt^2}+2\frac{dr}{dt}\frac{df}{dt}$$
is the derivative of the expression (\ref{eqn:angularmomentum}). Therefore, the left side of the equation (\ref{eqn:complexfixedequation}) becomes
\begin{equation*}
\frac{d^2Z}{dt^2}=e^{\textit{i}f}\left(r\left(\frac{df}{dt}\right)^2\left(\frac{d^2z}{df^2}+2\textit{i}\frac{dz}{df}\right)-z\frac{\lambda}{r^2}\right).
\end{equation*}

On the other hand, using the expression (\ref{eqn:pulsatingchange}), the right-side of (\ref{eqn:complexfixedequation}) becomes
\begin{equation*}
-G\sum_{i=1}^{3}\frac{m_{i}(Z-Z_{i})}{\rho_{i}^3}=-\frac{Ge^{\textit{i}f}}{r^2}\sum_{i=1}^{3}\frac{m_{i}(z-z_{i})}{\vert z-z_{i}\vert^3}.
\end{equation*}
If we define $\mu_{i}=m_{i}/M$, where $M$ is the total mass of the system, the equation (\ref{eqn:complexfixedequation}) in pulsating coordinates is given by the expression:
\begin{equation}\label{eqn:pulsating_preliminar_eqn}
r\left(\frac{df}{dt}\right)^2\left(\frac{d^2z}{df^2}+2\textit{i}\frac{dz}{df}\right)=z\frac{\lambda}{r^2}-\frac{GM}{r^2}\sum_{i=1}^{3}\frac{\mu_{i}(z-z_{i})}{\vert z-z_{i}\vert^3}.
\end{equation}

From equations (\ref{eqn:keplerellipse}) and (\ref{eqn:angularmomentum}) we obtain
$$\frac{1}{r^2}=\frac{r\left(\frac{df}{dt}\right)^2}{\lambda(1+\eps\cos(f))},$$ so the equation (\ref{eqn:pulsating_preliminar_eqn}) becomes
\begin{equation*}
r\left(\frac{df}{dt}\right)^2\left(\frac{d^2z}{df^2}+2\textit{i}\frac{dz}{df}\right)=\frac{r\left(\frac{df}{dt}\right)^2}{1+\eps\cos(f)}\left(z-\frac{GM}{\lambda}\sum_{i=1}^{3}\frac{\mu_{i}(z-z_{i})}{\vert z-z_{i}\vert^3}\right).
\end{equation*}

We  recall that the distances between the primaries are given by $\rho_{ij}=(\frac{GM}{\lambda})^{-1/3}$ so, in order to fix the length of the triangle equal to $1$, we choose $\lambda=GM$. Dividing the above equation by $r\left(\frac{df}{dt}\right)^2$ we obtain
\begin{equation}\label{eqn:pulsating_eqns}
\frac{d^2z}{df^2}+2\textit{i}\frac{dz}{df}=\frac{1}{1+\eps\cos(f)}\left(z-\sum_{i=1}^{3}\frac{\mu_{i}(z-z_{i})}{\vert z-z_{i}\vert^3}\right).
\end{equation}

In the new coordinates, the positions $z_i$ are the solutions of equation (\ref{eqn:cc_equation}) and consequently,  the positions are constant and the coordinates are given by the same expressions used in the circular version of the ER4BP:
\begin{align*}
		x_{1}&=\frac{-\vert K\vert\sqrt{m_{2}^{2}+m_{2}m_{3}+m_{3}^{2}}}{K} & y_{1}&=0\\
	    x_{2}&=\frac{\vert K\vert[(m_{2}-m_{3})m_{3}+m_{1}(2m_{2}+m_{3})]}{2K\sqrt{m_{2}^{2}+m_{2}m_{3}+m_{3}^{2}}} &
		y_{2}&=\frac{-\sqrt{3}m_{3}}{2m_{2}^{3/2}}\sqrt{\frac{m_{2}^{3}}{m_{2}^{2}+m_{2}m_{3}+m_{3}^{2}}} \\
		x_{3}&=\frac{\vert K\vert}{2\sqrt{m_{2}^{2}+m_{2}m_{3}+m_{3}^{2}}} & y_{3}&=\frac{\sqrt{3}}{2\sqrt{m_{2}}}\sqrt{\frac{m_{2}^{3}}{m_{2}^{2}+m_{2}m_{3}+m_{3}^{2}}},
\end{align*}
where $K=m_{2}(m_{3}-m_{2})+m_{1}(m_{2}+2m_{3})$, and we have reconsidered the mass ratios $\mu_i$ as $m_i$ in such a way that they satisfy $m_{1}+m_{2}+m_{3}=1$.

If we denote by $\dot{z}=\frac{dz}{df}$ and separate the real and imaginary parts of $z=x+\textit{i}y$ in (\ref{eqn:pulsating_eqns}), we obtain the equations in Cartesian coordinates
\begin{equation}\begin{split}\label{eqn:pulsating_cartesian_eqns}
\ddot{x}-2\dot{y}=\Omega_{x},\\
\ddot{y}+2\dot{x}=\Omega_{y},
\end{split}
\end{equation}
where
$$\Omega(x,y;m_{1},m_{2},m_{3})=\frac{1}{1+\eps\cos(f)}\left(\frac{1}{2}(x^{2}+y^{2})+\sum_{i=1}^{3}\frac{m_{i}}{r_{i}}\right),
$$
and $r_{i}=\sqrt{(x-x_{i})^{2}+(y-y_{i})^{2}}$, for $i=1,2,3$.  Making  the transformation $\dot x=p_x+y$, $\dot y=p_y-x$,  we obtain that the equations (\ref{eqn:pulsating_cartesian_eqns}) are equivalent to the equations given by the Hamiltonian
\begin{equation}
H=\frac{1}{2}(p_{x}^{2}+p_{y}^{2})+yp_{x}-xp_{y}+\frac{1}{2}(x^{2}+y^{2})-\frac{1}{1+\eps\cos(f)}\left(\frac{1}{2}(x^{2}+y^{2})
+\sum_{i=1}^{3}\frac{m_{i}}{r_{i}}\right).
\label{pulsatinghamiltonian}
\end{equation}

When we let $\eps=0$ in \eqref{pulsatinghamiltonian}, we obtain the Hamiltonian of the circular restricted four-body problem $H_{CR4BP}$, which is a two-degree-of-freedom Hamiltonian.

As  \eqref{pulsatinghamiltonian} depends  on the `new time' $f$, it is a two-and-a-half degree of freedom  Hamiltonian. This is the minimal number of degrees of freedom for which Arnold diffusion can occur. However, we will not show Arnold diffusion for this Hamiltonian, but for its Hill approximation, obtained when $m_{3}\rightarrow 0$ and $m_1,m_2$ are sent to infinity. See Section \ref{sec:EH4BP}.

\begin{rem} Since a planar homographic solution of the general $N$-body problem is defined as $q_{i}(t)=r(t)e^{\textit{i}\theta(t)}a_i$ where the $a_i$'s describe a central configuration, and $r(t)$ and $\theta(t)$ satisfy  the equations (\ref{eqn:polarkepler}) and (\ref{eqn:angularmomentum}), the performed construction can be easily extended to generate elliptic  restricted problems for $N>3$ primaries moving in homographic solutions.
\end{rem}

\section{The limit case and equations of motion}
\label{sec:EH4BP}
In this section we  study the limit case when $m_{3}\rightarrow 0$ in the Hamiltonian of the ER4BP. For our purposes it is sufficient to consider small values of the eccentricity, however, the procedure can be easily extended for all values of eccentricity. Expanding the expression $1/(1+\eps\cos(f))$ in terms of the eccentricity in a neighborhood of $\eps=0$ we have
\begin{equation*}\label{expansioneccentricity}
\frac{1}{1+\eps\cos(f)}=1-\eps\cos(f)+\mathcal{O}(\eps^2).
\end{equation*}
Therefore the Hamiltonian (\ref{pulsatinghamiltonian}) becomes
\begin{equation}\label{hamiltoniansmalle}
H_{ER4BP}=H_{CR4BP}+\eps\cos(f)\left(\frac{1}{2}(x^{2}+y^{2})+\sum_{i=1}^{3}\frac{m_{i}}{r_{i}}\right)+\mathcal{O}(\eps^2),
\end{equation}
where $H_{CR4BP}$ is the Hamiltonian of the circular restricted four body problem.

\begin{teo} \label{main theorem}
Let $H_{CH4BP}$ be the Hamiltonian of the circular Hill four-body problem, given by
\begin{equation}\begin{split}\label{circularhillhamiltonian}
H_{CH4BP}=&\frac{1}{2}(p^{2}_{x}+p^{2}_{y})+yp_{x}-xp_{y}-U(x,y),
\end{split}\end{equation}
with \begin{equation*}\label{eqn:U} U(x,y)=-\frac{1}{8}x^2+\frac{3\sqrt{3}}{4}(1-2\mu)xy+\frac{5}{8}y^2
+\frac{1}{\sqrt{x^2+y^2}}\end{equation*}
 being the gravitational potential relative to the co-rotating frame, $m_{1}=1-\mu$ and $m_{2}=\mu$.

Then, for small values of the eccentricity $\eps$, and after the conformally symplectic scaling
\[(x,y,p_{x},p_{y})\rightarrow m_{3}^{1/3}(x,y,p_{x},p_{y}),\]
 the limiting   Hamiltonian \eqref{hamiltoniansmalle} as $m_{3}\rightarrow0$ is
\begin{equation}\label{elliptichillhamiltonian}
H_{EH4BP}=H_{CH4BP}+\eps\cos(f)\left(\frac{1}{2}(x^{2}+y^{2})+U(x,y)\right)+O(\eps^2).
\end{equation}
 The system associated to   $H_{EH4BP}$ yields an approximation of the motion of the infinitesimal mass in an $O(m_3^{1/3})$-neighborhood of $m_3$,
and will be referred to as the Elliptic Hill's Four-Body Problem (EH4BP).
\end{teo}

\begin{proof} If we ignore the higher order terms in the Hamiltonian (\ref{hamiltoniansmalle}) we can write it as
\begin{equation*}\label{perturbation}
H=H_{CR4BP}+\eps\cos(f)H_{1}.
\end{equation*}

In \cite{Burgos_Gidea_15} we performed  the symplectic  coordinate changes
\begin{equation}\label{eqn:symplectic_transformation}
x\rightarrow x+x_{3},\quad y\rightarrow y+y_{3},\quad p_{x}\rightarrow p_{x}-y_{3},\quad p_{y}\rightarrow p_{y}+x_{3},
\end{equation}
and then  the conformally symplectic scaling with multiplier $m_{3}^{-2/3}$
\begin{equation}\label{eqn:symplectic_scaling} x\rightarrow m_{3}^{1/3}x,\quad y\rightarrow m_{3}^{1/3}y,\quad p_{x}\rightarrow m_{3}^{1/3}p_{x},\quad p_{y}\rightarrow m_{3}^{1/3}p_{y},\end{equation}
which transform the Hamiltonian $H_{CR4BP}$ into $H_{CH4BP}$.

We now make the same coordinate transformation again,
using the formulas from \cite{Burgos_Gidea_15}, to obtain the effect of
 those transformations on $H_1$.

After the coordinate change \eqref{eqn:symplectic_transformation}, the term $H_1$ becomes
\[
H_{1}=\frac{1}{2}(x^{2}+y^{2})+x_{3}x+y_{3}y+\sum_{i=1}^{3}\frac{m_{i}}{\bar{r}_{i}},
\]
where $\bar{r}_{i}^{2}=(x+x_{3}-x_{i})^2+(y+y_{3}-y_{i})^2:=(x+\bar{x}_{i})^2+(y+\bar{y}_{i})^2$, for $i=1,2,3$;
 we have omitted the constant term $\frac{1}{2}(x_{3}^2+y_{3}^2)$.
We expand the terms $\frac{1}{\bar{r}_{1}}$ and $\frac{1}{\bar{r}_{2}}$ in Taylor series around the new origin of coordinates. Ignoring the constant terms, we obtain
\begin{equation*}\begin{split}\frac{1}{\bar{r}_{1}}=\sum_{k\ge1}P_{k}^{1}(x,y),\\
\frac{1}{\bar{r}_{2}}=\sum_{k\ge1}P_{k}^{2}(x,y),\end{split}
\end{equation*}
where $P_{k}^{j}(x,y)$ is a homogenous polynomial of degree $k$ for $j=1,2$.

Then we perform  the conformally symplectic scaling \eqref{eqn:symplectic_scaling} obtaining
\begin{equation*}\begin{split}\label{scaled}
H_{1}=&\frac{1}{2}(x^{2}+y^{2})+\frac{1}{\sqrt{x^{2}+y^{2}}}+
m_{3}^{-1/3}(x_{3}x+y_{3}y+m_{1}P_{1}^{1}+m_{2}P_{1}^{2})\\
&+\sum_{k\ge2}m_{3}^{\frac{k-2}{3}}m_{1}P_{k}^{1}(x,y)+\sum_{k\ge2}m_{3}^{\frac{k-2}{3}}m_{2}P_{k}^{2}(x,y),
\end{split}
\end{equation*}
where
\begin{equation*}\begin{split}P_{1}^{1} =\left(\frac{x_{1}-x_{3}}{\bar{r}_{1}^{3}}x+\frac{y_{1}-y_{3}}{\bar{r}_{1}^{3}}y\right),\quad  P_{1}^{2} =\left(\frac{x_{2}-x_{3}}{\bar{r}_{2}^{3}}x+\frac{y_{2}-y_{3}}{\bar{r}_{2}^{3}}y\right ),\end{split}\end{equation*}

The term
\[
m_{3}^{-1/3}(x_{3}x+y_{3}y+m_1P_{1}^{1}+m_2 P_{1}^{2})+\sum_{k\ge2}m_{3}^{\frac{k-2}{3}}m_{1}P_{k}^{1}(x,y)+\sum_{k\ge2}m_{3}^{\frac{k-2}{3}}m_{2}P_{k}^{2}(x,y)
\]
was already computed for the circular case $H_{CR4BP}$ in \cite{Burgos_Gidea_15}, so if we use those computations we can write
\[
H_{1}=\frac{1}{2}(x^{2}+y^{2})+\frac{1}{\sqrt{x^{2}+y^{2}}}+m_{1}P_{2}^{1}+m_{2}P_{2}^{2}+\mathcal{O}(m_{3}^{1/3}).
\]
Finally, taking the limit $m_{3}\rightarrow0$ we obtain
\[
H_{1}=\frac{1}{2}(x^{2}+y^{2})+\frac{1}{\sqrt{x^2+y^2}}-\frac{1}{8}x^2+\frac{3\sqrt{3}}{4}(1-2\mu)xy+\frac{5}{8}y^2,
\]
as desired.
\end{proof}

We can rotate the coordinate system so that the axis from the center of mass to $m_3$ becomes the $x$-axis:

\begin{cor} \label{main cor}
The Hamiltonian \eqref{circularhillhamiltonian} can be written,  after a coordinate rotation which has eigenvalues $\lambda_1=\frac{3}{2}(1-d)$ and $\lambda_2=\frac{3}{2}(1+d)$, where $d=\sqrt{1-3\mu+3\mu^2}$, in the form
\begin{equation*}\begin{split}\label{elliptichillhamiltonian_rot_}
H_{CH4BP}=&\frac{1}{2}(p^{2}_{x}+p^{2}_{y})+yp_{x}-xp_{y}-U_{rot}(x,y), \textrm{ with }
\end{split}\end{equation*}
\begin{equation}\label{eqn:U_rot} U_{rot} =-ax^2-by^2
+\frac{1}{\sqrt{x^2+y^2}}\end{equation}
where $a=  \frac{1}{2}(1-\lambda_2)$ and $b= \frac{1}{2}(1-\lambda_1)$.

The effective potential and the Hamiltonian of the $CH4BP$ are given by
\begin{equation*}
\Omega_{eff}=\frac{1}{2}( x^2+ y^2)+U_{rot}(x,y)=\frac{1}{2}(\lambda_2 x^2+\lambda_1 y^2)+\frac{1}{\sqrt{x^2+y^2}},
\end{equation*}
\begin{equation*}\begin{split}\label{elliptichillhamiltonian_eff}
H_{CH4BP}=&\frac{1}{2}(p_{x}+y)^2 +\frac{1}{2}(p_{y}-x)^2 -\Omega_{eff}(x,y).
\end{split}\end{equation*}

The Hamiltonian of the $EH4BP$ can be written as
\begin{equation}\label{elliptichillhamiltonian_rot}
H_{EH4BP}=H_{CH4BP}+\eps\cos(f)\left(\frac{1}{2}(x^{2}+y^{2})+U_{rot}(x,y)\right)+O(\eps^2).
\end{equation}
\end{cor}

\section{Equilibrium points, Hill regions, and invariant manifolds in the CH4BP}
\label{sec:equilibria}

We assume $\mu=0.00095$, representing the relative mass of Jupiter in the Sun-Jupiter system.

For the case when $\eps=0$, we obtain the  circular Hill Four-Body Problem (CH4BP) which is described by the Hamiltonian $H_{CH4BP}$ given by \eqref{elliptichillhamiltonian}. After applying a rotation of coordinates, as in Corollary \ref{main cor}, the Hamiltonian of  the CH4BP reduces  to:
\begin{equation*}\label{elliptichillhamiltonian_rot_e0}
H_{CH4BP}=\frac{1}{2}(p^{2}_{x}+p^{2}_{y})+yp_{x}-xp_{y}-U_{rot}(x,y).
\end{equation*}
The equations of motion possess the following symmetries:
\begin{itemize}
\item  $S(x,y,\dot{x},\dot{y},t)=(x,-y,-\dot{x},\dot{y},-t)$, with  respect to the $x$-axis;
\item  $S'(x,y,\dot{x},\dot{y},t)=(-x,y,\dot{x},-\dot{y},-t)$, with  respect to the $y$-axis;
\item $S''=S\circ S'(x,y,\dot{x},\dot{y},t)=(-x,-y,-\dot{x},-\dot{y},t)$, the composition of $S$ and $S'$, with respect to the origin.
\end{itemize}

\begin{figure}
  \centering
\begin{tabular}{cc}
  \includegraphics[width=0.32 \textwidth]{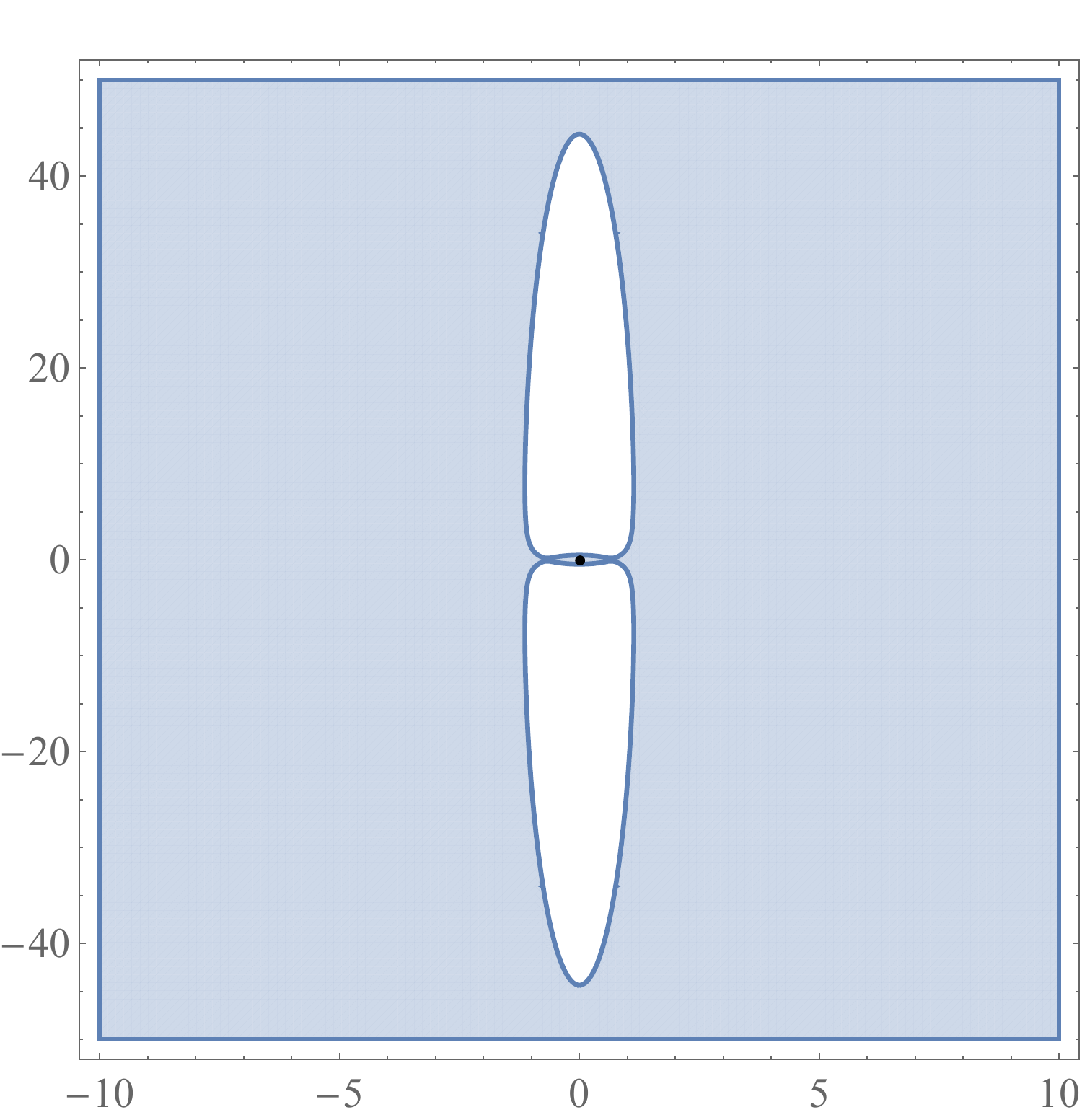}& \includegraphics[width=0.32 \textwidth]{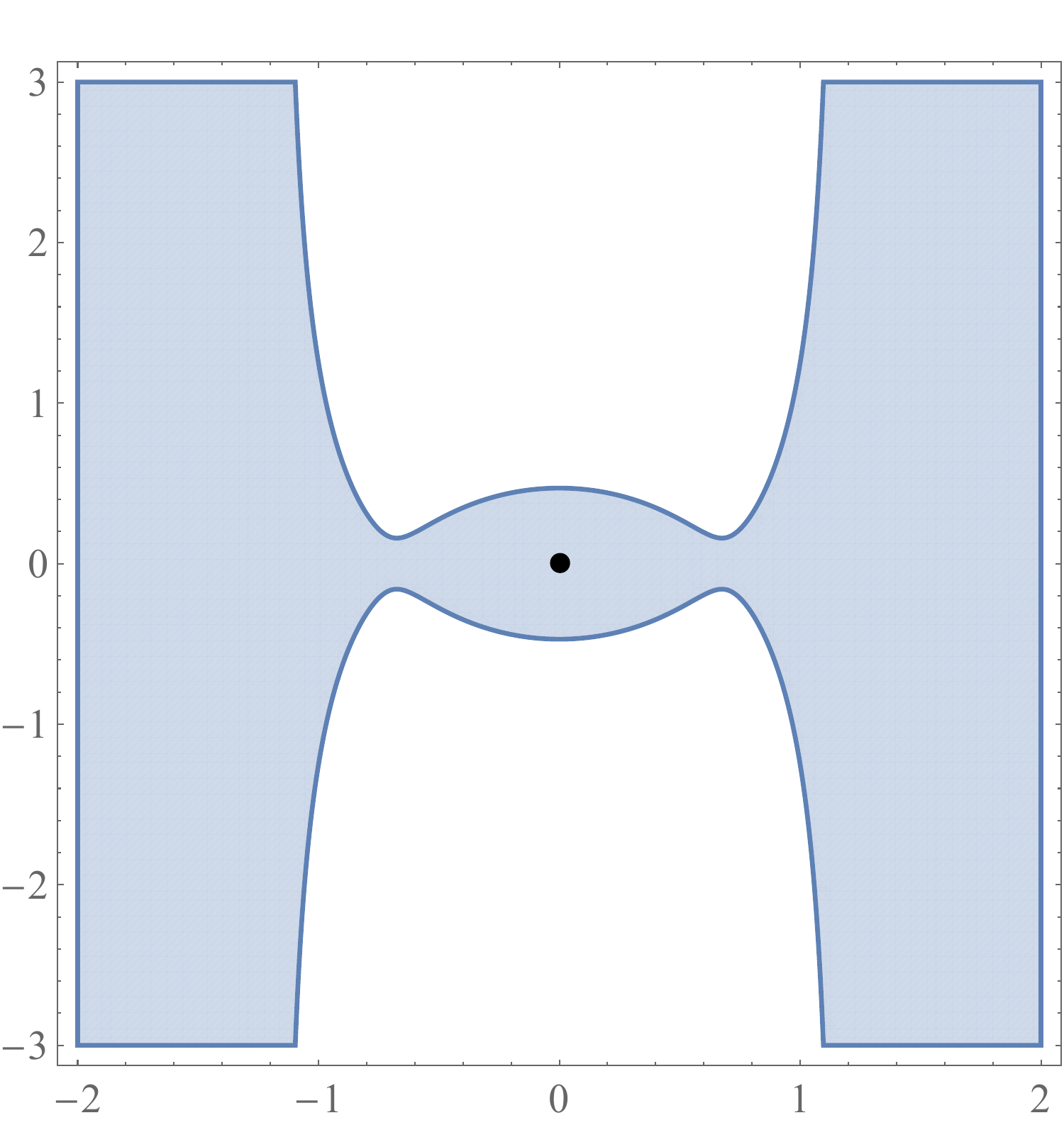}
\end{tabular}
 \caption{Hill's region for the CH4BP for $\mu= 0.00095$ and $h=-2.125$. The second figure is a magnification of the first one.}\label{limithillregions}
\end{figure}

The equations of motion for the circular case have 4 equilibrium points given by the $(x,y)$-coordinates   \[L_{1}=\left(\frac{1}{\sqrt[3]{\lambda_{2}}},0\right),\, L_{2}=\left(-\frac{1}{\sqrt[3]{\lambda_{2}}},0\right),\,
L_{3}=\left(0,\frac{1}{\sqrt[3]{\lambda_{1}}}\right),\,
L_{4}=\left(0,-\frac{1}{\sqrt[3]{\lambda_{1}}}\right).\]
The linear stability of  $L_1$ and $L_2$ is  of center-saddle  type, for all $\mu$, while that of  $L_3$ and $L_4$ is of center-center type, for $\mu$ less or equal than some critical value $\mu_{cr}$, and of complex-saddle type otherwise; see \cite{Burgos_Gidea_15}. The value of the energy $H_{CH4BP}$ at the equilibrium point $L_1$ is $h_{L_1}=-2.16286$.
In Fig. \ref{limithillregions}  we show  the Hill regions
\[ \{(x,y)\,|\,\Omega_{eff}\geq -h\}\]
for the value $h=-2.125$ of $H_{CH4BP}$,
where we  observe the \textit{inner region}  around the tertiary connected with the \textit{outer region} through two `bottlenecks', which become wider as we increase the energy level $h$.

We compute numerically the stable $W^{s}$ and unstable $W^{u}$ manifolds of the Lyapunov orbits for the saddle-center equilibrium points $L_{1}$ and $L_{2}$,
using the $(x,y,\dot x,\dot y)$ coordinates, where $\dot x =p_x+y$ and $\dot y =p_y-x$.
We also compute the intersection of these manifold with the section \[ \mathscr{S}_x=\{(x,y,\dot{x},\dot{y})\in\mathbb{R}^{4}\,\vert\, y=0\}.\]
In the $(x,y)$-configuration space, $y=0$ corresponds  to intersections with the $x$-axis.

\begin{figure}
  \centering
  \includegraphics[width=0.32 \textwidth]{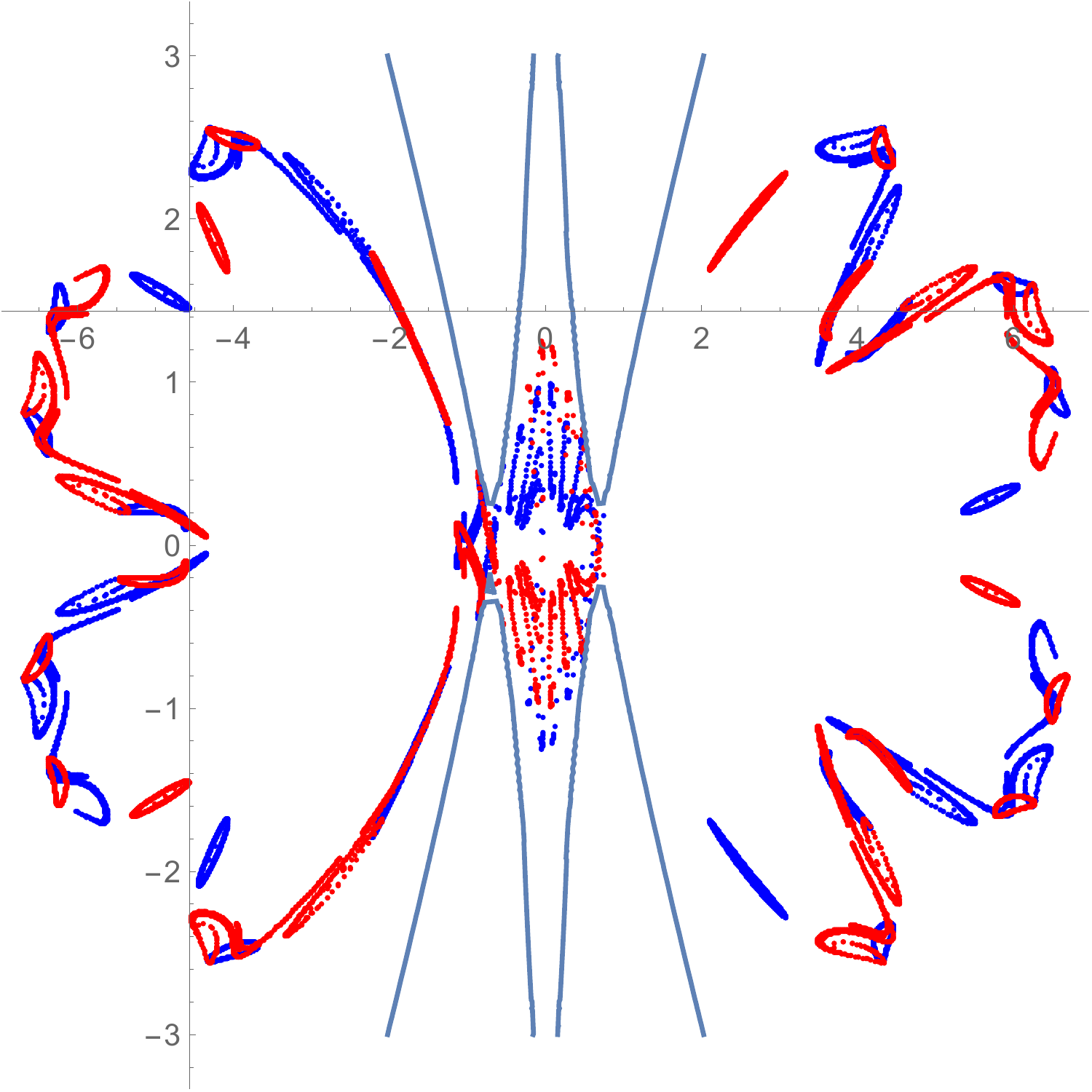}
  \includegraphics[width=0.32 \textwidth]{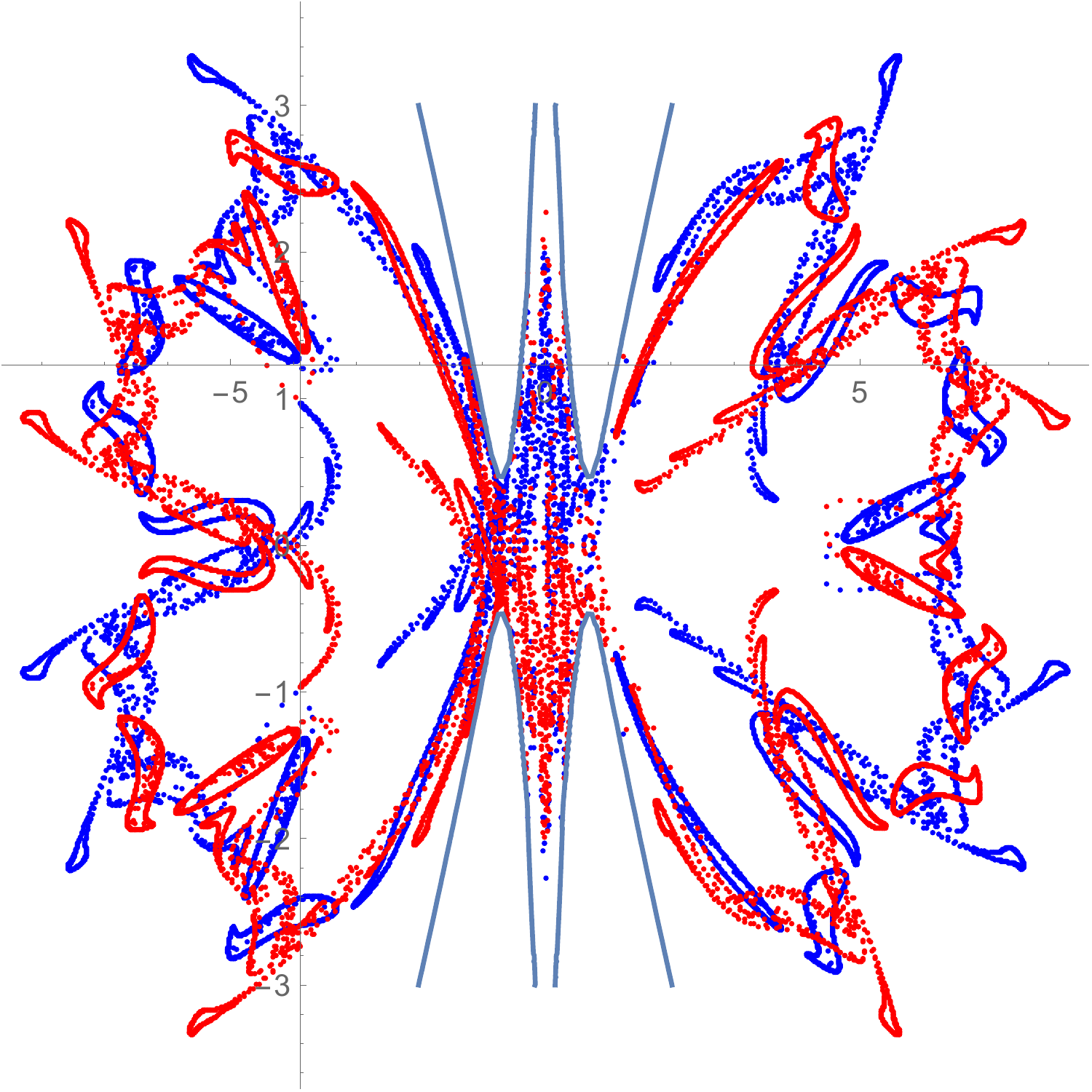}
  \includegraphics[width=0.32 \textwidth]{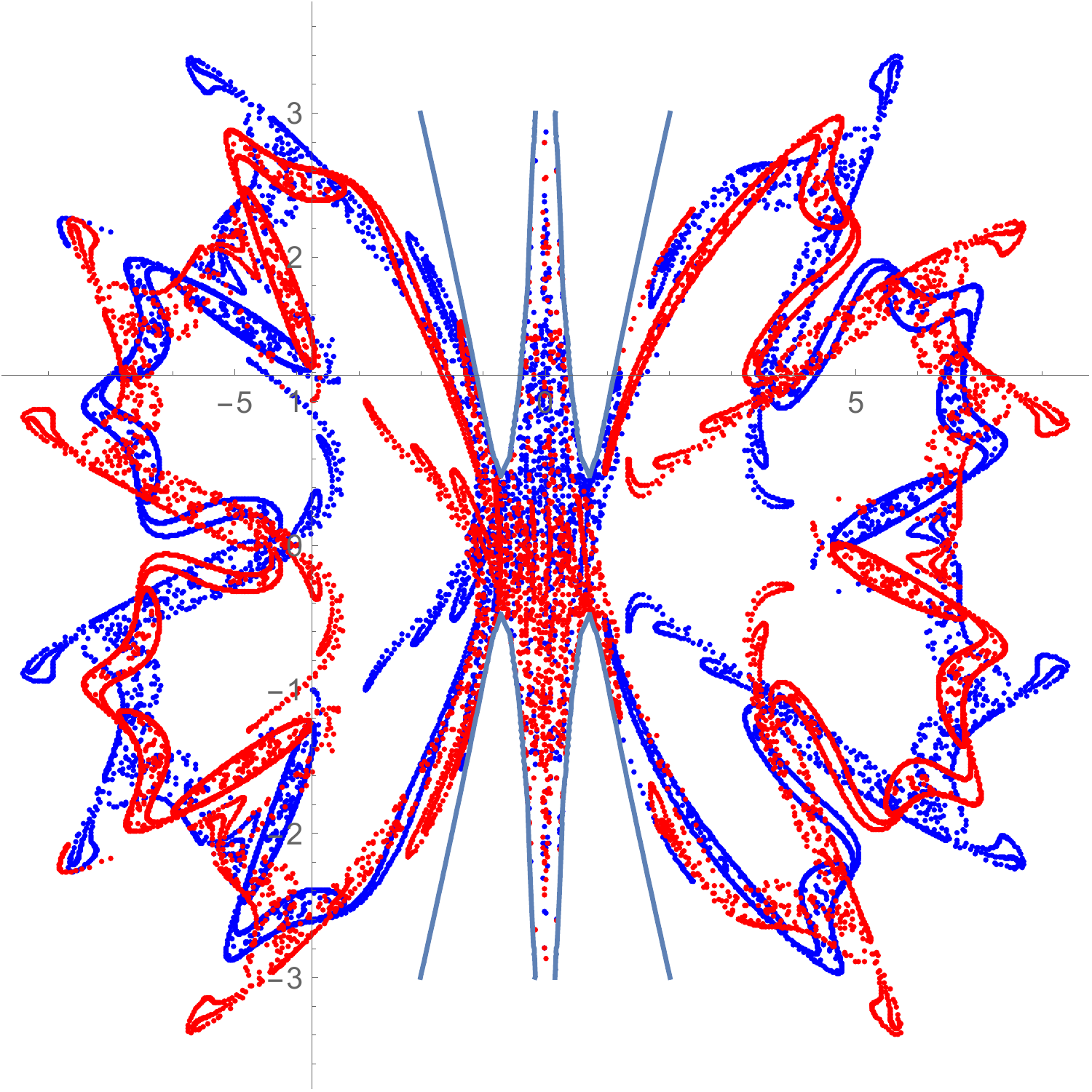}
      \caption{Stable (red) and unstable (blue) manifolds after 25 cuts with the Poincar\'e section $\mathscr{S}_x=\{y=0\}$, and tangency curve (grey) in the plane $(x,\dot{x})$, corresponding to trajectories that intersect $\mathscr{S}_x$ tangentially: $h=-2.15$ (left); $h=-2.075$ (center); $h=-2.0$ (right).}\label{cuts}
\end{figure}

\begin{figure}
  \centering
  \includegraphics[width=0.32\textwidth]{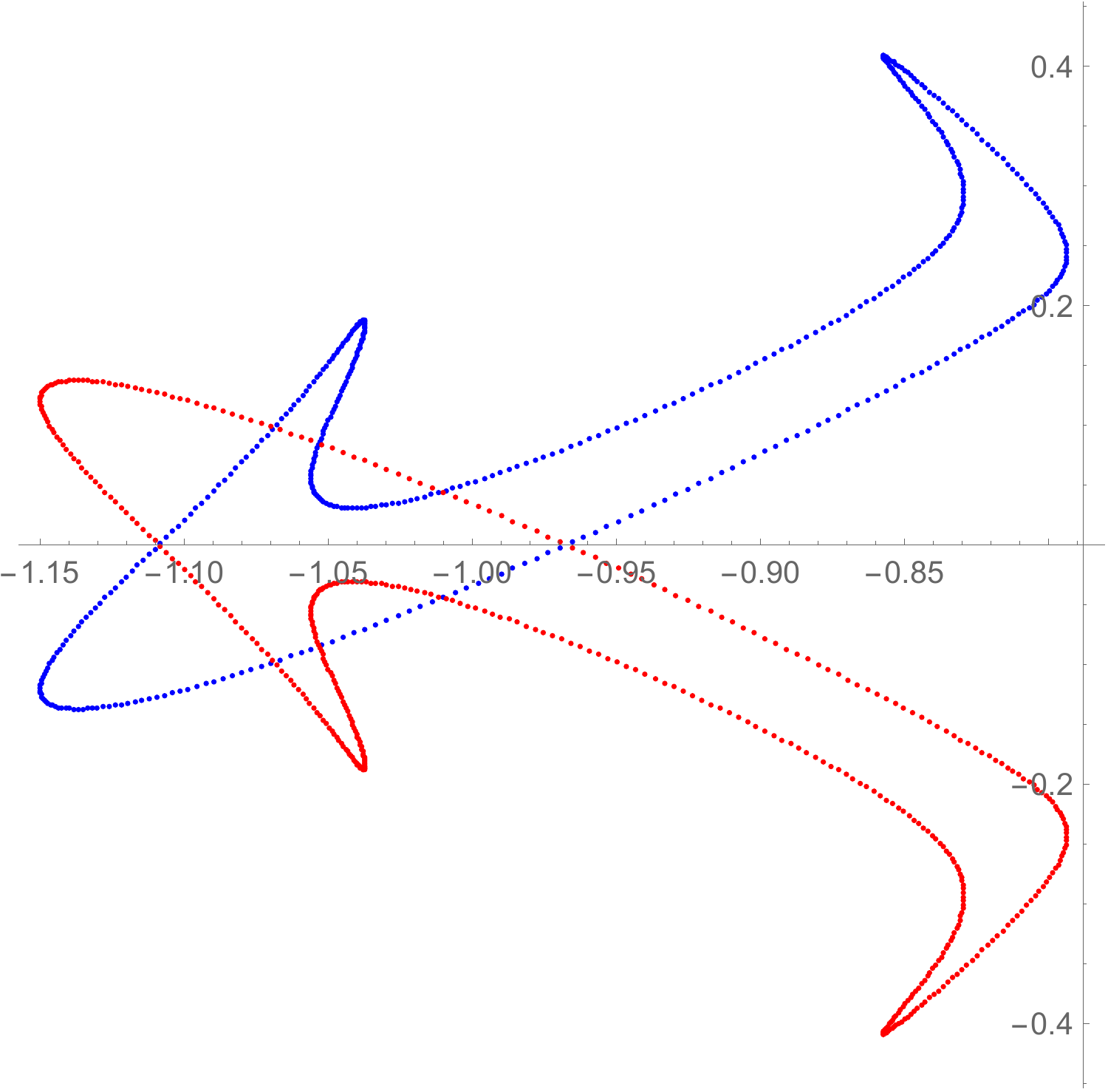}
  \includegraphics[width=0.32 \textwidth]{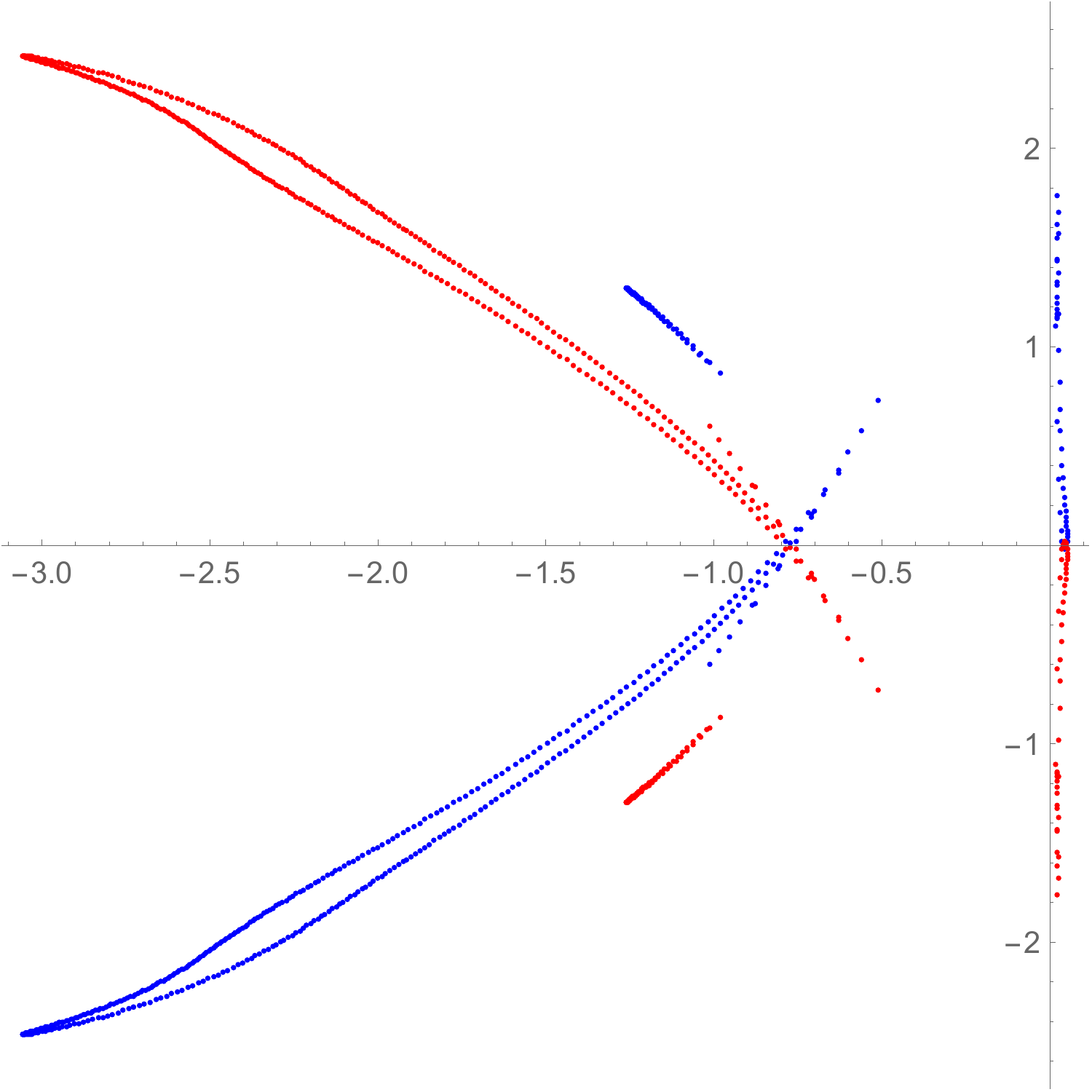}
  \includegraphics[width=0.32 \textwidth]{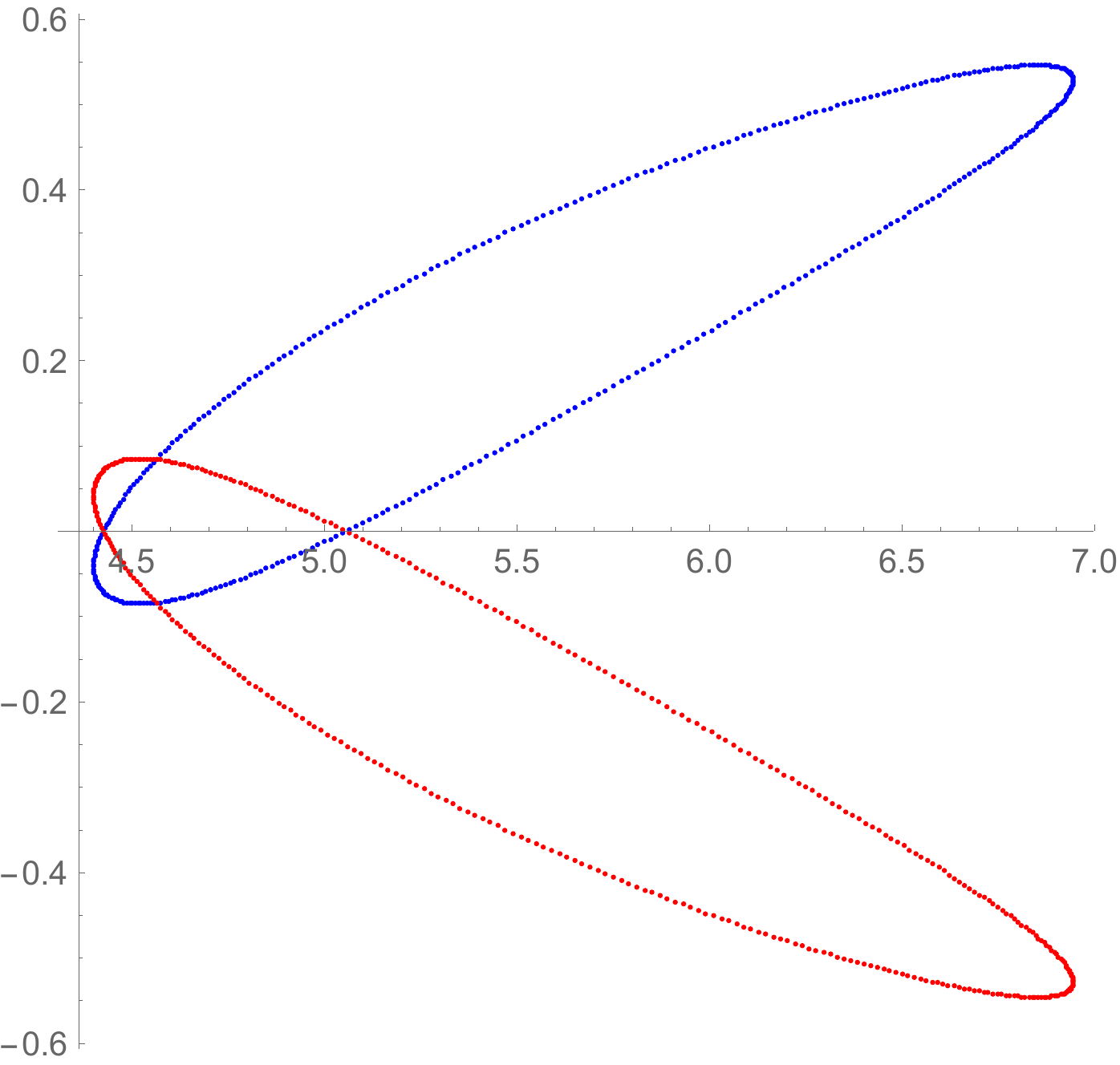}
      \caption{First intersections between the invariant manifolds for:  inner region $h=-2.15$ (left); inner region  $h=-2.075$ (center); outer region $h=-2.0$ (right).}
      \label{firstcuts}
\end{figure}

\begin{figure}[h]
  \centering
  \includegraphics[height=0.32\textwidth,width=0.32\textwidth]{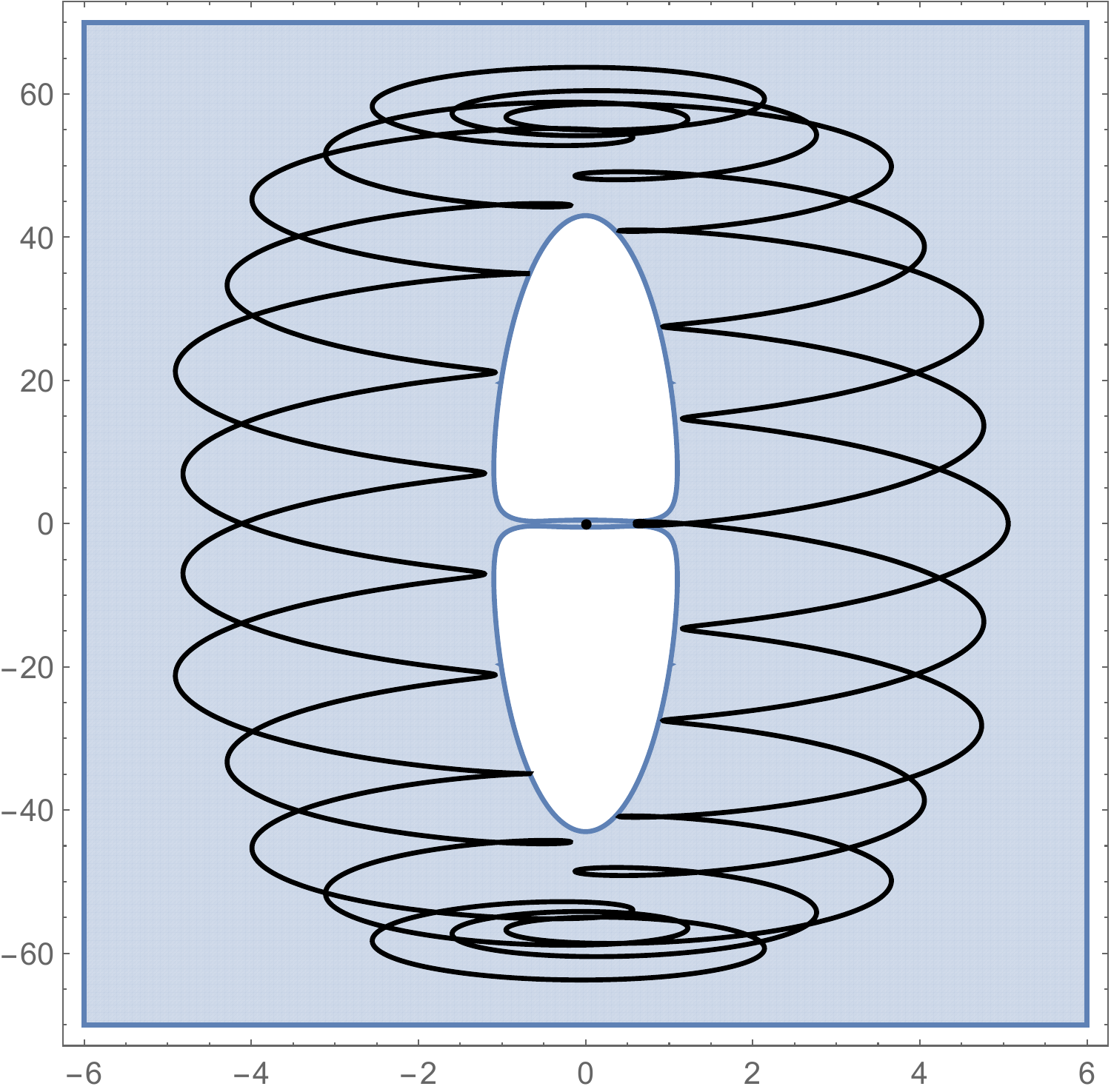}
  \includegraphics[height=0.32\textwidth,width=0.32\textwidth]{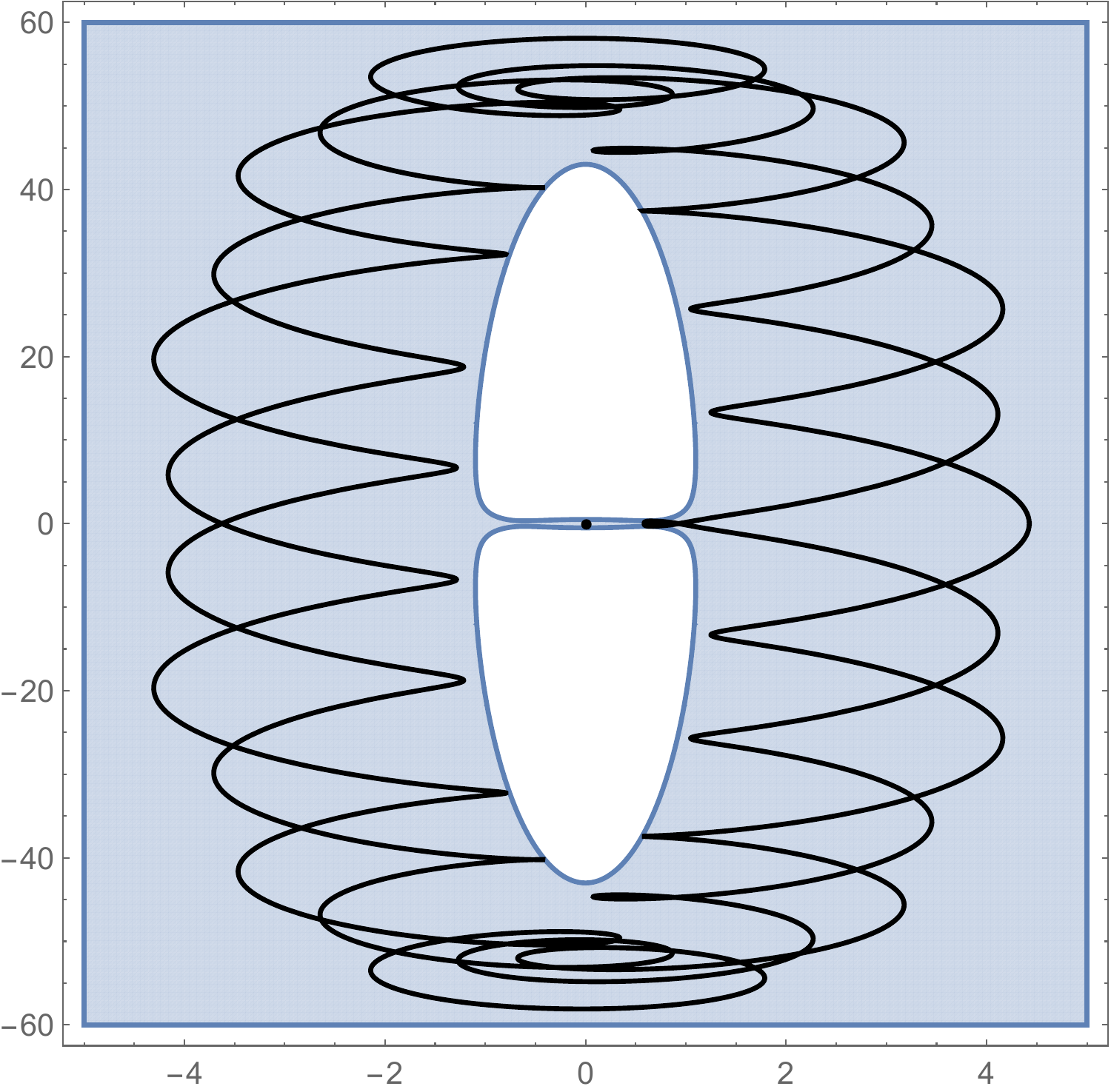}
  \includegraphics[width=0.32\textwidth]{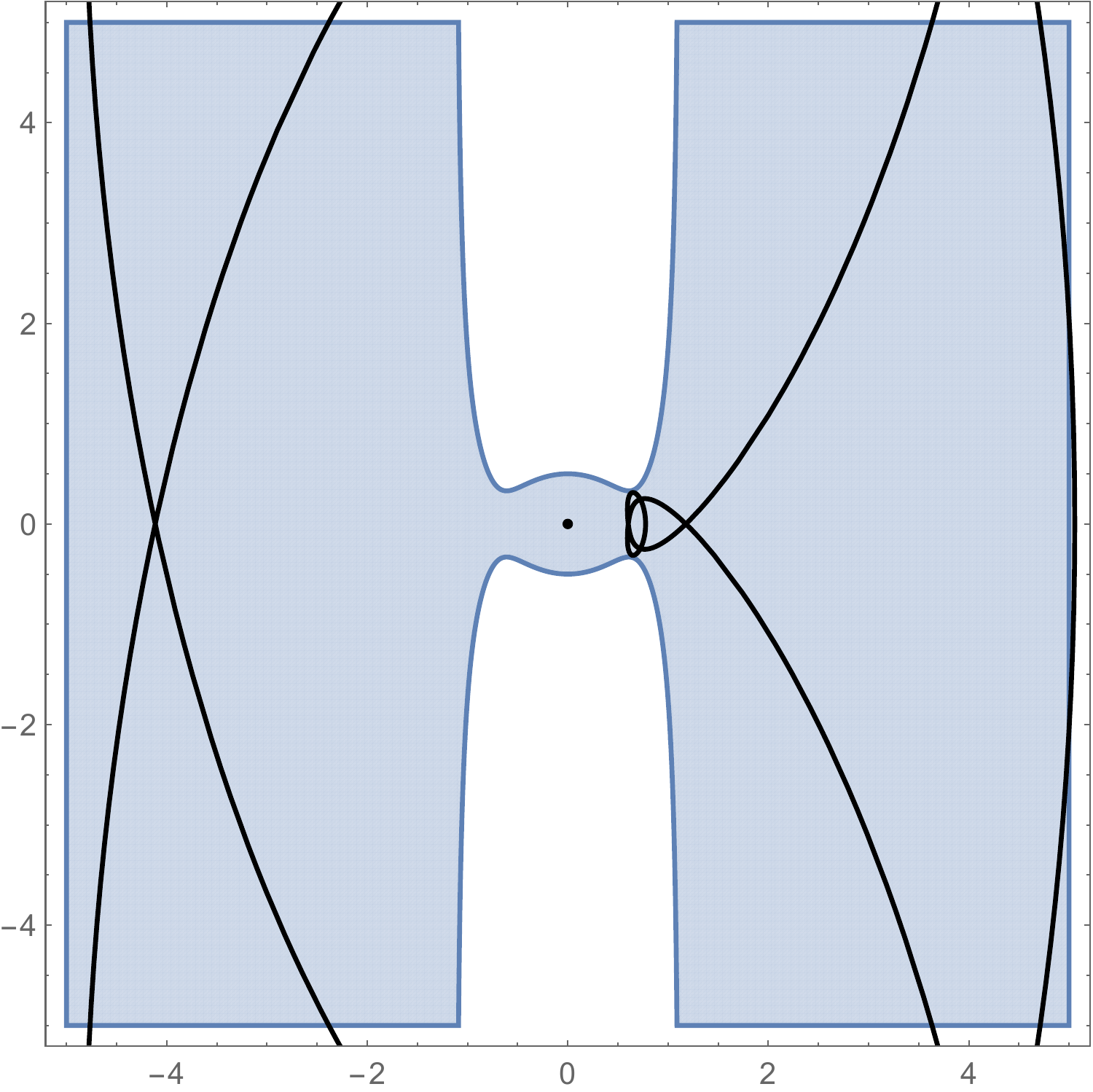}
\caption{Homoclinic orbits corresponding to the points $(5.05682901163445,0)$ and $(4.42411675944339,0)$ in the $(x,\dot{x})$-plane for $h=-2.0$. Magnification of the homoclinic orbit for $x_{0}=5.05682901163445$.}
\label{zoomhomoclinic1}
\end{figure}

We show the  evolution of the invariant manifolds in the plane $(x,\dot{x})$ for some values of the energy level $h$ in Fig. \ref{cuts}. It is notable that the invariant manifolds in the unbounded outer region do not escape to infinity, as in the lunar Hill's problem. In Fig. \ref {firstcuts} we show the first intersections between the invariant manifolds in the inner and outer region, where we note the presence of symmetric intersections $(\dot{x}=0)$ which correspond to homoclinic orbits symmetric with respect to the $x$-axis. In Fig. \ref{zoomhomoclinic1} we show two symmetric homoclinic orbits in the outer region. The computations were performed using a linear approximation around the periodic orbits and integrating numerically the vector field with tolerances for the error of order of machine epsilon.


\section{The mechanism of diffusion in the elliptic problem}
\label{sec:diffusion_EH4BP}
In Section \ref{sec:EH4BP} we showed that the elliptic problem for small values of the eccentricity can be written as
a perturbation of the  circular Hill problem, with the Hamiltonian given by \eqref{elliptichillhamiltonian_rot}.
Simplifying the notation, $H_\eps=H_{EH4BP}$, $H_0=H_{CH4BP}$, and $t=f$, we have:
\begin{equation}\label{perturbedsystem}
H_\eps(x,y,p_{x},p_{y},t)=H_{0}(x,y,p_{x},p_{y})+\eps G(x,y,t)+O(\eps^2),
\end{equation}
where $\eps$ is the eccentricity,
the time variable $t$ is $2\pi$-periodic, and the first order term of the perturbation is
\[G(x,y,t):=\cos(t)H_{1}(x,y),\]
where
$$H_{1}=\frac{1}{2}(x^{2}+y^{2})+U_{rot}(x,y),$$
with $U_{rot}(x,y)$ given by \eqref{eqn:U_rot}.
This is a  time-dependent perturbation of  $H_0$, with perturbation parameter $\eps$.


Below we formulate the main theoretical result on Arnold diffusion, noting that some of the
notions that appear in the statement will be defined in the later sections.

\begin{teo}\label{teo:main}
Consider the flow $\Phi^t_\eps$ of \eqref{perturbedsystem}.

Assume that the following conditions are satisfied for all  energy levels $h$ within some range $[h_\alpha, h_\beta]$ of $H_0$:
\begin{itemize}
\item[(i)] The unperturbed system for $\eps=0$ possesses two NHIMs $\Lambda^i_0$, consisting of families of periodic orbits $\lambda^i(h)$ around $L_i$, for $i=1,2$.
\item [(ii)] For $\eps=0$, each NHIM can be parameterized in terms of symplectic action-angle coordinates $(I,\theta)$, with the coordinate $I$  uniquely determined by  $h$ for both NHIMs.
\item[(iii)]  The corresponding stable and unstable manifolds $W^s(\Lambda^j_0)$, $W^u(\Lambda^i_0)$, intersects transversally for $i,j=1,2$.
\item [(iv)] There exists a collection of scattering maps associated to the transverse homoclinic connections $W^u(\Lambda^i_0)\pitchfork W^s(\Lambda^i_0)$, $i=1,2$,  as well as   a collection of scattering maps  associated to the transverse heteroclinic  connections  $W^u(\Lambda^i_0)\pitchfork W^s(\Lambda^j_0)$, $i\neq j$, such that, each scattering map expressed in the appropriate action-angle coordinate is of the form
    \[ (I,\theta)\mapsto (I,\theta+\Delta(I)).\]
\item [(v)]
Each perturbed  scattering map  in the $(I,\theta)$-coordinates  is of the form
\begin{equation*}
\begin{split}
\sigma_\eps(I,\theta)=&(I,\theta+\Delta(I))+\eps J\nabla{S}(I,\theta+\Delta(I))+O(\eps^2),
  \end{split}
\end{equation*}
for some Hamiltonian function $S=S(I,\theta)$, where $J=\begin{pmatrix} 0 &-1\\1& 0\end{pmatrix}$.

\item [(vi)]   Assume that the perturbed scattering maps satisfy
either the  conditions \eqref{eqn:domains}, \eqref{eqn:many_scatterings},
and \eqref{eqn:sum_integrals}  from  Theorem~\ref{prop:mechanism}, or the conditions   \eqref{eqn:many_scatterings}  and \eqref{eq:scatter-grad-a1} from Theorem \ref{th:mechanism-main}.
\end{itemize}


Then, there exist  $\eps_0>0$ and $C>0$, such that, for all $0<\eps<\eps_0$, there exists a diffusing orbit  $\Phi^t_\eps(z)$ of the perturbed system  and $T=T(\eps)>0$  such that
\[ | H_0(\Phi^T_\eps(z))-H_0(z)|> C.\]

Moreover,   there exist
such  diffusing orbits  that closely follow homoclinic connections, that is, they lie inside a small neighborhood of $W^u(\Lambda^1_0)\cap W^s(\Lambda^1_0)$,
or of $W^u(\Lambda^2_0)\cap W^s(\Lambda^2_0)$,  as well as diffusing orbits  that closely follow heteroclinic connections, that is, they lie inside a small neighborhoods of $(W^u(\Lambda^2_0)\cap W^s(\Lambda^1_0))\cup (W^u(\Lambda^1_0)\cup W^s(\Lambda^2_0)$.
\end{teo}

The NHIMs assumed in (i) will be described in Section \ref{sec:unperturbed_NHIM}.
The action-angle coordinates assumed in (ii) will be introduced in Section \ref{sec:coordinate}.
In Section \ref{sec:scattering}, we will first recall the  scattering map and the  formula for Hamiltonian function that generates it, and then we provide formulas for the scattering map in the case of the EH4BP.
In Section  \ref{sec:mechanism} and Section \ref{sec:mechanism_two} we provide two results on diffusion -- Theorem \ref{prop:mechanism} and  Theorem \ref{th:mechanism-main} --  from which  Theorem \ref{teo:main} immediately follows.

The numerical verification of the conditions of Theorem \ref{teo:main} is done in Section~\ref{sec:numerical_verification}.

To obtain diffusing orbits,  we will use a two-dynamics approach, by combining   the \emph{outer dynamics}, along homoclinic/heteroclinic  orbits to the NHIMs, encoded by the scattering map, with the  \emph{inner dynamics}, given by the restriction of the flow to the NHIMs.


\begin{rem} It is clear that the change of coordinates from the inertial non-rotating system to the pulsating system involves the perturbation parameter (similar changes of coordinates can be found in \cite{szebehely1967theory,belbruno2004capture}).
The main point is to obtain a system of the form (\ref{perturbedsystem}) in such a way that $H_{0}$ possesses a normal hyperbolic invariant manifold and transverse intersections between its stable and unstable manifolds. The model considered here satisfies these requirements. The effect of the parameter $\eps$ in the geometry of the coordinates could be relevant only if we consider the original non-rotating inertial system, which is not the case in the current work.\end{rem}

\subsection{The unperturbed NHIMs}
\label{sec:unperturbed_NHIM}
We consider the unperturbed system corresponding to $\eps=0$.
For each  center-saddle equilibrium point $L_i$, $i=1,2$, we consider a family of Lyapunov orbits
$\lambda_i(h)$ associated to energy levels $h$ within some suitable energy range $[h_\alpha,h_\beta]$ (to be specified in Section \ref{sec:Lyapunov}).

Due to the symmetry $S'(x,y,p_x,p_y)=(-x,y,p_x,-p_y)$, the Lyapunov orbits $\lambda_i(h)$, $i=1,2$, form a symmetric pair.
The Lyapunov orbits about $L_1$ and  $L_2$ are both traveled clockwise.
The existence of Lyapunov orbits about the points $L_i$, $i=1,2$, follows from the Lyapunov Center Theorem \cite{moser1958generalization}, for energy levels near the critical energy $h_{L_1}=H_0(L_1)=H_0(L_2)$. These orbits can be computed by numerical continuation for higher energy levels, up to collision with $m_3$.
Their existence can be established rigorously via computer assisted proofs \cite{Capinski_2012}.
In this paper we have non-rigorous computations of the Lyapunov orbits for the specific energy levels considered (see Section \ref{sec:Lyapunov}).

The family of Lyapunov orbits for the  energy range $[h_\alpha,h_\beta]$ determines
\begin{equation*}\label{eqn:NHIM_i}\begin{split}
  \Lambda^i_0=&\bigcup_{h\in[h_\alpha,h_\beta]}\lambda_i(h), \textrm { for } i=1,2,
\end{split}
\end{equation*}
which is a $2$-dimensional NHIM with boundary for the unperturbed flow $\Phi^t_0$ of $H_0$.

\subsection{Coordinate systems on the NHIMs and inner dynamics}
\label{sec:coordinate}
Below we provide some parameterizations for the NHIMs.
We focus on $\Lambda^1_0$; a similar parameterization can be obtained for $\Lambda^2_0$.
We refer to the restriction of the flow $\Phi^t_0$ to $\Lambda^1_0$ (or to $\Lambda^2_0$) as the \emph{inner dynamics}.

Each Lyapunov periodic orbit  $\lambda_1(h)$  intersects the $x$-axis in precisely two points, each of the form  $(q,p)=(x,0,0,p_{y})$, one on the left of $L_1$ and the other on  the right  of $L_1$.
Let $q^*=(x^*,0,0,p^*_{y})$ be the intersection between the Lyapunov orbit and the $x$-axis on the left of $L_1$, which is characterized by $p^*_{y}>0$.
The  orbit $\lambda_1(h)$ is uniquely determined by the coordinate $x^*=x^*(h)$, as the coordinate   $p^*_{y}=p^*_{y}(h)$ follows implicitly from the energy condition.
Note that  $h_\alpha<h_\beta$ implies $I_\alpha<I_\beta$ and $x^*(h_\beta)<x^*(h_\alpha)$ (since $x^*$ is the intersection point on the left of $L_1$).

Let $T(h)$ be the period of  $\lambda_1(h)$. Since the energy level is uniquely determined by $x^*$, we can also write $T(h)=T(x^*)$ for $x^*=x^*(h)$.
Then the Lyapunov periodic orbit is given by
\[\lambda_1(h)=\{\Phi^t_0(x^*(h),0,0,p^*_{y}(h))\,|\, t\in \mathbb{R}/(T(h)\mathbb{Z})  \}.\]

Then the NHIM $\Lambda^1_0$ can be described 
as
\begin{equation*}\label{eqn:NHIM_1}\begin{split}
  \Lambda^1_0=&\bigcup_{h\in[h_\alpha,h_\beta]}\lambda_1(h)\\
  =&\{\Phi^t_0(x^*(h),0,0,p^*_{y}(h))\,|\, t\in \mathbb{R}/(T(h)\mathbb{Z}) \textrm { for }  h\in [h_\alpha,h_\beta]\},
\end{split}
\end{equation*}
and can be parameterized by $(x^*,t^*)$ for $x^*\in [x^*(h_\beta), x^*( h_\alpha)]$ and $t^* \in  \mathbb{R}/ (T(h)\mathbb{Z}) $, where  $x^*=x^*(h)$ for $h\in[h_\alpha,h_\beta]$.
The  coordinate system $(x^*,t^*)$ is not symplectic on $\Lambda^1_0$.

We define symplectic coordinates $(I,\theta)$ on $\Lambda^1_0$ as follows.
Let $\omega_{\mid\Lambda^1_0} =(dp\wedge dq)_{\mid\Lambda^1_0}:=(dp_x\wedge dx+dp_y\wedge dy)_{\mid\Lambda^1_0}$ be the induced symplectic form on
 $\Lambda^1_0$. Since the induced form is non-degenerate, it  endows $\Lambda^1_0$ with a symplectic structure  \cite{DelshamsLS08a}.
 Let
\[\theta =\frac{t^*}{T(x^*)}\textrm { and } I=\int_{\lambda^1_0(x^*)} pdq.\]
The action
$I = I(h)$ is uniquely determined by the energy level $h$, and corresponds to a
unique Lyapunov orbit.
We then have
\[dI\wedge d\theta=\omega_{\mid\Lambda^1_0}.\]
The condition $h\in[h_\alpha,h_\beta]$ is equivalent to $I\in [I_\alpha,I_\beta]$ for $I_\alpha=I(h_\alpha)=I(x^*(h_\alpha)) $ and  $I_\beta=I(h_\beta)=I(x^*(h_\beta)) $.
Each $I\in [I_\alpha,I_\beta]$ is in unique correspondence with  some  $h\in [h_\alpha,h_\beta]$ as well as with some $x^*\in [x^*(h_\beta),x^*(h_\alpha)]$, so we can write $T(I)=T(x^*(I))$.

We  define the reference manifold
\begin{equation}
\label{eqn:reference_manifold}
\mathbb{A}=\{(I,\theta)\,|\, I\in [I_\alpha,I_\beta],\, \theta\in \mathbb{T}^1\}
\end{equation}
and the parameterization of $\Lambda^1_0$ given by
\[k_{0}:\mathbb{A}\to\Lambda^1_0, \textrm { where }  k_{0}(I,\theta)=\Phi_0 ^{\theta T(I)}(x^*(I),0,0, p^*_{y}(I)).\]
The parameterization $k_0$ induces a flow $R^t$ on $\mathbb{A}$ given by
\[R^t(I,\theta)=(I,\theta+ t \omega(I)) \textrm{ where } \omega(I)=\frac{1}{T(I)} .\]
We have the following conjugacy
\[ \begin{tikzcd}
\mathbb{A} \arrow{r}{R^t} \arrow[swap]{d}{k_{0}} & \mathbb{A} \arrow{d}{k_{0}} \\%
\Lambda^1_0 \arrow{r}{\Phi^t_0}& \Lambda^1_0
\end{tikzcd}
\]

When we  consider both NHIMs $\Lambda^1_0$,  $\Lambda^2_0$, we have the corresponding $(x^*_1,t^*_1)$ and $(x^*_2,t^*_2)$-coordinates, and the action-angle coordinates $(I_1, \theta_1)$ and $(I_2,\theta_2)$.
By the symmetry of the Lyapunov orbits $x^*_2=-x^*_1$.
The angle $\theta_1$ is measured clockwise starting from the point $(x^*_1, 0, 0, p^*_{1,y})$, and the angle $\theta_2$ is measured clockwise starting from the point $(x^*_2, 0, 0, p^*_{2,y})$. 
Since the action $I$ is uniquely determined by the energy level $H_0=h$, the action coordinate is the same $I_1=I_2:=I$ for both NHIMs.

\subsection{The scattering map}
\label{sec:scattering}
In Sections \ref{sec:unperturbed_NHIM}  we have identified two NHIMs, as well as  homoclinic and heteroclinic connections between them.
Moving along the homoclinic / heteroclinic orbits represents the \emph{outer dynamics}.
The outer dynamics can be described by  scattering maps.
First we review the  definition of the scattering map and its properties, following \cite{DelshamsLS08a}.
Second, we  describe the scattering map for the unperturbed problem, that is, in the context of the CH4BP.

\subsubsection{The scattering map: definition and properties}
\label{sec:scattering_def}
Consider the general case of a normally hyperbolic invariant manifold $\Lambda$ for a flow $\Phi^t$ on some smooth manifold $M$. We assume that $M$ and $\Phi^t$ are sufficiently regular; explicit  regularity conditions can be found in \cite{DelshamsLS08a}.  Let $W^s(\Lambda)$, $W^u(\Lambda)$ be the stable and unstable manifolds of $\Lambda$.
Define the wave maps  as  canonical projections along fibers, i.e.,
\begin{equation*}
\begin{split}
  \Omega^+:W^s(\Lambda)\to\Lambda,\\
  \Omega^-:W^u(\Lambda)\to\Lambda,
\end{split}
\end{equation*}
where
$\Omega^+$ assigns to $z\in W^s(\Lambda)$ its stable foot point $z^{+}=\Omega^{+}(z)$, uniquely defined by
$z\in W^s(z^+)$, and, similarly,  $\Omega^-$ assigns to $z\in W^u(\Lambda)$ its unstable footpoint $z^{-}=\Omega^{-}(z)$,  uniquely defined by
$z\in W^u(z^-)$.

The wave maps  satisfy the following equivariance property
\begin{equation}\label{eqn:equivariance_wave}
 \Omega^{\pm} \circ\Phi^t=\Phi^t\circ\Omega^{\pm} .
\end{equation}

A  homoclinic channel $ \Gamma$ is a submanifold in $W^u( \Lambda)\cap
W^s( \Lambda)$ that  satisfies the following strong transversality conditions for all  $z\in\Gamma$:
\begin{equation}
\label{eqn:strong transversality}
\begin{split}
T_z\Gamma=&T_zW^s(\Lambda)\cap
T_zW^u(\Lambda),
\\
T_zM=&T_z\Gamma \oplus T_zW^u(z^-)\oplus  T_xW^s(z^+),
\end{split}\end{equation}
and such that  \[\Omega^\pm_{\mid \Gamma}:\Gamma\to \Omega^\pm(\Gamma) \textrm{ is a diffeomorphism. }\]

Then, the scattering map associated to the homoclinic channel $\Gamma$ is the mapping
\begin{equation*}\begin{split}
\sigma=&\sigma^\Gamma:\Omega^{-}(\Gamma)\subseteq \Lambda\to \Omega^{+}(\Gamma)\subseteq \Lambda,\\
\sigma=&\Omega^+_{\mid \Gamma}\circ (\Omega_{\mid \Gamma}^{-})^{-1}.
\end{split}
\end{equation*}
The map $\sigma$ is a diffeomorphism   from its domain onto its image in $\Lambda$. It is not globally defined, in general.

The regularity of $\sigma$ is the same as the regularity of $W^u( \Lambda)$ and $W^s( \Lambda)$, which in turn depends on the regularity of the system as well as on the hyperbolic rates. For precise conditions on regularity, see \cite{DelshamsLS08a}. Since we  assume  that our system is sufficiently regular, we also have that $\sigma$ is sufficiently regular.

We have that $\sigma( z^-)=z^+$ if and only if
\begin{equation*}\label{eqn:scattering_homoclinic}
d(\Phi^{-T^-}(z),\Phi^{-T^-}(z^-))\to~0,\textrm{ and  }d(\Phi^{T^+}(z),\Phi^{T^+}(z^+))\to~0
\end{equation*} as $T^-,T^+\to+\infty$, respectively, for some uniquely defined $z\in\Gamma$.

The equivariance of the wave maps yields the following equivariance property of the scattering map:
\begin{equation*}\label{eqn:equivariance_scattering}
 \Phi^t\circ\sigma^\Gamma=\sigma^{\Phi^t(\Gamma)}\circ \Phi^t .
\end{equation*}

An important property of the scattering map is that  $\sigma$ is symplectic, and, moreover, is exact symplectic provided that  $M$ is endowed with an exact form \cite{DelshamsLS08a}.

\begin{rem}\label{rem:monodromy} If  the strong transversality condition \eqref{eqn:strong transversality} is satisfied at a point $z_0$, then it is satisfied in a small neighborhood $\Gamma_0$ of $z_0$ in $W^u( \Lambda)\cap W^s( \Lambda)$. One can take $\Gamma_0$ as a homoclinic channel and define the associated scattering map, and then extend  $\Gamma_0$ to a maximal domain on which the wave maps $\Omega^\pm_{\mid \Gamma_0}$ are diffeomorphisms. Extending $\Gamma_0$ beyond its maximal domain may result  in failure of monodromy. Monodromy here means that when the scattering map is applied to $(I,\theta)$ and to $(I,\theta+1)$ it gives the same value,  in other words the scattering map is  well defined as a map on the annulus, since $\theta$ is defined mod $1$; this does not happen in general, see \cite{CapinskiGL17}.

For example, consider the  homoclinic orbit $\Phi^t(z_0)$; each point $z$ on this homoclinic orbit is of the form  $z=\Phi^t(z_0)$ for some $t$. The equivariance property \eqref{eqn:equivariance_wave} implies
\begin{equation*}\label{eqn:equivariance-homoc}\Omega^{\pm}(z)=\Omega^ {\pm}(\Phi^t(z_0))=\Phi^t\circ\Omega^{\pm}(z_0)=\Phi^t(z_0^\pm).\end{equation*}
Therefore, it is  possible to have a pair of homoclinic points
$z_1=\Phi^{t_1}(z_0)$, $z_2=\Phi^{t_2}(z_0)$  such that $\Omega^{-}(z_1)=\Omega^{-}(z_2)$ and $\Omega^{+}(z_1)\neq\Omega^{+}(z_2)$, so $\Omega^+ \circ (\Omega ^{-})^{-1}$ is not a well-defined map. Restricting to a suitable $\Gamma_0$ will ensure the monodromy property.

An alternative  is to define a \emph{scattering correspondence} rather than a scattering map, which allows  for assigning to a given point $z^-$ multiple points $z^+$; see \cite{gidea_marco_2017}.
\end{rem}

The scattering map can be defined in a similar way for heteroclinic connections between two NHIMs $\Lambda^1$ and $\Lambda^2$.
Assuming that $W^u(\Lambda^2)$ intersects transversally $W^s(\Lambda^1)$, we can define  wave  maps
\begin{equation*}
\begin{split}
  \Omega^{1,+}:W^s(\Lambda_1)\to\Lambda_1,\\
  \Omega^{2,-}:W^u(\Lambda_2)\to\Lambda_2,
\end{split}
\end{equation*}
and a  heteroclinic channel $\Gamma \subseteq  W^u(\Lambda^2)\cap W^s(\Lambda^1)$ in a similar fashion to the homoclinic case.
Then we define
\[\sigma =\Omega_{\mid \Gamma}^{1,+}\circ (\Omega_{\mid \Gamma}^{2,-})^{-1}.\]
The scattering map in the heteroclinic case enjoys similar properties to the one for the scattering map in the homoclinic case.

We now consider the scattering map in the case of a perturbed system.
Consider a  family of NHIMs $\Lambda_\eps$   for $\Phi^t_\eps$ on $M$, for  $0\leq \eps<\eps_0$. We will implicitly assume that $\Phi^t_\eps$ and $\Lambda_\eps$   are sufficiently regular in the variables and in the parameter $\eps$. Here we view $\Phi^t_\eps$ as a perturbation of $\Phi^t_0$. Below, we will also consider the special case when the flow $\Phi^t_\eps$  is the flow of a  perturbed Hamiltonian system
$H_\eps$.

Assume that $\Lambda_\eps$ has a smooth parameterization in terms of some reference manifold $\mathbb{A}$, that is, there exists a  family of diffeomorphisms (in a suitable regularity class)  \[k_\eps :\mathbb{A} \to \Lambda_\eps= k_\eps(\mathbb{A}), \textrm{ for  $0\leq \eps<\eps_0$.}\]

Assume that $W^u(\Lambda_\eps)$ and $W^s(\Lambda_\eps)$ intersect transversally
along a homoclinic channel  $\Gamma_\eps$. Let $\sigma_\eps:\Omega^{-}(\Gamma_\eps)\subseteq \Lambda_\eps\to \Omega^{+}(\Gamma_\eps)\subseteq \Lambda_\eps$ be the scattering map associated to $\Gamma_\eps$.

This is induces a locally defined diffeomorphism  on $\mathbb{A}$ given by
\[s_\eps=k_\eps^{-1}\circ\sigma_\eps\circ k_\eps.\]
We will also refer to $s_\eps$ as the scattering map.
Note that for all $\eps$ the map $s_\eps$ is a diffeomorphism on some subset  of the same manifold $\mathbb{A}$, in contrast with $\sigma_\eps$, which is
a diffeomorphism on some subset  of a manifold $\Lambda_\eps$ that depends on $\eps$. The advantage of using the scattering map $s_\eps$ induced on the reference manifold $\mathbb{A}$ is that one can compare $s_\eps$'s for different parameter values $\eps$.

It was shown in \cite{DelshamsLS08a} that the perturbed scattering map $s_\eps$ is smooth and depends smoothly on the parameter $\eps$, and can be expanded in powers of $\eps$ as:
\begin{equation}
\label{eqn:scattering_expansion}\begin{split}
s_\eps(I,\theta)=&s_0(I,\theta)+ \eps J \nabla S  \circ s_0(I,\theta)+O(\eps^2)\\
=&s_0(I,\theta)+\eps\left(-\frac{\partial S}{\partial \theta},\frac{\partial S}{\partial I}\right)\circ s_0 (I,\theta)+O(\eps^2)
\end{split}\end{equation}
for some Hamiltonian function $S$ defined on some domain in $\mathbb{A}$. Above, $s_0$ is the scattering map induced by $\sigma_0$, corresponding to the unperturbed system.

In the perturbative setting, $S$ can be computed explicitly.
In particular, in the case  of a perturbed Hamiltonian system
$H_\eps$,
$S$  can be computed as a Melnikov integral (see \cite{DelshamsLS08a}):
\begin{equation}\label{eqn:S0_z_0^+}\begin{split}
S(k_0^{-1}(z^+))=&\int_{-\infty}^{0}\left[\frac{\partial H_\eps}{\partial \eps}_{\mid\eps=0}\circ\Phi_0^{t} (z)  - \frac{\partial H_\eps}{\partial \eps}_{\mid\eps=0}\circ\Phi_0^{t} (z^-)\right]\,dt\\
&+ \int_{0}^{+\infty}\left[\frac{\partial H_\eps}{\partial \eps}_{\mid\eps=0} \circ\Phi_0^{t}  (z) -\frac{\partial H_\eps}{\partial \eps}_{\mid\eps=0} \circ\Phi_0^{t}   (z^+)\right ]\,dt,
\end{split}
\end{equation}
where $z$ is a homoclinic point for the unperturbed system.

The integrals in \eqref{eqn:S0_z_0^+} are  improper integrals  whose integrand is given by the difference between the perturbation evaluated on homoclinic orbits of the unperturbed system and  the perturbation evaluated on the asymptotic orbits  on the unperturbed NHIM. The  improper integrals converge exponentially fast, so they can be  efficiently  computed via numerical methods.

In the context of Section \ref {sec:coordinate}, if $k_0^{-1}(z^+)=(I,\theta)$, then we can write $S$ in terms of $(I,\theta)$ as
\begin{equation}\label{eqn:S0_I_theta}\begin{split}
S(I,\theta)=&\int_{-\infty}^{0}\left[\frac{\partial H_\eps}{\partial \eps}_{\mid\eps=0}\circ\Phi_0^{t} \circ (\Omega^-)^{-1}\circ\sigma_0^{-1}\circ k_0(I,\theta) \right.\\
&\qquad\qquad \left.- \frac{\partial H_\eps}{\partial \eps}_{\mid\eps=0}\circ\Phi_0^{t} \circ\sigma_0^{-1}\circ k_0(I,\theta)\right]\,dt\\
&+ \int_{0}^{+\infty}\left[\frac{\partial H_\eps}{\partial \eps}_{\mid\eps=0} \circ\Phi_0^{t}  \circ (\Omega^+)^{-1} \circ k_0(I,\theta)\right.\\
&\qquad\qquad\left.-\frac{\partial H_\eps}{\partial \eps}_{\mid\eps=0} \circ\Phi_0^{t}   \circ k_0(I,\theta)\right ]\,dt,
\end{split}
\end{equation}

It is clear that \eqref{eqn:S0_I_theta} is equivalent to \eqref{eqn:S0_z_0^+} since, for $k_0^{-1}(z^+)=(I,\theta)$, we have
\begin{equation*}
\begin{split}
 &k_0(I,\theta)=z^+,\\
 &\sigma_0^{-1}\circ k_0(I,\theta)=z^-,\\
 &(\Omega^-)^{-1}\circ\sigma_0^{-1}\circ k_0(I,\theta)=
 (\Omega^+)^{-1} \circ k_0(I,\theta)=z.
\end{split}
\end{equation*}

\subsubsection{The unperturbed scattering map associated to homoclinic  connections for the CH4BP}
\label{sec:scattering_hom}
First, we show that, in general, the unperturbed scattering map associated to a homoclinic channel is a phase shift in the angle variable.

For a given range of energies $h\in [h_\alpha,h_\beta]$, let $z=z(h)$
be a transverse homoclinic point,
$z^\pm$ be its stable and unstable foot points, $\theta^\pm$ their corresponding angle coordinates,
and  $T=T(h)$  the period of the Lyapunov orbit $\lambda(h)$.

We can define a homoclinic channel associated to this homoclinic point $z$ by
\begin{equation}\label{eqn:unperurbed_homoclininc_channel}
  \Gamma=\bigcup_{\substack{ h\in[h_\alpha,h_\beta]\\ \theta\in (\theta^-+\theta^*, \theta^-+1+\theta^*) }}\Phi_0^{(\theta-\theta^-)T} (z),
\end{equation}
for some  $\theta^*\in\mathbb{R}$. (As we will see in Lemma \ref{lem:Delta}, the corresponding scattering map is globally defined and, in fact,
in \eqref{eqn:unperurbed_homoclininc_channel} we can choose the range of $\theta$ to be $[\theta^-+\theta^*, \theta^-+1+\theta^*]$.)

The flow $\Phi^t_0$ acts on the fiber $W^u(z^-)=W^u(k_0(I,\theta^-))$ by shifting its base point
$k_0(I,\theta^-)$ to $k_0(R^t(I,\theta^-))=k_0(I,\theta^-+\frac{t}{T(I)})$, and on the fiber $W^s(z^+)=W^s(k_0(I,\theta^+))$ by shifting its base point
$k_0(I,\theta^+)$ to $k_0(R^t(I,\theta^+))=k_0(I,\theta^++\frac{t}{T(I)})$.
The flow $\Phi^t_0$ takes $z\in W^u(z^-)\cap W^s(z^+)$ to $\Phi^t_0(z)\in W^u(\Phi^t_0(z^-))\cap W^s(\Phi^t_0(z^+))$.
In particular, for $t=(\theta-\theta^-)T$, we have
\begin{equation*}\begin{split}
 \Phi^{(\theta-\theta^-)T}_0(W^u(z^ -))=&W^u(k_0(I,\theta)),\\
  \Phi^{(\theta-\theta^-)T}_0(W^s(z^+))=&W^s(k_0(I,\theta+\theta^+-\theta^-)).
\end{split}\end{equation*}
By the $S$-symmetry of the flow on the Lyapunov orbit we also have
 \begin{equation*}
 \begin{split}
 \Phi_0^{(1+\theta-\theta^-)T} (W^{u}(z^-)) =&\Phi_0^{(\theta-\theta^-)T}\circ \Phi_0^T (W^{u}(z^-))=\Phi_0^{(\theta-\theta^-)T}(W^{u}( \Phi_0^T (z^-))\\
 =&\Phi_0^{(\theta-\theta^-)T}(W^{u}( z^-)).
\end{split}
\end{equation*}
and similarly
 \begin{equation*}
 \begin{split}
 \Phi_0^{(1+\theta-\theta^+)T} (W^{s}(z^+) )=\Phi_0^{(\theta-\theta^+)T}(W^{s}( z^+)).
\end{split}
\end{equation*}
The following result from \cite{gidea2022melnikov} gives a complete description of the unperturbed scattering map.

\begin{lem} \label{lem:Delta}
Let $z=z(h)$ be a transverse homoclinic point, $\Gamma$ be a homoclinic channel defined as in \eqref{eqn:unperurbed_homoclininc_channel}, $z^\pm=\Omega_\Gamma^\pm(z)$, and  $(I,\theta ^\pm)$ be the action-angle coordinates of $z^\pm$.

Then the unperturbed scattering map is globally defined, and is given in the action-angle $(I,\theta)$-coordinates by
a shift in the angle coordinate, i.e.,
\[s_0(I,\theta)=(I,\theta+\Delta )=(I,\theta+(\theta^+-\theta^-)) , \]
where $I=I(h)$ and $\Delta=\Delta(h)$ depend on the energy level $h$.
\end{lem}

\begin{proof}
The points $z^\pm$ are  on the Lyapunov orbit $\lambda(h)$.
We have that $\sigma(z^-)=z^+$ is given in local coordinates
by \[s_0(I,\theta^-)=(I,\theta^+):=(I,\theta^-+\Delta),\]  for $\Delta=\Delta(h)=\theta^+-\theta^-$.

If $z_1\in\Gamma\cap \{H_0(z)=h\}$ is another  homoclinic point, then  $z_1=\Phi^{t_1}_0(z)$ for some $t_1$. Let $(I,\theta_1^\pm)$ be the action-angle coordinates of $z_1^\pm$, respectively.
By the equivariance property \eqref{eqn:equivariance_wave}, $z_1^\pm=\Omega_\Gamma^\pm(z_1)=\Omega_\Gamma^\pm(\Phi^{t_1}_0(z))=\Phi^{t_1}_0(\Omega_\Gamma^\pm(z))=\Phi^{t_1}_0(z^\pm)$, therefore
$\theta_1^\pm=\theta^\pm +\frac{1}{T}t_1$, so $\theta_1^+-\theta_1^-=\theta^+-\theta^-$. Thus
\[ s_0(I,\theta_1^-)= (I,\theta_1^-+\Delta).\]

As the mapping $s_0(I,\theta)$ is periodic of period $1$ in the angle variable, the homoclinic channel $\Gamma$ in \eqref{eqn:unperurbed_homoclininc_channel} can be extended to $\theta\in[\theta^-+\theta^*, \theta^-+1+\theta^*]$  so that $s_0$ is globally defined.
\end{proof}

The fact that the scattering map is a shift in the angle coordinate was verified numerically in \cite{canalias2006scattering}.

In Section \ref{sec:homoclinic connections} we will choose some symmetric homoclinic points  and associate to them homoclinic channels defined as in \eqref{eqn:unperurbed_homoclininc_channel}. We will then compute numerically the angle shifts for the corresponding scattering maps.

For a symmetric homoclinic point $z$, by the $S$-symmetry, the angle coordinates of the stable and unstable foot-points $z^\pm$ satisfy
\[\theta^+=-\theta^-.\] We then have:

\begin{cor}\label{cor:Delta_symmetric}
If $z$ is a symmetric homoclinic point, then $\Delta=-2\theta^-$, so, the unperturbed scattering map is given by
\[s_0(I,\theta)=(I,\theta-2\theta^-). \]
\end{cor}

\subsubsection{The unperturbed NHIMs and scattering map in the extended phase space}
\label{sec:unperturbed_NHIM_extended}
Since the perturbation $G(t)=\eps\cos(t)H_1+O(\eps^2)$ is time-periodic, it is convenient to consider time as an additional variable which we denote by~$s$.
The extended phase space is $M\times \mathbb{T}^1$. The equations of motion in the extended phase space are augmented with an equation for the evolution of the new variable $\frac{ds}{dt} =1$. Denoting by $\tilde\Phi^t_0$ the flow in the extended space, we have
$\tilde{\Phi}^t_0(z,s)=(\Phi^t_0(z),s+t)$.

In the extended space, the NHIM $\Lambda^i_0$ corresponds to the NHIM
\[\tilde\Lambda^i_0=\Lambda^i_0\times \mathbb{T}^1.\]

Similarly, any homoclinic/heteroclinic channel $\Gamma$ gives rise to a homoclinic/heteroclinic channel in the extended space
\begin{equation}\label{eqn:unperturbed homoclinic channel extended}
 \tilde\Gamma_0=\Gamma_0\times \mathbb{T}^1
\end{equation}

The scattering map corresponding to $\sigma_0$ in the extended space is given by
\[\tilde\sigma_0(z,s)=(\sigma_0(z),s)\]

%
%

\subsubsection{The perturbed NHIMs}
\label{sec:perturbed_NHIM_extended}
We now consider the effect of the perturbation $G(t)=\eps\cos(s)H_1+O(\eps^2)$ on the NHIMs in the extended phase space.

The standard theory of normally hyperbolic invariant manifolds,
\cite{Fenichel71,HirschPS77} shows the persistence of compact  NHIMs  (without boundary) under small perturbations.
In the case of manifolds with
boundary, the theory only guarantees the persistence of a normally hyperbolic manifold that is  locally invariant but not necessarily unique
(see \cite{BatesLZ00,berger2013geometrical,eldering2013normally}).
The proof in that case involves extending the vector
field in such a way that the manifold we consider is an invariant manifold
without boundary. Then, applying the result of persistence of an invariant
manifold without boundary, one obtains the existence of a locally invariant
manifold. The persistent manifold  depends on the extension considered, and hence is not unique.
However, all orbits that remain in a small neighborhood of
the manifold and away from its boundary remain present in all extensions
that do not modify the dynamics in that neighborhood.

This theory implies that  there exists $\eps_0$ such that the NHIM $\tilde{\Lambda}^i_0$, $i=1,2$,
persists as a  normally hyperbolic manifold $\tilde{\Lambda}^i_\eps$, which is locally invariant under the flow $\tilde{\Phi}^t_\eps$, for all $0<\eps <\eps_0$.
Moreover, there exists  smooth parameterizations of these manifolds
(see \cite{DelshamsLS08a}):

\[\tilde{k}_\eps:\mathbb{A}\times\mathbb{T}^1\to\tilde{\Lambda}_\eps.\]

Choosing $\eps_0$  sufficiently small ensures that the transverse homoclinic/heteroclinic channels for the unperturbed system also persist as homoclinic/heteroclinic channels $\tilde{\Gamma}_\eps$  for $\tilde{\Phi}^t_\eps$. Therefore, we have a scattering map  $\tilde{\sigma}_\eps$ associated to each channel. To compute the perturbed scattering map, we use \eqref{eqn:scattering_expansion}.
Since we work in the extended space the formula \eqref{eqn:S0_I_theta} for $S$ becomes

\begin{equation}\label{eqn:S0_extended}\begin{split}
S(I,\theta,s)=&\int_{-\infty}^{0}\left[\left(\frac{\partial H_\eps}{\partial \eps}_{\mid\eps=0}\circ\Phi_0^{t} \circ (\Omega^-)^{-1}\circ\sigma_0^{-1}\circ k_0(I,\theta),s+t\right) \right.\\
&\qquad\qquad \left.- \left(\frac{\partial H_\eps}{\partial \eps}_{\mid\eps=0}\circ\Phi_0^{t} \circ\sigma_0^{-1}\circ k_0(I,\theta),s+t\right)\right]\,dt\\
&+ \int_{0}^{+\infty}\left[\left(\frac{\partial H_\eps}{\partial \eps}_{\mid\eps=0} \circ\Phi_0^{t}  \circ (\Omega^+)^{-1} \circ k_0(I,\theta),s+t\right)\right.\\
&\qquad\qquad\left.-\left(\frac{\partial H_\eps}{\partial \eps}_{\mid\eps=0} \circ\Phi_0^{t}   \circ k_0(I,\theta),s+t\right)\right ]\,dt .
\end{split}
\end{equation}

In the subsequent sections, we will compute the  perturbed scattering map for the EH4BP.

\subsubsection{The perturbed scattering map in the EH4BP}
\label{sec:perturbed_scattering}

In what follows, we will always assume that $z$ is a symmetric homoclinic point, and associate to it a homoclinic channel defined as in
\eqref{eqn:unperurbed_homoclininc_channel}.

With the notation from Section~\ref{sec:scattering_hom} we have
\begin{equation*}\begin{split} k_0^{-1}(\Omega^{-}_{0}(\Phi^{(\theta-\theta^-)T(I)} (z)))&=(I,\theta),\\
k_0^{-1}( \Omega^{+}_{0}(\Phi^{(\theta-\theta^-)T(I)} (z)))&= (I,\theta+\Delta),\\
s_0(I,\theta)&=(I,\theta+\Delta),\\
\Delta&=-2\theta^- .
\end{split}\end{equation*}
The last equation follows from Corollary \ref{cor:Delta_symmetric}.

We then consider that homoclinic channel \eqref{eqn:unperturbed homoclinic channel extended}  for the unperturbed problem in the extended phase space.
The corresponding scattering map is given by
\[\tilde{s}_0(I,\theta,s)=(I,\theta+\Delta,s).\]

To compute the perturbed scattering map $\tilde{\sigma}_\eps$ associated to the homoclinic channel $\tilde\Gamma_\eps$,
we use \eqref{eqn:scattering_expansion} where
$S$ is given by \eqref{eqn:S0_extended}.

We  recall  that
\begin{equation*}
\begin{split}
\frac{\partial H_\eps}{\partial \eps}_{\mid \eps=0}(z,s)=&G(z,s)=\cos(s)H_1(z)\\
(\Omega^{-}_{0})^{-1}(k_0( s_0^{-1} (I,\theta)))=&(\Omega^{-}_{0})^{-1}(k_0(I,\theta-\Delta))=\Phi^{(\theta+\theta^-)T(I)}(z),\\
s_0^{-1} (I,\theta)= &(I,\theta-\Delta)=(I,\theta+2\theta^-),\\
(\Omega^{+}_{0})^{-1}(k_0  (I,\theta))=&\Phi_0^{(\theta+\theta^-)T(I)}(z).
\end{split}
\end{equation*}

Substituting in \eqref{eqn:S0_I_theta} we obtain
\begin{equation}\label{eqn:S0-1}\begin{split}
S(I,\theta,s)
=&\int_{-\infty}^{0}\cos(t+s)\left[H_1(\Phi_0^{t+(\theta+\theta^-)T(I)}(z))\right.\\
&\qquad\qquad \left.-H_1(\Phi_0^{t}(k_0(I,\theta+2\theta^-)))\right]dt\\
&+ \int_{0}^{+\infty}\cos(t+s) \left[H_1(\Phi_0^{t+(\theta+\theta^-)T(I)}(z))\right.\\
&\qquad\qquad\left.-H_1(\Phi_0^{t}(k_0(I,\theta)))\right ]dt .
\end{split}
\end{equation}

For our particular case of the EH4BP we have the following:

\begin{prop}\label{prop:S0-2} In the EH4BP, the   derivative with respect to the angle variable of the Hamiltonian that generates the scattering map in the extended phase space is given by the expression
\begin{equation}
\label{eqn:S0_2}
\begin{split}
\frac{d}{d\theta}S(s_0(I,\theta),s)&= T(I)\cos(s)\left [ H_1(\Phi_0^{(\theta-\theta^-)T(I)}(k_0(I,  \theta^+)))\right.\\
&\qquad\qquad\qquad\qquad\left.-H_1(\Phi_0^{(\theta-\theta^-)T(I)}(k_0(I,  \theta^-)))\right]\\
&+T(I) \int_{-\infty}^{0}\sin(t+s) \left[H_1(\Phi_0^{t+(\theta-\theta^-)T(I)}(z))\right.\\
&\qquad\qquad\qquad\qquad\qquad \left. -H_1(\Phi_0^{t+(\theta-\theta^-)T(I)}(k_0(I,  \theta^-)))\right]dt\\
&+T(I)\int_{0}^{+\infty}\sin(t+s)\left[H_1(\Phi_0^{t+(\theta-\theta^-)T(I)}(z))\right.\\
&\qquad\qquad\qquad\qquad\qquad \left. -H_1(\Phi_0^{t+(\theta-\theta^-)T(I)} (k_0(I,  \theta^+)))\right] dt.\end{split}
\end{equation}
\end{prop}

\begin{proof}
We will use the fact that the perturbation is a separable function, i.e., $G(z,s)=\cos(s)H_1(z)$.
We will also use that \[\Phi^t(k_0(I,\theta))=k_0(R^t(I,\theta))=k_0(I,\theta+t\omega(I)).\]

We perform the change of variable $\tau=t+(\theta+\theta^-)T(I)$ in \eqref{eqn:S0-1}.
We note that, under this change of variable, we have
\begin{equation*}\begin{split}\Phi_0^t(k_0(I, \theta+2\theta^-))=&\Phi_0^{\tau-(\theta+\theta^-)T(I)}(k_0(I, \theta+2\theta^-))\\=&\Phi_0^{\tau}(k_0(I,\theta+2\theta^--(\theta+\theta^-)T(I)\omega(I)))\\
=&\Phi_0^{\tau}(k_0(I,  \theta^-))
\end{split}\end{equation*}
since $T(I)\omega(I)=1$.
Similarly,
\begin{equation*}\begin{split}\Phi_0^t(k_0(I, \theta))=&\Phi_0^{\tau}(k_0(I,\theta -(\theta+\theta^-)T(I)\omega(I)))\\=&\Phi_0^{\tau}(k_0(I,-\theta^-))=\Phi_0^{\tau}(k_0(I, \theta^+)).
\end{split}\end{equation*}
The last equality is by the $S$-symmetry.

We obtain
\begin{equation}\label{eqn:proof1}
\begin{split}
S(I,\theta,s)=&\int_{-\infty}^{(\theta+\theta^-)T(I)}\cos(s+\tau-(\theta+\theta^-)T(I))
\left[H_1(\Phi_0^{\tau}(z))\right.\\
&\qquad\qquad\qquad\qquad\qquad\qquad \left. -H_1(\Phi_0^{\tau}(k_0(I,  \theta^-)))\right]d\tau\\
+&\int_{(\theta+\theta^-)T(I)}^{+\infty}\cos(s+\tau-(\theta+\theta^-)T(I))\left[H_1(\Phi_0^{\tau}(z))\right.\\
&\qquad\qquad\qquad\qquad\qquad\qquad\qquad \left. -H_1(\Phi_0^{\tau}(k_0 (I, \theta^+)))\right] d\tau.
\end{split}
\end{equation}

We recall the following property
\begin{equation}\label{eqn:formula}\begin{split}
\frac{d}{d\theta}\int_{a(\theta)}^{b(\theta)}f(x,\theta)dx&=
\int_{a(\theta)}^{b(\theta)}\frac{\partial}{\partial\theta}f(x,\theta)dx\\
&\qquad +
f(b(\theta),\theta)b'(\theta)-f(a(\theta),\theta)a'(\theta).
\end{split}\end{equation}

%

We apply \eqref{eqn:formula} to \eqref{eqn:proof1}, first, for $a(\theta)=-\infty$ and $b(\theta)=(\theta+\theta^-)T(I)$, and second
for  $a(\theta)=(\theta+\theta^-)T(I)$ and $b(\theta)=+\infty$. Since the integrands in \eqref{eqn:proof1} converge  to $0$ exponentially fast as $t\to \pm\infty$,
the terms corresponding to the infinite limits  on the right hand side of \eqref{eqn:formula} vanish.
We obtain:
\begin{equation*}\label{eqn:proof2}
\begin{split}
\frac{d}{d\theta}S(I,\theta,s)=&\quad T(I)\cos(s)\left [ H_1(\Phi_0^{(\theta+\theta^-)T(I)}(z))-H_1(\Phi_0^{(\theta+\theta^-)T(I)}(k_0(I,  \theta^-)))\right]\\
&-T(I)\cos(s)\left [ H_1(\Phi_0^{(\theta+\theta^-)T(I)}(z))-H_1(\Phi_0^{(\theta+\theta^-)T(I)}(k_0(I,  \theta^+)))\right]\\
&+T(I) \int_{-\infty}^{(\theta+\theta^-)T(I)}\sin(s+\tau-(\theta+\theta^-)T(I))
\left[H_1(\Phi_0^{\tau}(z))\right.\\
&\qquad\qquad\qquad\qquad\qquad\qquad \left. -H_1(\Phi_0^{\tau}(k_0(I,  \theta^-)))\right]d\tau\\
&+T(I)\int_{(\theta+\theta^-)T(I)}^{+\infty}\sin(s+\tau-(\theta+\theta^-)T(I))\left[H_1(\Phi_0^{\tau}(z))\right.\\
&\qquad\qquad\qquad\qquad\qquad\qquad\qquad \left. -H_1(\Phi_0^{\tau}(k_0(I, \theta^+)))\right] d\tau.
\end{split}
\end{equation*}

After cancelations and changing  the variable back $t=\tau-(\theta+\theta^-)T(I)$ we obtain
\begin{equation*}\label{eqn:proof3}
\begin{split}
\frac{d}{d\theta}S(I,\theta,s)=&\quad T(I)\cos(s)\left [ H_1(\Phi_0^{(\theta+\theta^-)T(I)}(k_0(I,  \theta^+)))\right.\\
&\qquad\qquad\qquad\left. -H_1(\Phi_0^{(\theta+\theta^-)T(I)}(k_0(I,  \theta^-)))\right]\\
&+T(I) \int_{-\infty}^{0}\sin(s+t) \left[H_1(\Phi_0^{t+(\theta+\theta^-)T(I)}(z))\right.\\
&\qquad\qquad\qquad\qquad\qquad  \left. -H_1(\Phi_0^{t+(\theta+\theta^-)T(I)}(k_0(I, \theta^- )))\right]dt\\
&+T(I)\int_{0}^{+\infty}\sin(s+t)\left[H_1(\Phi_0^{t+(\theta+\theta^-)T(I)}(z))\right.\\
&\qquad\qquad\qquad\qquad\qquad \left. -H_1(\Phi_0^{t+(\theta+\theta^-)T(I)} (k_0(I, \theta^+)))\right] dt.
\end{split}
\end{equation*}

Last, we calculate $ \frac{d}{d\theta}S (s_0(I,\theta),s)=\frac{d}{d\theta}S(I,\theta-2\theta^-,s)$, obtaining

\begin{equation*}\label{eqn:proof4}
\begin{split}
\frac{d}{d\theta}S(s_0(I,\theta),s)=&\quad T(I)\cos(s)\left [ H_1(\Phi_0^{(\theta-\theta^-)T(I)}(k_0(I,  \theta^+)))\right.\\
&\qquad\qquad \qquad \qquad\left.-H_1(\Phi_0^{(\theta-\theta^-)T(I)}(k_0(I,  \theta^-)))\right]\\
&+T(I) \int_{-\infty}^{0}\sin(s+t)
\left[H_1(\Phi_0^{t+(\theta-\theta^-)T(I)}(z))\right.\\
&\qquad\qquad\qquad\qquad\qquad \left. -H_1(\Phi_0^{t+(\theta-\theta^-)T(I)}(k_0(I,  \theta^-)))\right]dt\\
&+T(I)\int_{0}^{+\infty}\sin(s+t)\left[H_1(\Phi_0^{t+(\theta-\theta^-)T(I)}(z))\right.\\
&\qquad\qquad\qquad\qquad\qquad \left. -H_1(\Phi_0^{t+(\theta-\theta^-)T(I)} (k_0(I,  \theta^+)))\right] dt.
\end{split}
\end{equation*}
\end{proof}

\begin{rem} We can obtain in a similar fashion a formula for $\frac{\partial S}{\partial I}$; see \cite{gidea2022melnikov}. However, we do not need it in this paper. \end{rem}

\subsubsection{Reduction to a Poincar\'e section}

We switch from   continuous-time (flow) dynamics to discrete-time dynamics
by taking the  time-$2\pi$ map of the extended flow.
We  choose the Poincar\'e section
\[ \Sigma=\left\{(z,s)\,|\, s=-{\pi}/{2}\right\} .\]

We denote by $f_\eps$ the  time-$2\pi$ map of the extended flow.

Then \[\Lambda_\eps:=\tilde\Lambda_\eps\cap \Sigma\] is a NHIM for $f_\eps$.
Any scattering map for the extended flow  $\tilde\sigma_\eps^{\tilde\Gamma_\eps}$ induces a scattering map on $\Lambda_\eps$ given by
\[\sigma^{\Gamma_\eps}_\eps(z)=\tilde\sigma_\eps^{\tilde\Gamma_\eps}\left(z,-{\pi}/{2}\right)\textrm { for } z\in\Omega^-(\tilde{\Gamma}_\eps)\cap \Lambda_\eps,\]
where $\Gamma_\eps=\tilde\Gamma_\eps\cap \Sigma$.
%

The following is a consequence of Proposition~\ref{prop:S0-2}:

\begin{cor}
\label{prop:S0-3}
Assume  a homoclinic channel as in \eqref{eqn:unperurbed_homoclininc_channel} associated to a symmetric homoclinic point $z$. 
Then the   derivative with respect to the angle variable of the Hamiltonian that generates the scattering map on $\Lambda_\eps$ is given by
\begin{equation*}\label{sac-map-new-form-}
		\begin{split}
			 \frac{d}{d\theta}S(s_0(I,\theta))=&-T(I) \int_{-\infty}^{0}\cos(t) \left[H_1(\Phi_0^{t+(\theta-\theta^-)T(I)}(z)) \right.\\
&\qquad\qquad\qquad\qquad\quad \left. -H_1(\Phi_0^{t+\theta T(I)}(q^*(I)))\right]dt\\
			&-T(I)\int_{0}^{+\infty}\cos(t)\left[H_1(\Phi_0^{t+(\theta-\theta^-)T(I)}(z)) \right. \\
&\qquad\qquad\qquad\qquad\quad  \left. -H_1(\Phi_0^{t+(\theta-2\theta^-)T(I)} (q^*(I)))\right] dt ,\end{split}
\end{equation*}
where  $q^*(I)=\left(x^*(I),0,0,p_y^*(I)\right)$ is the initial point of the periodic orbit of energy level corresponding to $z$.
\end{cor}
\begin{proof}

First, let us get a new expression for the \eqref{eqn:S0_2}. We have that $k_0(I, \theta)=\Phi_0^{\theta T(I)}(x^*(I),0,0,p_y^*(I))$ is a parameterization of the periodic orbit $\lambda(I)$ with initial point $q^*(I)=\left(x^*(I),0,0,p_y^*(I)\right)$. Then $k_0(I, \theta^\pm)=\Phi_0^{\theta^\pm T(I)}(q^*(I))$, and we can rewrite \eqref{eqn:S0_2} as

	\begin{equation*}
		\begin{split}
			\frac{d}{d\theta}S(s_0(I,\theta),s)&= T(I)\cos(s)\left [ H_1(\Phi_0^{(\theta-\theta^-)T(I)+\theta^+ T(I)}(q^*(I)))\right.\\
			&\qquad\qquad\qquad\qquad\left.-H_1(\Phi_0^{(\theta-\theta^-)T(I)+\theta^-T(I)}(q^*(I)))\right]\\
			&+T(I) \int_{-\infty}^{0}\sin(t+s) \left[H_1(\Phi_0^{t+(\theta-\theta^-)T(I)}(z))\right.\\
			&\qquad\qquad\qquad\qquad\qquad \left. -H_1(\Phi_0^{t+(\theta-\theta^-)T(I)+\theta^- T(I)}(q^*(I)))\right]dt\\
			&+T(I)\int_{0}^{+\infty}\sin(t+s)\left[H_1(\Phi_0^{t+(\theta-\theta^-)T(I)}(z))\right.\\
			&\qquad\qquad\qquad\qquad\qquad \left. -H_1(\Phi_0^{t+(\theta-\theta^-)T(I)+\theta^+ T(I)} (q^*(I)))\right] dt.\end{split}
	\end{equation*}
	Since $\theta^+=-\theta^-$, the  derivative with respect to the angle variable of the Hamiltonian that generates the scattering map in the extended phase space can be written as
	\begin{equation*}\label{sac-map-new-form}
		\begin{split}
			\frac{d}{d\theta}S(s_0(I,\theta),s)=& T(I)\cos(s)\left [ H_1(\Phi_0^{(\theta-2\theta^-)T(I)}(q^*(I)))-H_1(\Phi_0^{\theta T(I)}(q^*(I)))\right]\\
			&+T(I) \int_{-\infty}^{0}\sin(t+s) \left[H_1(\Phi_0^{t+(\theta-\theta^-)T(I)}(z)) \right. \\
&\qquad\qquad \left. -H_1(\Phi_0^{t+\theta T(I)}(q^*(I)))\right]dt\\
			&+T(I)\int_{0}^{+\infty}\sin(t+s)\left[H_1(\Phi_0^{t+(\theta-\theta^-)T(I)}(z)) \right. \\
&\qquad\qquad \left.  -H_1(\Phi_0^{t+(\theta-2\theta^-)T(I)} (q^*(I)))\right] dt.\end{split}
	\end{equation*}
	
	For $s=-\pi/2$ we obtain
	\begin{equation}\label{sac-map-new-form}
		\begin{split}
			 \frac{d}{d\theta}S(s_0(I,\theta))=&-T(I) \int_{-\infty}^{0}\cos(t) \left[H_1(\Phi_0^{t+(\theta-\theta^-)T(I)}(z))\right. \\
&\qquad\qquad \left. -H_1(\Phi_0^{t+\theta T(I)}(q^*(I)))\right]dt\\
			&-T(I)\int_{0}^{+\infty}\cos(t)\left[H_1(\Phi_0^{t+(\theta-\theta^-)T(I)}(z))\right. \\
&\qquad\qquad  \left . -H_1(\Phi_0^{t+(\theta-2\theta^-)T(I)} (q^*(I)))\right] dt.\end{split}
	\end{equation}
\end{proof}
	
%

For the numerical computations in Section \ref{sec:numerical_verification} it will be convenient to use expressions for $S$ and for $-\frac{\partial S}{\partial\theta}\circ \sigma_0$ in terms of the (non-symplectic) coordinates $(x^*,\theta)$, which follow immediately from  \eqref{eqn:S0-1} and \eqref{sac-map-new-form}.

\begin{cor} We have
\begin{equation}\label{eqn:S0-x-star}\begin{split}
S  (x^*,\theta)=&\int_{-\infty}^{0}\sin(t) \left[H_1(\Phi_0^{t+(\theta+\theta^-)T(x^*)} (z)))\right.\\
&\qquad\qquad\quad \left.-H_1(\Phi_0^{t+(\theta+ 2\theta^-)T(x^*)}(q^*(x^*)))\right]\,dt\\
&+ \int_{0}^{+\infty}\sin(t) \left[H_1(\Phi_0^{t+(\theta+\theta^-)T(x^*)}(z))\right.\\
&\qquad\qquad\quad\quad\left.-H_1(\Phi_0^{(t+\theta)T(x^*)}(q^*(x^*)))\right ]\,dt.
\end{split}
\end{equation}
and
\begin{equation*}\label{eqn:partial-S0-x-star}
		\begin{split}
			 \frac{d}{d\theta}S(s_0(x^*,\theta))=&-T(x^*) \int_{-\infty}^{0}\cos(t) \left[H_1(\Phi_0^{t+(\theta-\theta^-)T(x^*)}(z))\right. \\
&\qquad\qquad \left. -H_1(\Phi_0^{t+\theta T(x^*)}(q^*(x^*))\right]dt\\
			&-T(x^*)\int_{0}^{+\infty}\cos(t)\left[H_1(\Phi_0^{t+(\theta-\theta^-)T(x^*)}(z))\right. \\
&\qquad\qquad  \left . -H_1(\Phi_0^{t+(\theta-2\theta^-)T(x^*)} (q^*(x^*))\right] dt,\end{split}
\end{equation*}
where  $q^*(x^*)=\left( x^* ,0,0, p_y^*(x^*)\right )$ is the initial condition of a periodic orbit.
\end{cor}

In Section \ref{sec:pertubed_scattering_hom} we will compute   $-\frac{d}{d\theta}S(s_0(x^*,\theta))$ for  scattering maps $\sigma_\eps$ associated to  homoclinic channels of the form
\begin{equation}\label{eqn:unperurbed_homoclininc_channel_1}
  \Gamma_\eps=\bigcup_{\substack{ x^*\in[x^*(h_\beta), x^*(h_\alpha)]\\ \theta\in(\theta^-+\theta^*, \theta^-+1+\theta^*) }}\Phi_\eps^{(\theta-\theta^-)T} (z),
\end{equation}
for some constant $\theta^*\in\mathbb{R}$. The range of the angle $\theta$ is restricted to an open interval  of length $1$ in order to ensure the monodromy property of the scattering map. See Remark \ref{rem:monodromy}.

	To compute $-\frac{d}{d\theta}S(s_0(x^*,\theta))$ we take a finite partition $\{\theta_j\}_{j=0}^n\subset[\theta^- -1,\theta^-+1]$ with $\theta_j=\frac{1}{n}j$, for $j=-n,\ldots, n$, and for every $\theta_j$ we compute $-\frac{d}{d\theta}S(s_0(x^*,\theta_j))$.
	
Then  we  will restrict $\theta$ to certain intervals of the form $(\theta^-+\theta^*, \theta^-+1+\theta^*)$ corresponding to some particular choices of homoclinic channels and scattering maps. As we shall see, different choices of $\theta^*$ yield different scattering maps with different properties.

\subsection{Mechanism of diffusion based on one scattering map and the Birkhoff Ergodic Theorem}
\label{sec:mechanism}

We describe a general mechanism of diffusion based on iterating a single scattering map.
To avoid notation overload, we will identify $\Lambda_\eps$ through $k_\eps$ with a neighborhood of  $\mathbb{A}:=[I_\alpha,I_\beta]\times \mathbb{T}^1$
(by slightly changing $I_\alpha,I_\beta$ in \eqref{eqn:reference_manifold}).
We will also identify the scattering map $\sigma_\eps$  with its action-angle coordinate representation $s_\eps$.
The main result of this section is  Theorem \ref{prop:mechanism}.
We start by formulating some general assumptions for this result.

We denote by $\mu=dI\wedge d\theta$ the Liouville  measure on $\mathbb{A}$.

We recall that the scattering map is in general  not globally defined. As we shall see in  Section \ref{sec:pertubed_scattering_hom},
in the EH4BP, if we extend the scattering map  the whole annulus it will fail monodromy, since moving $z^-$  a full turn around a non-trivial circle in the annulus $\mathbb{A}$ yields different values for $z^+=\sigma_\eps(z^-)$ (see Remark \ref{rem:monodromy}).
To avoid the lack of monodromy, and hence ensure that $\sigma_\eps$ is well defined, we have to restrict the domain of the scattering map to $\mathbb{A}\setminus \mathcal{C}$,  where $\mathcal{C}$ is some curve that  crosses all non-trivial circle in  $\mathbb{A}$,  e.g.,
$\mathcal{C}=\{\theta=\theta_0\}$ for some $\theta_0$.
We now makes this more precise.

Assume that we have a  scattering map
$\sigma_\eps:\Omega^-(\Gamma)\to \Omega^+(\Gamma)$, such that
\begin{equation}\label{eqn:domains}
\Omega^-(\Gamma)\cap \mathbb{A}=\mathbb{A}\setminus \mathcal{C},
\end{equation}
where
$\mathcal{C}$ is a $C^1$-curve  that runs from the lower boundary of the annulus to the upper boundary, i.e., the graph of a     $C^1$-function in  $I\in [I_\alpha, I_\beta]\mapsto\mathcal{C} (I)\in\mathbb{A}$ with
 $\mathcal \mathcal{C}(I_\alpha)\in \{I_\alpha\}\times\mathbb{T}^1$ and $\mathcal \mathcal{C}(I_\beta)\in \{I_\beta\}\times\mathbb{T}^1$ (with an abuse of notation, $\mathcal{C}$ stands for both the function and its graph).
By this assumption, the scattering map $\sigma_\eps$ is  defined on all $\mathbb{A}$ except for the curve $\mathcal{C}$, which is a set of $\mu$-measure zero.

We assume that the scattering map $\sigma_\eps$ expressed in the action-angle coordinates  is of the form
\begin{equation}
\label{eqn:many_scatterings}
 \begin{split}
  \sigma_\eps(I,\theta)=
   &(I,\theta+\Delta(I))+\eps\left(-\frac{\partial S}{\partial \theta}, \frac{\partial S}{\partial I} \right)(I,\theta+\Delta(I))+O(\eps^2),\\
 \end{split}
\end{equation}
for $(I,\theta)\in \mathbb{A}\setminus\mathcal{C}$.

We  also assume  the  following non-degeneracy condition:
\begin{equation}\label{eqn:sum_integrals}
   \int_{\mathbb{A}} \left[ -\frac{\partial S }{\partial \theta}(I,\theta+\Delta (I))\right]\, d\mu \neq 0 .
\end{equation}
Since $S$ is defined a.e. on $\mathbb{A}$, the above integral is well defined in the sense of Lebesgue.

We will see in Section \ref{sec:verification_mechanism}
that these conditions are met in the EH4BP.

Let
\[\sigma_0(I,\theta)=(I,\theta+\Delta(I)).\]

As a consequence of the smooth dependence of the scattering map $\sigma_\eps$ on the variables and on the parameter $\eps$,  the function $S$ in \eqref{eqn:many_scatterings} is $C^1$-bounded,
i.e., for some $M>0$ we have
\begin{equation}
\label{eqn:bounds_scatterings}
 \| S \|_{C^1}\le M ,
\end{equation}
and  the error term  in \eqref{eqn:many_scatterings} is uniformly bounded, that is, there exists $C_0>0$ such that for all $(I,\theta)$ we have
\begin{equation}\label{eqn:Oeps2}
\left \vert
  \sigma_\eps(I,\theta)-\left(
    (I,\theta+\Delta(I))+\eps\left(-\frac{\partial S}{\partial \theta}, \frac{\partial S}{\partial I} \right)(I,\theta+\Delta(I))\right )\right\vert <C_0\eps^2.
\end{equation}

\begin{teo}\label{prop:mechanism}
Assume the conditions \eqref{eqn:domains}, \eqref{eqn:many_scatterings}, 
and \eqref{eqn:sum_integrals}.

Then there exists $\eps_0>0$ and $C>0$, such that for each $0<\eps<\eps_0$ there exists a (pseudo-) orbit $(I_n,\theta_n)$ of $\sigma_\eps$, for $n=0,\ldots,N$,   such that
\[ \left\vert I_N- I_0 \right\vert >C>0.\]
\end{teo}

\begin{proof}
Suppose that  the  integral  in \eqref{eqn:sum_integrals} is positive. (A similar argument can be done when the integral is negative.)
There exists $C_1>0$ such that:
\begin{equation}\label{eqn:sum_integrals_C1}
  \int_{\mathbb{A}}  \left[ -\frac{\partial S }{\partial \theta}(I,\theta+\Delta (I)) \right] \, d\mu>C_1>0.
\end{equation}

Define $\mathbb{A}'= [I'_\alpha,I'_\beta]\times\mathbb{T}^1$ for some $I'_\alpha \gtrsim  I_\alpha$ and $I'_\beta\lesssim I_\beta$.
That is, the annulus $\mathbb{A}'$ is inside and slightly smaller than the  annulus $\mathbb{A}$.
Choose $C_2$ slightly smaller than $C_1$, i.e., $0<C_2\lesssim C_1$ and impose that  $I'_\alpha$ and $I'_\beta$ are such that

\begin{equation}\label{eqn:condition_M}
 M\cdot\mu(\mathbb{A}\setminus\mathbb{A}') < C_1-C_2 .
\end{equation}

The bounds \eqref{eqn:bounds_scatterings} and \eqref{eqn:condition_M}  imply that
\begin{equation*}\label{eqn:condition_A_prime}
\int_{\mathbb{A}\setminus\mathbb{A}'} \left\vert -\frac{\partial S }{\partial \theta}\circ \sigma_0(I,\theta) \right\vert \, d\mu<C_1-C_2.
\end{equation*}

Choose $C_4$ such that
\begin{equation}\label{eqn:C}
  0<C_4\lesssim C_3:=\frac{C_2}{\mu(\mathbb{A})},
\end{equation}
 and let
\begin{equation}\label{eqn:eps0}
  \eps_0=\frac{C_3-C_4}{1+C_0},
\end{equation}
where $C_0$ is the constant from \eqref{eqn:Oeps2}. Note that  the constants $C_1$, $C_2$, $C_3$, $C_4$ do not dependent on $\eps$.

We denote an orbit of $\sigma_\eps$ by
\[(I_n,\theta_n), \textrm{ for }n=0,\ldots,N.\]
To distinguish it from a `true orbit' of the system, we will refer to this as a  (pseudo-)orbit.

We will argue that there exist  (pseudo-)orbits of $\sigma_\eps$ along which the action $I$ increases by a constant $C$ independent of $\eps$.

Let  $\eps>0$ be fixed in $(0,\eps_0)$, where $\eps_0$ was defined in \eqref{eqn:eps0}.
We consider the sets
\[X_\eps^N=\bigcup_{0\le n\le N}  \sigma_\eps^n (\mathbb{A}'),\]
for $N\ge 0$.

We consider two cases.

\emph{Case 1.} Assume that  for some $N\ge 0$ we have that $X_\eps^N\not\subseteq \mathbb{A}$.  Then there are points $(I_0,\theta_0)\in \mathbb{A}'$ whose orbits leave
$\mathbb{A}$ at some moment. In this case, there exist $N$ such that
\[  I_N>I_\beta, \textrm{ or }   I_N<I_\alpha.\]
Since $I_\alpha'\le I_0\le I_\beta'$, this  pseudo-orbit  changes the action coordinate by at least
\begin{equation}\label{eqn:C5}
 C_5=\min\{|I_\alpha'-I_\alpha|, |I_\beta'-I_\beta|\},
\end{equation}
which is independent of $\eps$. In this case, we obtain the conclusion of  Theorem~\ref{prop:mechanism}.

\emph{Case 2.}   If the assumption in Case 1 does not hold, then the set
\[X_\eps =\bigcup_{n\ge 0}  \sigma_\eps^n (\mathbb{A}')\]
is forward invariant under $\sigma_\eps$ and
\begin{equation}
\label{eqn:X_inclusion}
\mathbb{A}'\subseteq X_\eps\subseteq \mathbb{A},
\end{equation}
hence of finite measure. See Fig.~\ref{fig:strip}.

\begin{figure}
\centering
\includegraphics[scale=0.35]{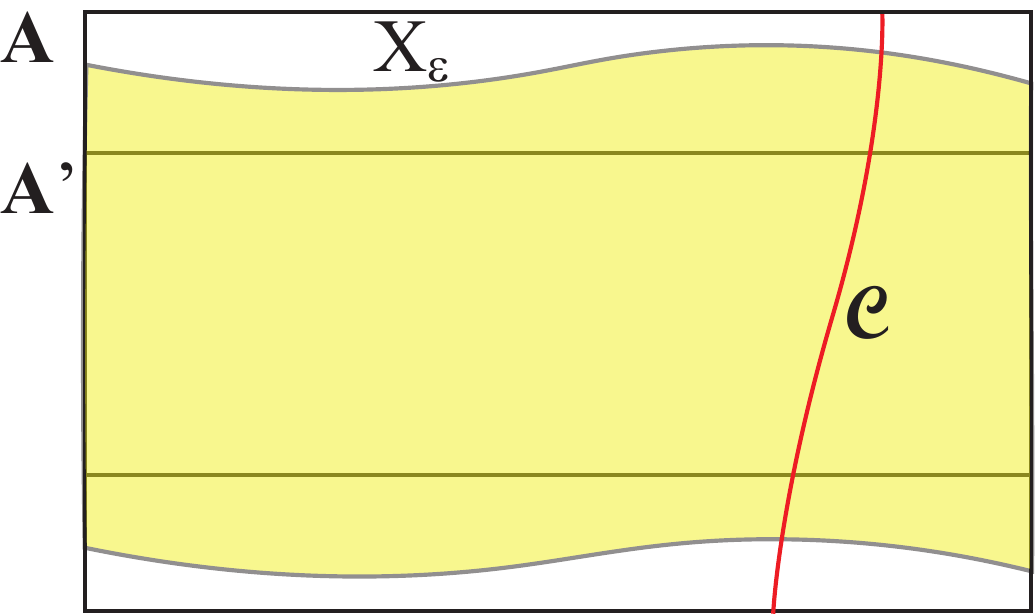}
	\caption{The annuli $\mathbb{A}$ and $\mathbb{A}'$, the curve $\mathcal{C}$, and the invariant set $X_\eps$.}\label{fig:strip}
\end{figure}

Since $\sigma_\eps(X_\eps)\subseteq X_\eps$,
an  orbit $(I_n,\theta_n)_{n\ge 0}$ starting at a.e. point $(I_0,\theta_0)$ in $X_\eps$ stays in $X_\eps$ for all future times.
In this second case, we will obtain diffusing orbits using the argument below.

Recall  Birkhoff's Ergodic Theorem. If $T:X\to X$ is a  measure preserving map, where $(X,\mathcal{B},\mu )$  is a finite measure space, and $f\in L^1(X)$, then
\[\frac{1}{N}\sum_{n=0}^{N-1} f\circ T^n (x) \to \bar f(x)  \textrm{ for }  \mu-\textrm{a.e. } x\in X\]
where $\bar f \in L^1(X)$ is $T$-invariant and satisfies
\[\int_X \bar f d\mu =\int_X f d\mu.\]
Moreover, if $T$ is ergodic (i.e., any $T$-invariant set has either zero or full measure), then  the pointwise limit is a constant
$\bar f=\frac{1}{\mu(X)}\int_X f d\mu$.

For each $\eps\in(0,\eps_0)$, with $\eps_0$ given by \eqref{eqn:eps0}, we will apply Birkhoff's  theorem  in the case when $X:=X_\eps$, $\mu=dI\wedge d\theta$ is the Liouville  measure on $\mathbb{A}$ restricted to $X_\eps$,  $T=T_\eps:=\sigma_\eps$, and $f=-\frac{\partial{S}}{\partial \theta}\circ\sigma _0$.
By \eqref{eqn:X_inclusion}, $X_\eps$ is of finite measure, and by construction it is forward invariant under $\sigma_\eps$.
Since the scattering map $\sigma_\eps$  is symplectic, it is $\mu$-measure preserving.
By \eqref{eqn:domains} and \eqref{eqn:bounds_scatterings}, the function $f$  is defined a.e. on $\mathbb{A}$ and is an $L^1$-function.

By \eqref{eqn:X_inclusion} we have
\[\int_{X_\eps} f\,d\mu=\int_{\mathbb{A}} f\,d\mu- \int_{(\mathbb{A}\setminus \mathbb{A}')\setminus X_\eps} f\,d\mu.\]

Since by  \eqref{eqn:bounds_scatterings}, \eqref{eqn:X_inclusion}, and \eqref{eqn:condition_M} we have
\[ \int_{(\mathbb{A}\setminus \mathbb{A}')\setminus X_\eps} \left\vert f\right\vert \,d\mu
\le  \int_{(\mathbb{A}\setminus \mathbb{A}')\setminus X_\eps} M  \,d\mu
\le M\cdot \mu(\mathbb{A}\setminus \mathbb{A}')<C_1-C_2 ,\]
from \eqref{eqn:sum_integrals_C1}
it follows that
\begin{equation}\label{eqn:fC2}
   \int_{X_\eps} f\,d\mu > C_2 .
\end{equation}

Applying Birkhoff's theorem  to $T=\sigma_\eps$ and $f=-\frac{\partial{S}}{\partial \theta}\circ\sigma _0$, for each $\eps\in(0,\eps_0)$  we obtain
\begin{equation}\label{eqn:Birkhoff_eps}
 \frac{1}{N}\sum_{n=0}^{N-1} f \circ T_\eps ^n(I_0,\theta_0) \to \bar{f}_\eps (I_0,\theta_0) \textrm{ for } \mu-\textrm{a.e. } (I_0,\theta_0)\in X_\eps,
\end{equation}
for some $L^1$-function $\bar{f}_\eps$, satisfying
\[ \int_{X_\eps} \bar{f}_\eps \, d \mu=\int_{X_\eps} f\, d \mu.\]
Since $T_\eps$ may not be ergodic, when we apply Birkhoff's theorem  we obtain a pointwise limit that  is in general  not a constant but a function.
Let us denote by $Y_\eps$ the subset of $(I_0,\theta_0)$ in $X_\eps$ for which the  pointwise limit in \eqref{eqn:Birkhoff_eps} holds,
so $\mu(X_\eps-Y_\eps)=0$.

From \eqref{eqn:fC2} it follows
\begin{equation*}\label{eqn:fbar C2}
\int_{X_\eps} \bar{f}_\eps d\mu > C_2>0.
\end{equation*}

This implies that
\begin{equation*}\label{eqn:fbarI0theta0_C2}
 \mu\{(I,\theta)\,|\, \bar{f}_\eps(I,\theta)>{C_2}/{\mu(X_\eps)}\}>0.
\end{equation*}

Since $C_2/\mu(X_\eps)\ge C_2/\mu(\mathbb{A})>0$, then for $C_3=C_2/\mu(\mathbb{A})$, which is  independent of $\eps$, we have that
 there exists a  positive measure set $Z_\eps\subseteq Y_\eps$ of points $(I_0,\theta_0)$  for which
\begin{equation}\label{eqn:fbarI0theta0_C3}
\bar{f}_\eps(I_0,\theta_0)>C_3=\frac{C_2}{\mu(\mathbb{A})}>0.
\end{equation}

For every $\eps\in(0,\eps_0)$ and every $(I_0,\theta_0)\in Z_\eps$, there exists $\bar{N}=\bar{N}(\eps, (I_0,\theta_0))$ such that for $N\ge \bar{N}(\eps, (I_0,\theta_0))$ we have
\[
 \left\vert\frac{1}{N}\sum_{n=0}^{N-1} f \circ T_\eps ^n(I_0,\theta_0) - \bar{f}_\eps(I_0,\theta_0)\right\vert<\eps,
\]
so
\begin{equation}\label{eqn:Birkhoff_ineq}
\sum_{n=0}^{N-1} f\circ T_\eps ^n(I_0,\theta_0)  >N \left(\bar{f}_{\eps}(I_0,\theta_0) -\eps\right).
\end{equation}

For  every point $(I_n,\theta_n)$ of the pseudo-orbit, the corresponding change in action when moving to the next point   $(I_{n+1},\theta_{n+1})$ is
\[
\Delta I_n=I_{n+1}-I_{n}=-\eps\frac{\partial S }{\partial \theta}(I_n,\theta_n+\Delta (I_n))+O(\eps^2).
\]
By \eqref{eqn:Oeps2}, on $\mathbb{A}$ each error term  is bounded by $C_0\eps^2$.
Thus, for $N\ge \bar{N}(\eps, (I_0,\theta_0))$, the cumulative change in $I$ associated to  the points of the pseudo-orbit   is given by
\[\sum_{n=0,\ldots,N-1} \Delta I_n:=\sum_{n=0,\ldots,N-1} \left [-\eps\frac{\partial S }{\partial \theta}(I_n,\theta_n+\Delta (I_n))+O(\eps^2)\right].\]
Using  \eqref{eqn:fbarI0theta0_C3},  \eqref{eqn:Birkhoff_ineq}, and considering the bounds on the error terms \eqref{eqn:Oeps2}  we have
\[ \sum_{n=0,\ldots,N-1} \Delta I_n>\eps \sum_{n=0}^{N-1} f\circ T^n_\eps(I_0,\theta_0) - N C_0\eps^2>\eps N \left(\bar{f}_{\eps}(I_0,\theta_0) -(1+C_0)\eps\right).\]

By \eqref{eqn:fbarI0theta0_C3}, for points $(I_0,\theta_0)\in Z_\eps$ and $N>\bar{N}(\eps,(I_0,\theta_0))$ we have
\[\sum_{n=0,\ldots,N-1} \Delta I_n >\eps N (C_3-(1+C_0)\eps)>\eps N C_4,\]
for  $0<C_4 \lesssim C_3$ defined in \eqref{eqn:C}, which is independent of $\eps$.
The last inequality holds for each $\eps\in(0,\eps_0)$, where $\eps_0$ has been defined in \eqref{eqn:eps0}.
Since this holds for $N\ge N(\eps,(I_0,\theta_0))$ large enough, we can choose $N>1/\eps$, and we conclude that
\[ \sum_{n=0,\ldots,N-1} \Delta I_n  >C_4>0.\]

That is, there exists a positive measure set $Z_\eps$  of points $(I_0,\theta_0)$ in $X_\eps$ such that the change in $I$ along the corresponding pseudo-orbit $(I_n,\theta_n)_{n=0,\ldots,N}$ is at least $C_4>0$.

Putting together \emph{Case 1} and \emph{Case 2}, there exists  $C=\min\{C_4,C_5\}$ (where $C_5$ comes from \eqref{eqn:C5}), such that  $\left\vert I_N- I_0 \right\vert >C>0$. This concludes the proof.
\end{proof}

\begin{rem}\label{rem:birkhoff_rate}
It is well known that for general ergodic dynamical systems, the rate of convergence of the Birkhoff   sums can be arbitrarily slow,
as shown in \cite{krengel1978speed}.
Therefore, in the above argument we may need to choose  $N$ very large.
This approach does not immediately lead to estimates for the diffusion time.
\end{rem}

\begin{rem}
Theorem~\ref{prop:mechanism} only says that there is a pseudo-orbit along which $I$ changes by $O(1)$, but it does not say whether $I$ increases  or decreases along that pseudo-orbit. If we have more information on  the scattering maps, we can also obtain the existence of pseudo-orbits along which $I$ increases  (decreases) by $O(1)$. See Theorem \ref{th:mechanism-main}.
\end{rem}

\subsection{Mechanism of diffusion based on multiple scattering maps}
\label{sec:mechanism_two}

In this section we recall another mechanism of diffusion based on several scattering maps.

\begin{teo}{\cite[Theorem 10]{CapinskiGL17}}
\label{th:mechanism-main} Assume that we have
a family  of scattering maps $\sigma_{\eps}^{1},\ldots,\sigma_{\eps}^{L}$ on $\mathbb{A}=[I_\alpha,I_\beta]\times \mathbb{T}^1$, which are  of the form
\eqref{eqn:many_scatterings},  with corresponding generating Hamiltonian functions
$S^{1},\ldots,S^{L}$.

Assume that there exists $c>0$ such that for every
$\left( I,\theta \right) \in \mathbb{A}$ there exists $j\in\{1,\ldots,L\}$ such that
\begin{equation}
-\frac{\partial S^{j}}{\partial \theta }\circ \sigma_0^j(I,\theta)>c.
\label{eq:scatter-grad-a1}
\end{equation}

Then there exists $\eps_0>0$ and $I_\alpha<I_a<I_b<I_\beta$, such that for each $0<\eps<\eps_0$ there exists a (pseudo-) orbit $(I_n,\theta_n)$ of the iterated function system $\{\sigma^j_\eps\}_{j=1,\ldots, L}$, for $n=0,\ldots,N$,   such that
\[  I_0<I_a \textrm{ and } I_N>I_b.\]

Similarly, if there exists $c>0$ such that for every $\left( I,\theta \right) \in
\mathbb{A}$ there exists $j\in\{1,\ldots,L\}$ such that
\begin{equation*}
-\frac{\partial S^{j}}{\partial \theta }\circ \ \sigma_0^j(I,\theta)<-c<0,
\label{eq:scatter-grad-a2}
\end{equation*}%
then there exists $\eps_0,I_a,I_b$ as before, such that for each $0<\eps<\eps_0$ there exists a (pseudo-) orbit $(I_n,\theta_n)$, $n=0,\ldots, N$, such that
\[  I_0>I_b \textrm{ and } I_N<I_a.\]
\end{teo}

\subsection{Proof of Theorem \ref{teo:main}}

Applying Theorem \ref{prop:mechanism} (resp., Theorem \ref{th:mechanism-main}) for each $\eps<\eps_0$ there exists  a pseudo-orbit of the scattering map $\sigma_\eps$ (resp., of the iterated function system $\sigma^1_\eps,\ldots, \sigma^L_\eps$) along which the action $I$ changes by some amount independent of $\eps$.

We now recall  \cite[Theorem 3.7]{GideaLlaveSeara20-CPAM}, which says that any pseudo-orbit   of an iterated system  of scattering maps can be shadowed by a true orbit, provided that the there is a neighborhood $\mathcal{U}$ of the  pseudo-orbit such that almost every point
of $\mathcal{U}$ is recurrent under the inner map (the restriction of  $f_\eps$ to the NHIM). Thus, if there is  such a neighborhood $\mathcal{U}$, this result implies the existence of a true orbits shadowing the pseudo-orbit, thus achieving a change in the action $I$  by some constant  independent of $\eps$.

Otherwise,  there are orbits of the inner map that leave every such  neighborhood $\mathcal{U}$, thus traveling $O(1)$.  Therefore we obtain diffusion by the inner map  alone.

This  dichotomy  argument has  been employed in \cite{GideaLlaveSeara20-CPAM}.

\begin{rem}\label{rem:KAM}
If in Theorem \ref{th:mechanism-main} we assume \eqref{eq:scatter-grad-a1}, which implies the existence of pseudo-orbits of the scattering map along which $I$ increases by $O(1)$,
the dichotomy part of the above proof does not allow us to conclude the existence of true orbits along which $I$ increases by $O(1)$ (but only of true orbits along which $I$ changes by $O(1)$).

However, if we further assume that the inner map $f_\eps$ satisfies the conditions of the KAM Theorem, then it follows that the NHIM
contains an annulus   bounded by KAM tori, which is invariant under $f_\eps$ (as in \cite{CapinskiGL17}).  Then the Poincar\'e Recurrence Theorem applies to $f_\eps$ on that annulus.
In this case, we do not need to invoke the dichotomy argument, and we can conclude the existence of true orbits along which $I$ increases by $O(1)$.
\end{rem}

\section{Numerical verification of the conditions in Theorem \ref{teo:main}}
\label{sec:numerical_verification}
In Sections \ref{sec:transverse_connections}, \ref{sec:unperturbed_scattering_hom_het}, \ref{sec:pertubed_scattering_hom_het} we verify numerically the conditions \textrm{(i)-(v)} of Theorem \ref{teo:main}. That is, we establish numerically the existence of transverse homoclinic and heteroclinic connections to the NHIM  $\Lambda^1_0$ and of heteroclinic connections between  $\Lambda^1_0$ and $\Lambda^2_0$ for the  unperturbed system, we compute the  unperturbed scattering map associated to  the homoclinic/heteroclinic connections, and then we compute the effect of the perturbation on these maps.
In Section \ref{sec:verification_mechanism}  we verify numerically the conditions of  Theorem  \ref{prop:mechanism},
and in Section \ref{sec:mechanism-main} we verify numerically the conditions of Theorem \ref{th:mechanism-main}.
These imply  the condition \textrm{(vi)} of Theorem \ref{teo:main}. %
\subsection{Homoclinic and heteroclinic connections to the unperturbed NHIMs}
\label{sec:transverse_connections}

\subsubsection{Lyapunov orbits}
\label{sec:Lyapunov}
\begin{figure}
\centering
   \includegraphics[scale=0.35]{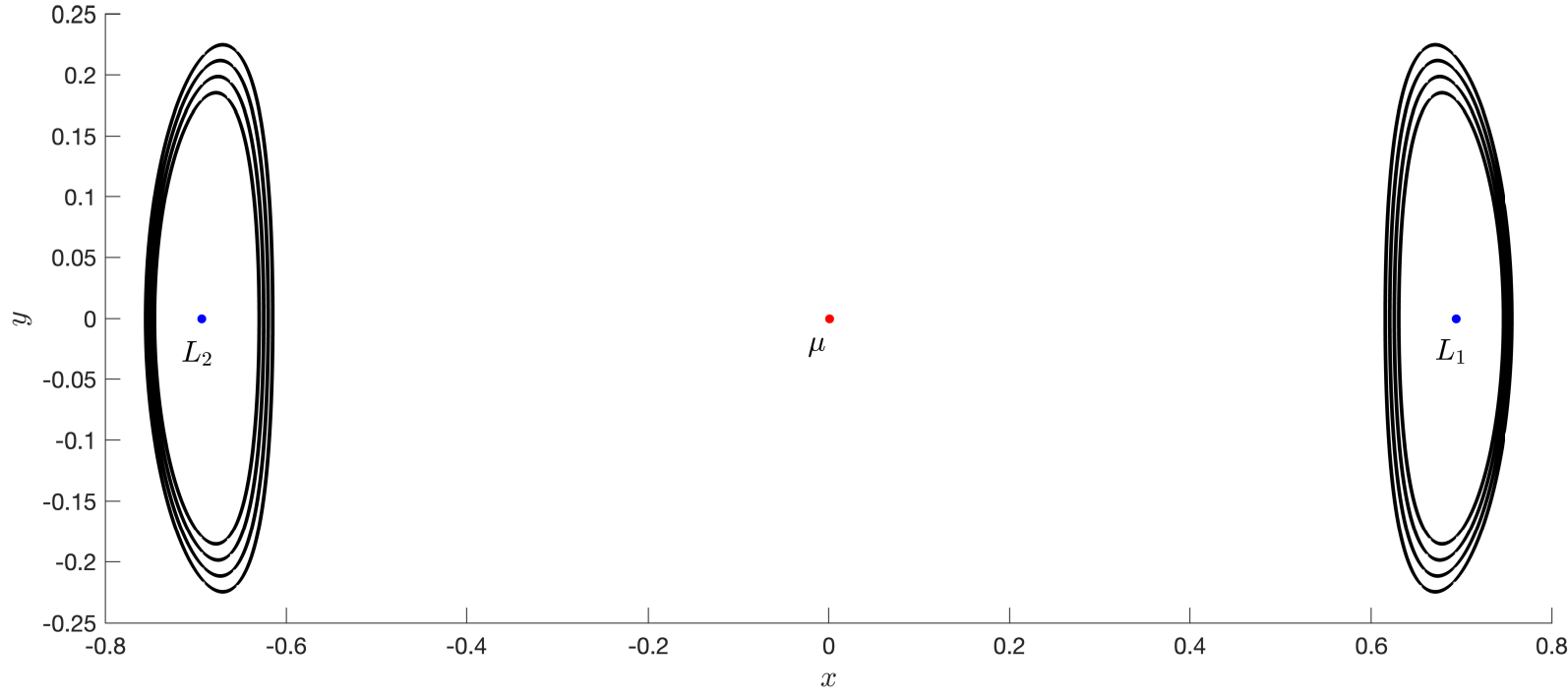}
	\caption{Families of periodic Lyapunov orbits emanating from $L_1$ and $L_2$ for $\pm x^*(h) \in\text{[0.615,0.63]} $.}\label{fig:families_periodic_orbits}
\end{figure}

We compute   the NHIMs $\Lambda^1_0$ near $L_1$, for  the energy range
\begin{equation}\label{eqn:energy_range_hom}
 [h_\alpha,h_\beta]=[ -2.10446079 ,-2.07715457  ]
\end{equation}
The corresponding range for $x^*(h)$ is $ [x^*(h_\beta), x^*(h_\alpha)]=[0.615,0.63 ]$.
In Table~\ref{Cuadro-initial-conditions-1}, we give the energy levels, the initial conditions and the periods for a family of  periodic orbits with $x^*(h)$ within this range.
The corresponding periodic orbits are shown in Fig.~\ref{fig:families_periodic_orbits}.

\subsubsection{Homoclinic connections}
\label{sec:homoclinic connections}
\begin{figure}
		\centering
        \subfigure[Projection on the $(x,y,p_x)$-coordinates.]{\includegraphics[scale=0.27]{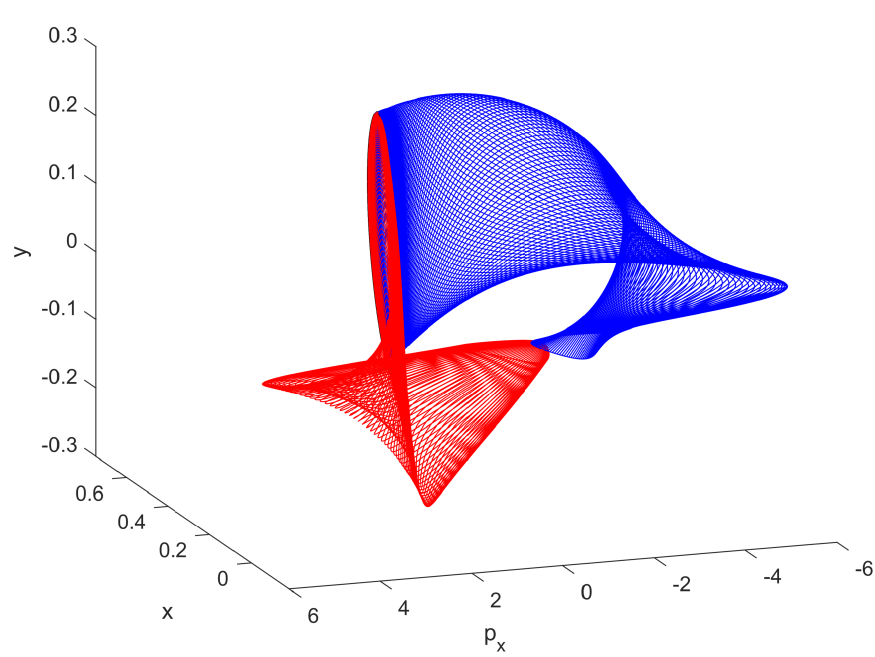}}
		\subfigure[Projection on the $(x,y)$-coordinates.]{\includegraphics[scale=0.26]{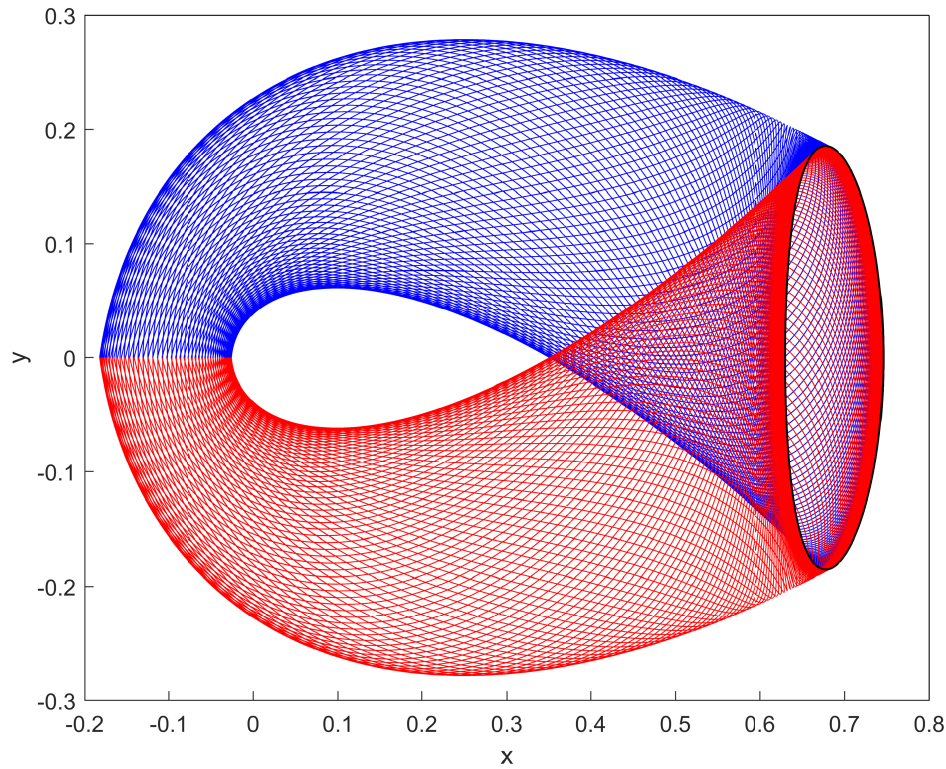}}
		\subfigure[Projection on the $(x,p_x)$-coordinates.]{\includegraphics[scale=0.26]{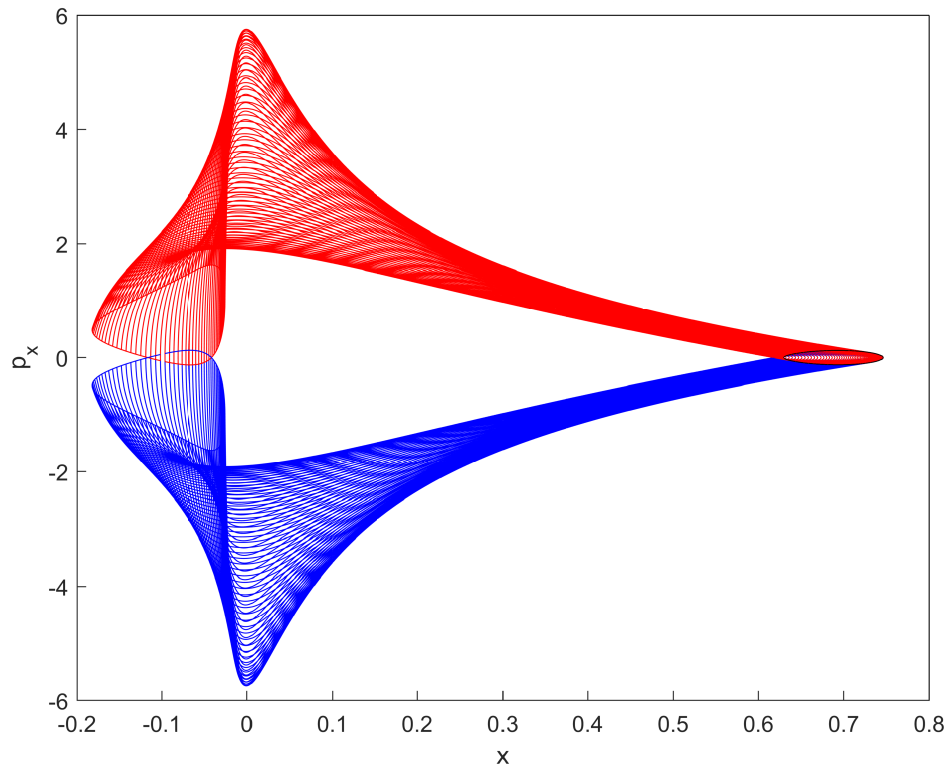}}
		\caption{Stable (red) and unstable (blue) manifold of the Lyapunov periodic orbit for $x^*(h_\alpha)=0.63$. \label{manifolds_homoclinics}}
\end{figure}

We compute the stable and unstable manifolds, $W^s(\lambda_1(h))$ and $W^u(\lambda_1(h))$, respectively,  associated to the Lyapunov orbits $\lambda_1(h)$, for $h\in[h_\alpha,h_\beta]$. As an example, we show the projections of the stable and unstable manifolds of the Lyapunov orbit corresponding to $x^*(h_\alpha)=0.63$,  onto the $(x,y,p_x)$-coordinates,  $(x,y)$-coordinates, and $(x,p_x)$ coordinates in Fig.~\ref{manifolds_homoclinics}.

\begin{figure}
	\centering
    \subfigure[Homoclinic points $z_2$ (left) and  $z_1$ (right).]{\includegraphics[scale=0.25]{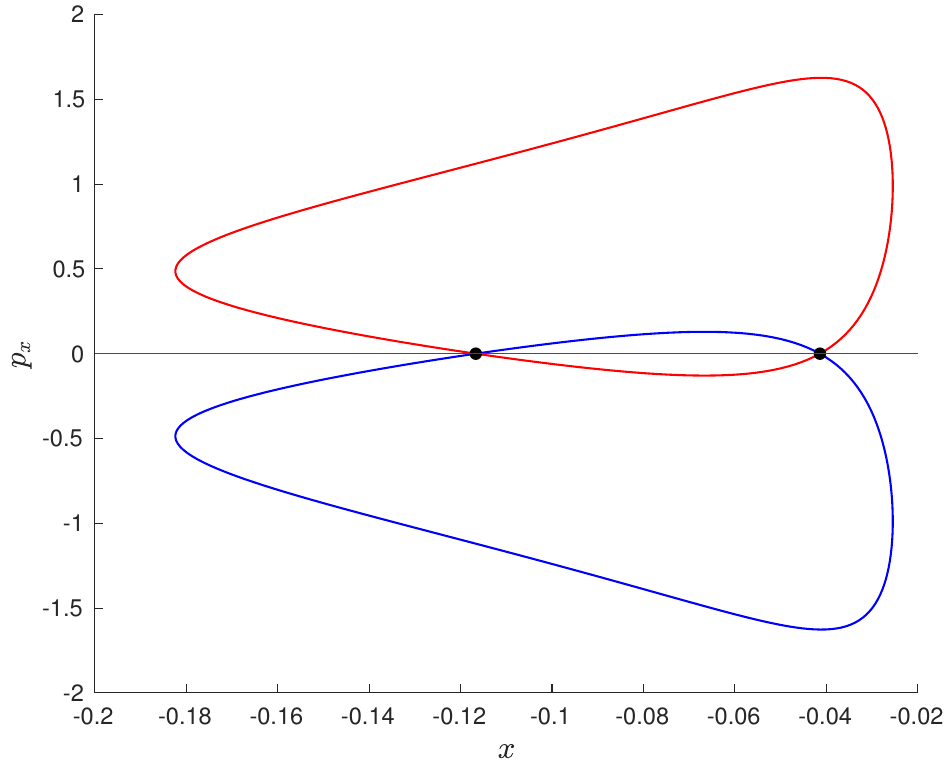}}
    \subfigure[Homoclinic orbit $\Phi^t_0(z_2)$.]{\includegraphics[scale=0.25]{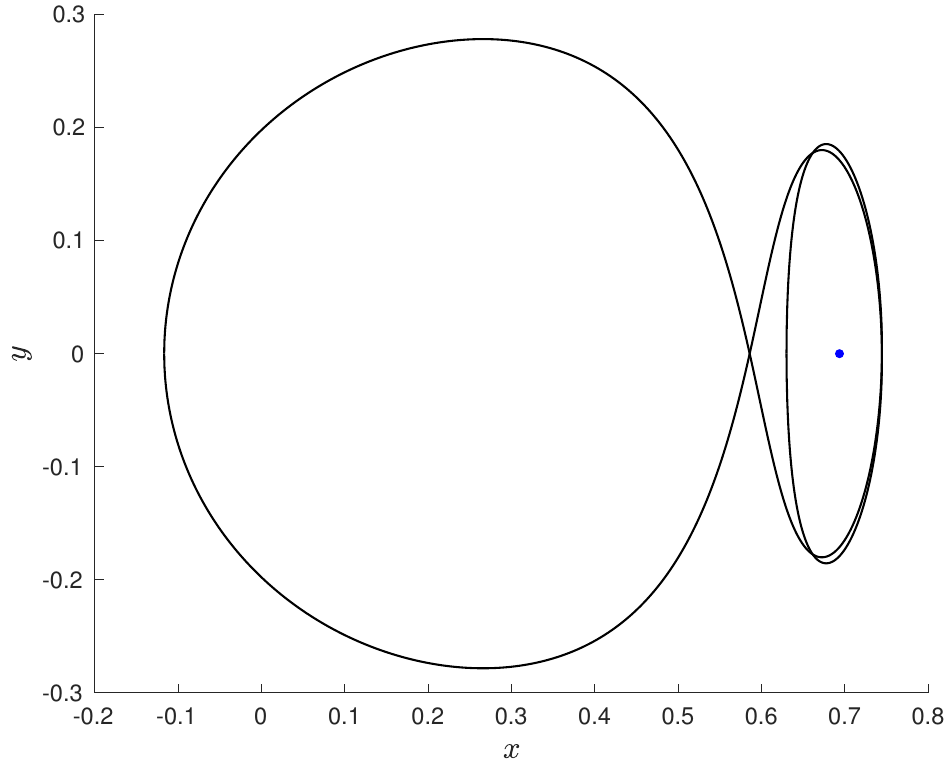}}
    \subfigure[Homoclinic orbit $\Phi^t_0(z_1)$.]{\includegraphics[scale=0.25]{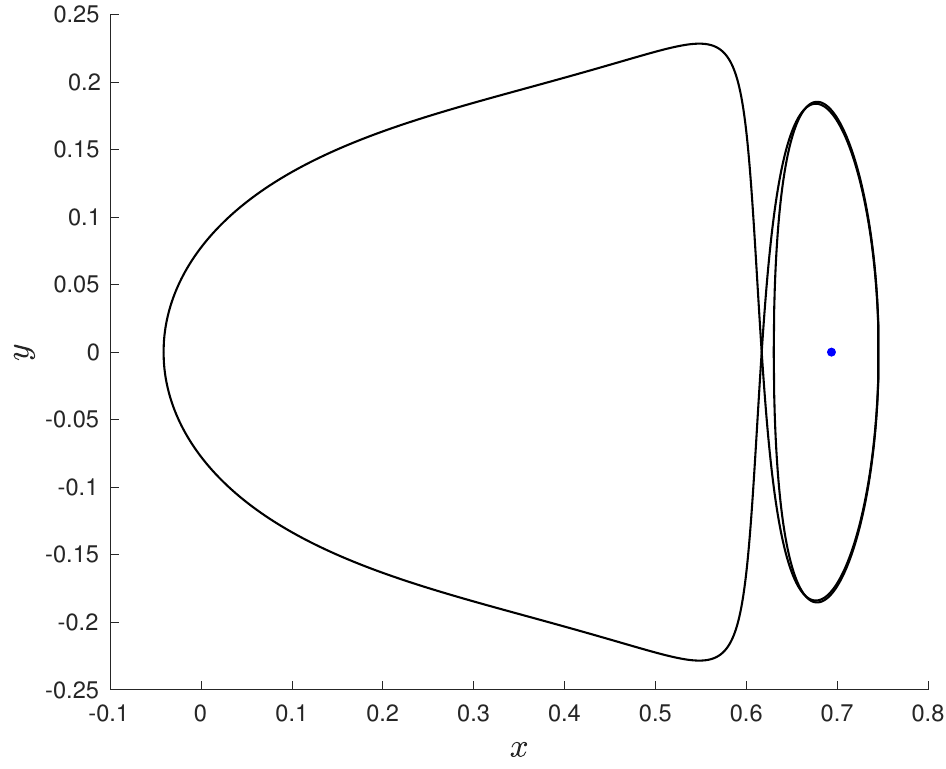}}
	\caption{Homoclinic points and orbits. 
\label{homoclinics orbits}}
\end{figure}

The manifolds are computed up to the their first intersections with the section $\mathscr{S}_x=\{y=0, p_y<0\}$. On this section, we can find transverse homoclinic points as  intersection between the stable and unstable manifold cuts with the section $\mathscr{S}_x$.
For example, the intersections of the stable and unstable manifold associated with the periodic orbit corresponding to $x^*(h_\alpha)=0.63$ are shown  in Fig.~\ref{homoclinics orbits}. Note that there are two symmetric homoclinic points, $z_1$ and $z_2$ (given explicitly in Table \ref{Cuadro-homoclinic-points-z1} and Table  \ref{Cuadro-homoclinic-points-z2} for $x(h_\alpha)=0.63$),  which give rise to two homoclinic orbits that are shown in Fig.~\ref{homoclinics orbits}.

\subsubsection{Heteroclinic connections}
\label{sec:heteroclinic_connections}
For studying the heteroclinic connections between $\Lambda^1_0$ and $\Lambda^2_0$, we consider the same energy range $[h_\alpha,h_\beta]$ given in \eqref{eqn:energy_range_hom}.

The corresponding range for $x^*(h)$ is $ [x^*_1(h_\beta), x^*_1(h_\alpha)]=[0.615, 0.63]$ and $ [x^*_2(h_\alpha), x^*_2(h_\beta)]=[-0.63,-0.615]$, with the positive sign corresponding to $L_1$ and the negative sign to $L_2$.
%
%

\begin{figure}
	\centering
	\subfigure[Projection on the $(x,y)$ plane.]{\includegraphics[scale=0.28]{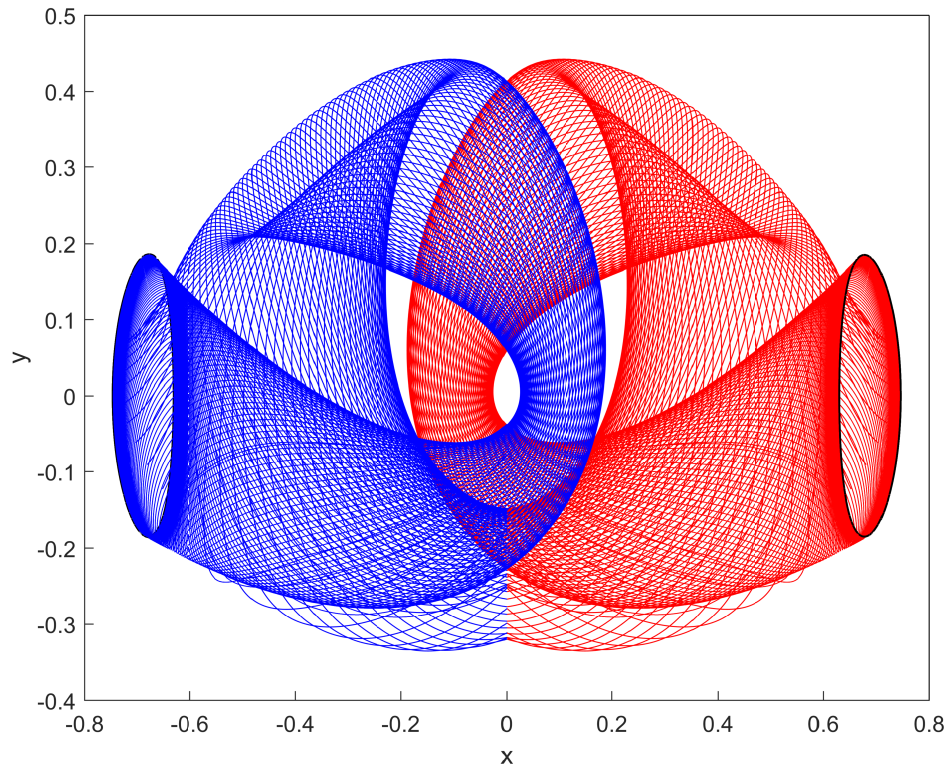}}\hspace{0.001cm}
	\subfigure[Projection on the $(x,p_y)$ plane.]{\includegraphics[scale=0.28]{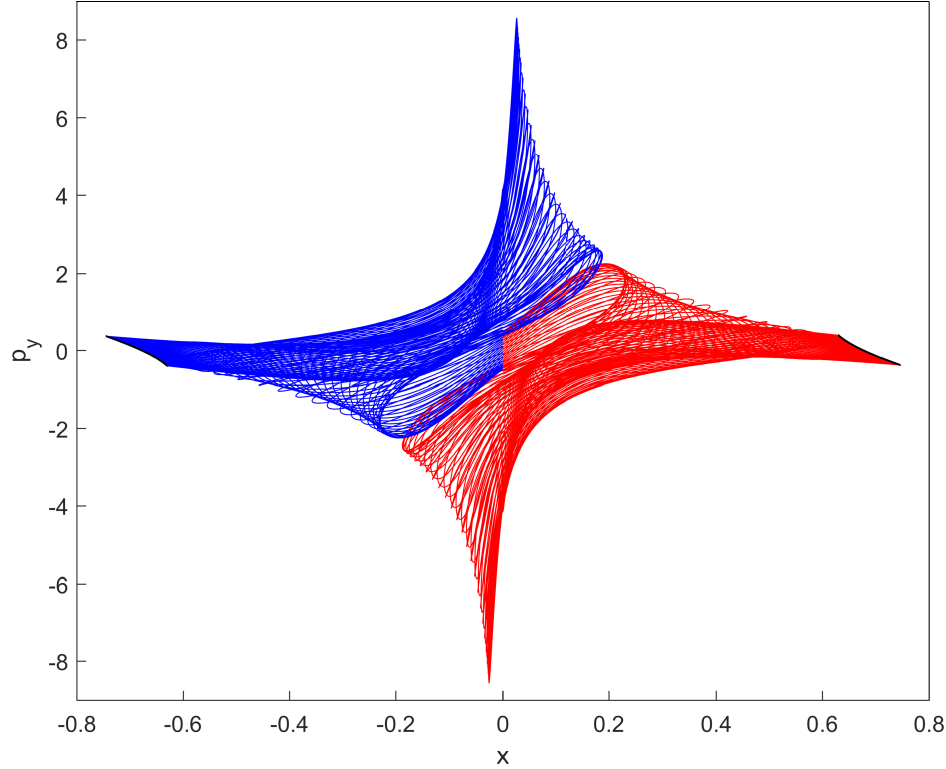}}
	\caption{Stable manifold (red) of the Lyapunov periodic orbit that start at $x_1^*(h_\alpha)=0.63$ and unstable manifold (blue) of the Lyapunov periodic orbit that start at $x_2^*(h_\alpha)=-0.63$. \label{fig:intersection heteroclinic orbit}}
\end{figure}

We compute the stable and unstable manifolds, $W^s(\lambda_1(h))$ and $W^u(\lambda_2(h))$  associated to Lyapunov orbits $\lambda_1(h)$ and $\lambda_2(h)$, respectively, for $h\in[h_\alpha,h_\beta]$. As an example, we show the projection of the stable manifold of the Lyapunov periodic orbit that start at $x_1^*(h_\alpha)=0.63$ and unstable manifold of the Lyapunov periodic orbit that start at $x_2^*(h_\alpha)=-0.63$,  onto the $(x,y)$-coordinates, and $(x,p_x)$ coordinates in Fig.~\ref{fig:intersection heteroclinic orbit}.

\begin{figure}
	\centering	
    \subfigure[Heteroclinic points $z_2$ (on the left) and  $z_1$ (on the right).]{\includegraphics[scale=0.25]{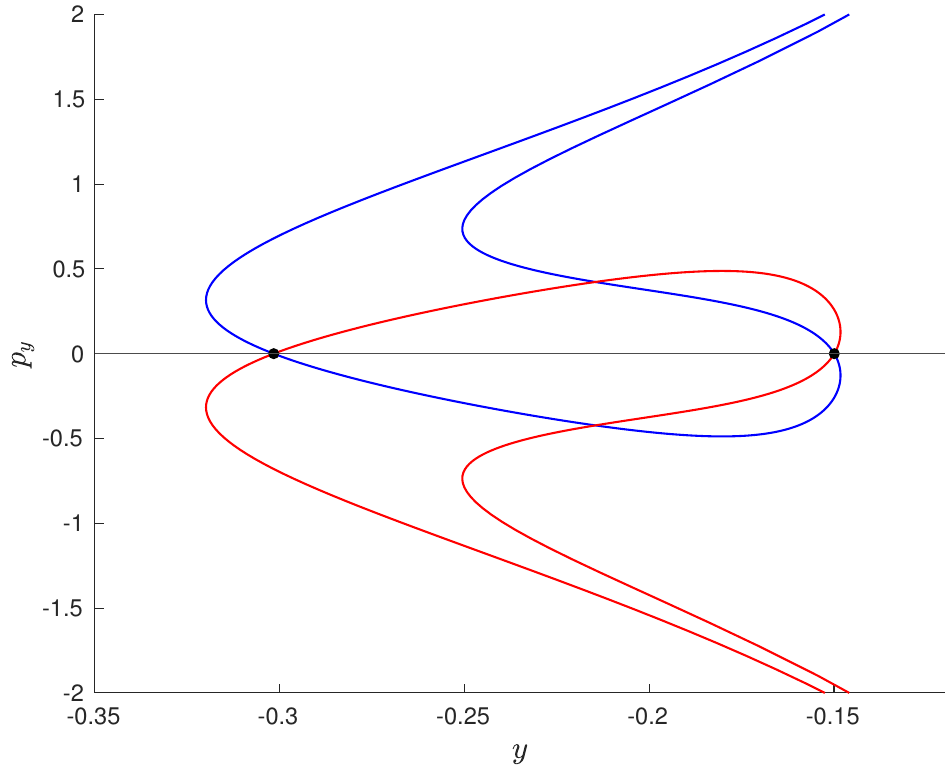}}
	\subfigure[Heteroclinic orbit $\Phi^t_0(z_1)$.]{\includegraphics[scale=0.25]{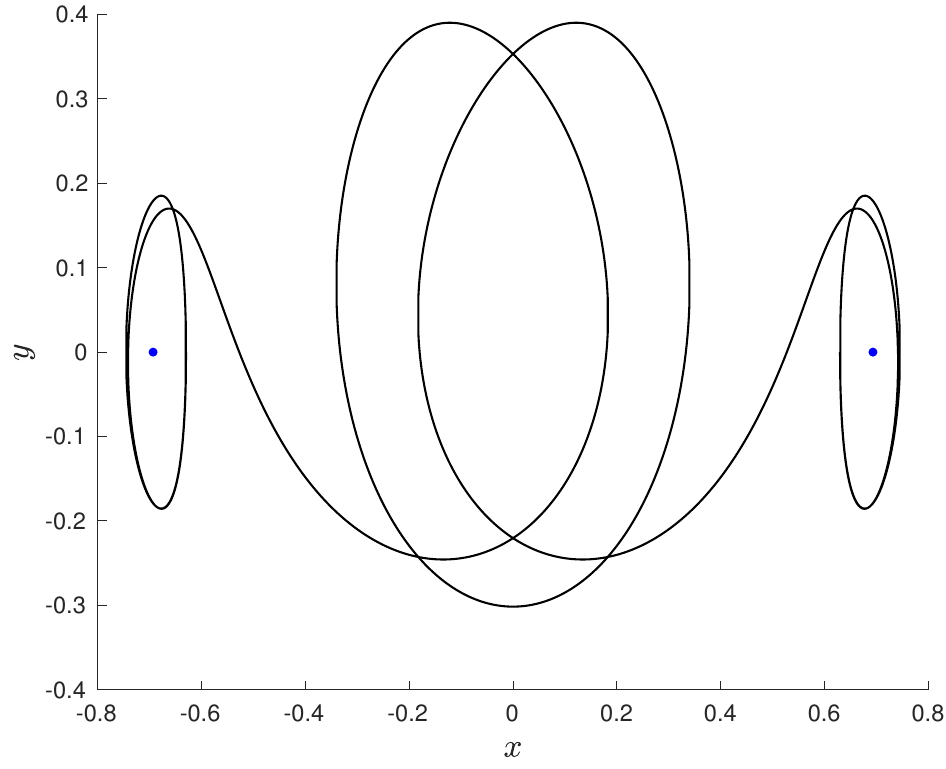}}
	\subfigure[Heteroclinic orbit $\Phi^t_0(z_2)$.]{\includegraphics[scale=0.25]{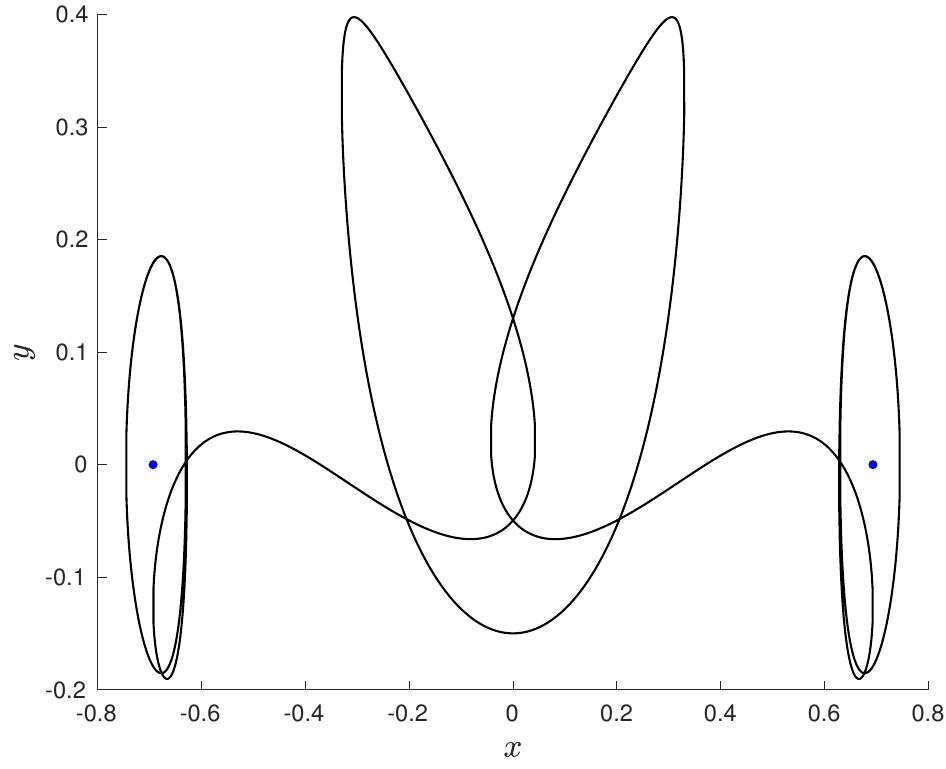}}
    \caption{Heteroclinic points and orbits.\label{fig:het_points_orbits}}
\end{figure}

The manifolds are computed up to the their second intersection with the section $\mathscr{S}_y=\{x=0,\, p_x>0 \}$. On this section, we can find a transverse intersection between the stable and unstable manifold, then we obtain two symmetric heteroclinic points, $z_1$ and $z_2$, from which two heteroclinic orbits are obtained. See Fig.~\ref{fig:het_points_orbits}. Due to the symmetry
$S(x, y, p_x , p_y, t) = (x,-y,-p_x , p_y ,-t)$ there are also   a symmetric heteroclinic connections given by the intersection of $W^u(\lambda_1(h))$ and $W^s(\lambda_2(h))$.
The corresponding heteroclinic points $\hat z_1,\hat z_2$ are the symmetric images under $S$ of $z_1,z_2$, respectively.
These heteroclinic points are given  explicitly in Table~\ref{Cuadro-heteroclinic-points_z_1} and Table~\ref{Cuadro-heteroclinic-points}.

\subsection{The unperturbed scattering map associated to homoclinic and heteroclinic connections}
\label{sec:unperturbed_scattering_hom_het}
\subsubsection{The unperturbed scattering map associated to homoclinic connections}
\label{sec:unperturbed_scattering_hom}
Consider the homoclinic points $z_i=z_i(h)$, $i=1,2$, defined in Section \ref{sec:homoclinic connections}, for the corresponding energy range
$h\in[h_\alpha,h_\beta]$.

\begin{figure}
    \centering
	\includegraphics[scale=0.3]{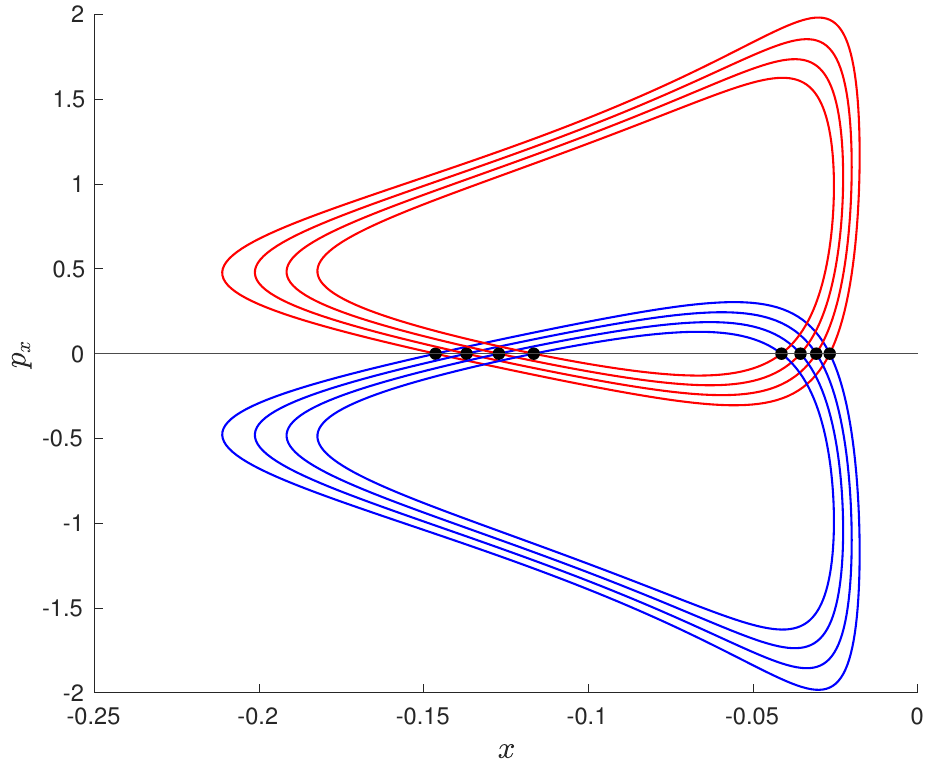}
	\caption{Homoclinic point  $z_2$ (on the left) and $z_1$ (on the right) 
for $x^*(h)=0.615,0.62,0.625,0.63$.\label{all_intersections_hom_inner} }
\end{figure}

For each homoclinic point $z_i$, let $z_i^+\in\Lambda^1_0$ be the foot-point of the stable fiber through $z_i$,
and  $z_i^-\in\Lambda^1_0$ be the foot-point of the unstable fiber through $z_i$.  In $(I,\theta)$-coordinates,
$z_i^+=k_0(I_i,\theta_i^+)$ and $z_i^-=k_0(I_i,\theta_i^-)$.

A sample of  values  $x^*(h)$ and of homoclinic points $z_i=z_i(h)$, for $i=1,2$, are given in Table \ref{Cuadro-homoclinic-points-z1} and Table \ref{Cuadro-homoclinic-points-z2}.
The homoclinic points  are shown in  Fig.~\ref{all_intersections_hom_inner}.
The corresponding foot-points $z_i^\pm$  and their angle coordinates  $\theta^\pm_i$
are given in Table \ref{z1^+_homoc} and \ref{z2^+_homoc}.

\begin{figure}
	\centering
	\subfigure[$\Gamma^{1}_{\textrm{hom}}$.]{\includegraphics[scale=0.35]{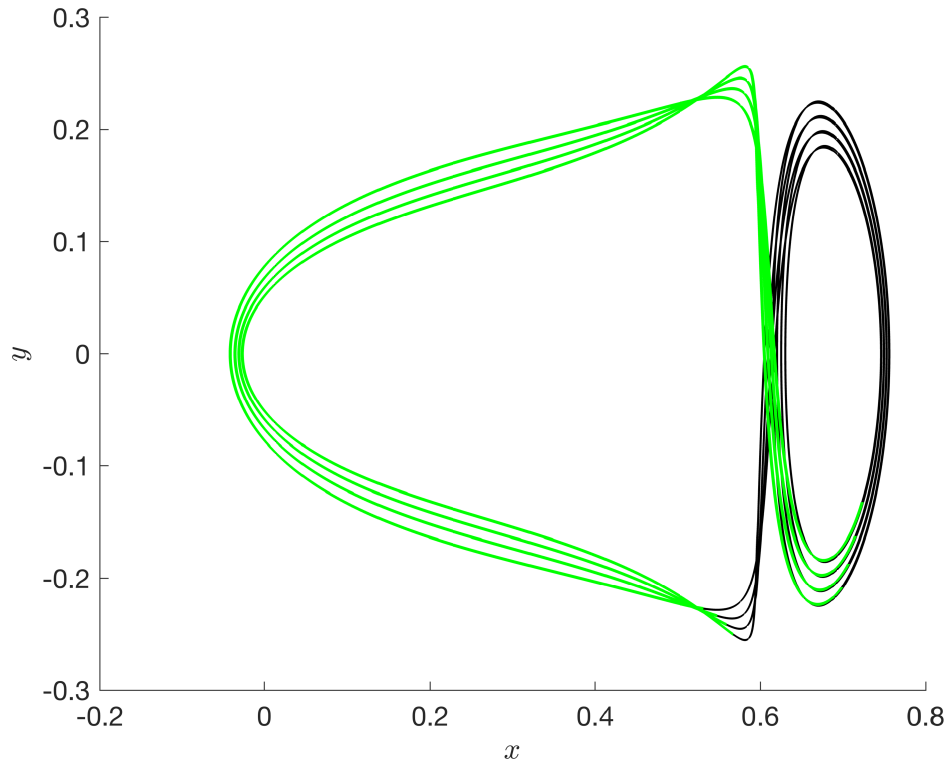}}\hspace{0.001cm}
     \subfigure[$\Gamma^{2}_{\textrm{hom}}$.]{\includegraphics[scale=0.35]{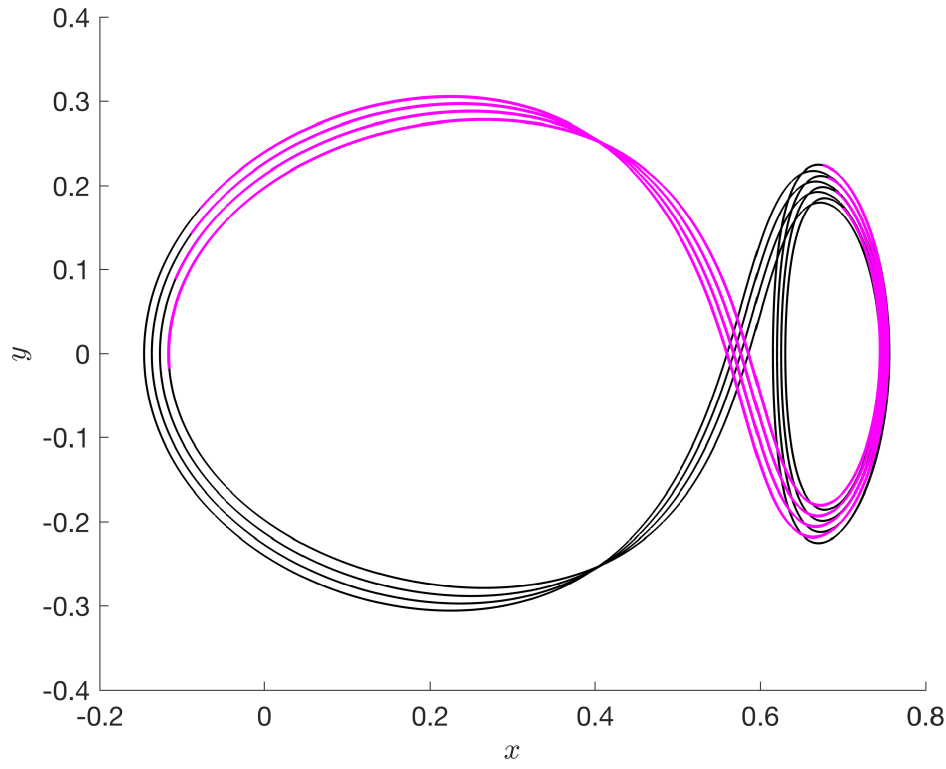}}\hspace{0.001cm}
	\caption{Homoclinic channels $\Gamma^{i}_{\textrm{hom}}$, $i \in\{1,2\}$.}  
\label{Hom_channel_inner}
\end{figure}

As in \eqref{eqn:unperurbed_homoclininc_channel_1}, to the homoclinic point $z_i$ we can associate a homoclinic channel  of the form
\begin{equation}\label{eqn:homoclinic_z2}
\Gamma^i_{\textrm{hom}}=\bigcup_{\substack{ h\in [h_\alpha,h_\beta]\\ \theta\in (\theta_i^-+\theta^*, \theta_i^-+1+\theta^*) }}\Phi_0^{(\theta-\theta_i^-)T} (z_i(h)),
\end{equation}
for $\theta^*\in\mathbb{R}$. By making different choices of $\theta^*$  we  obtain different homoclinic channels and associated  perturbed scattering maps with different properties.
A sample of such homoclinic channels is shown in Fig.~\ref{Hom_channel_inner}.

By Corollary \ref{cor:Delta_symmetric}, the unperturbed scattering map  in $(I,\theta)$-coordinates is a  phase-shift  given  by
\[\sigma^i_0(I,\theta)=(I, \theta+\Delta_i)=(I, \theta+(\theta_i^+-\theta_i^-))=(I, \theta -2\theta_i^-), \textrm{ for } i=1,2.\]
(As in Section \ref{sec:mechanism}, we identify a scattering map $\sigma^i_0$ with its representation $s^i_0$ in the action-angle coordinates.)

\subsubsection{The unperturbed scattering map associated to heteroclinic connections in the CH4BP}
\label{sec:scattering_het}

Consider now a heteroclinic point $z_i=z_i(h)$, $i=1,2$ defined in  Section \ref{sec:heteroclinic_connections}, for the corresponding energy range
$h\in[h_\alpha,h_\beta]$.

\begin{figure}
    \centering
	\includegraphics[scale=0.35]{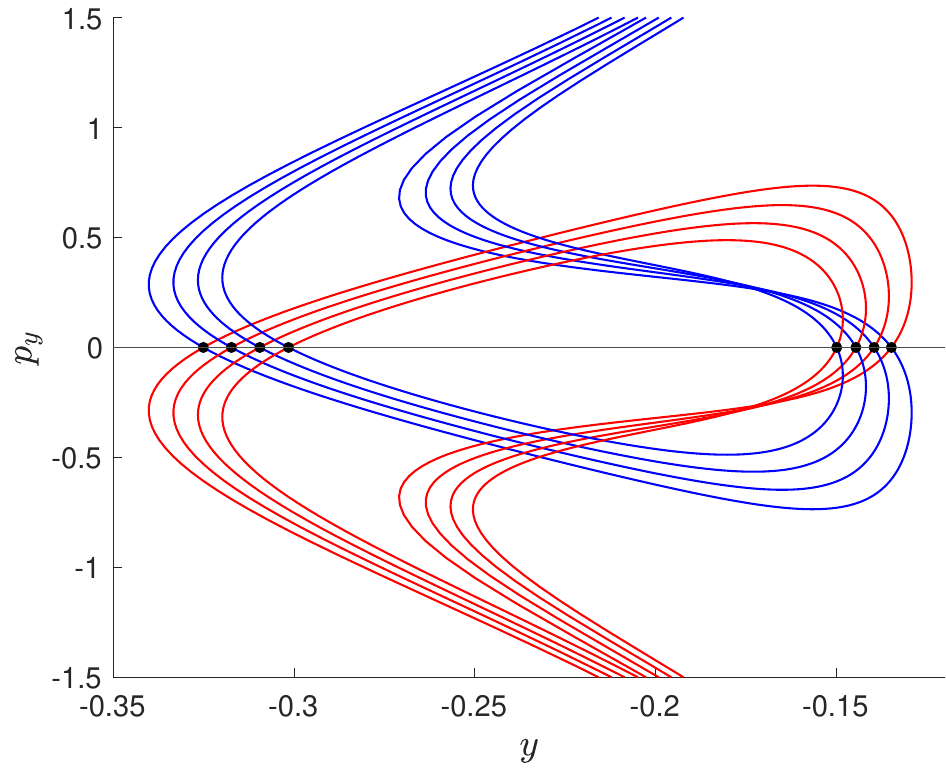}
	\caption{Heteroclinic point  $z_2$ (on the left) and $z_1$ (on the right) for $x^*(h)=0.615,0.62,0.625,0.63$. \label{heteroclinic channels}}
\end{figure}

A sample of values of  $x^*$ and of heteroclinic  points $z_i,\hat z_i$, for $i=1,2$, are given in Table~\ref{Cuadro-heteroclinic-points_z_1} and Table~\ref{Cuadro-heteroclinic-points}.
The heteroclinic points are shown in Fig.~\ref{heteroclinic channels}.
The corresponding foot-points $z_i^\pm$ and their angle coordinates  $\theta^\pm_i$
are given in Table~\ref{z_1^+ hetero} and Table~\ref{z_2^+ hetero}.

\begin{figure}
	\centering
	\subfigure[{$\Gamma^1_{\textrm{het}}$}]{\includegraphics[scale=0.3]{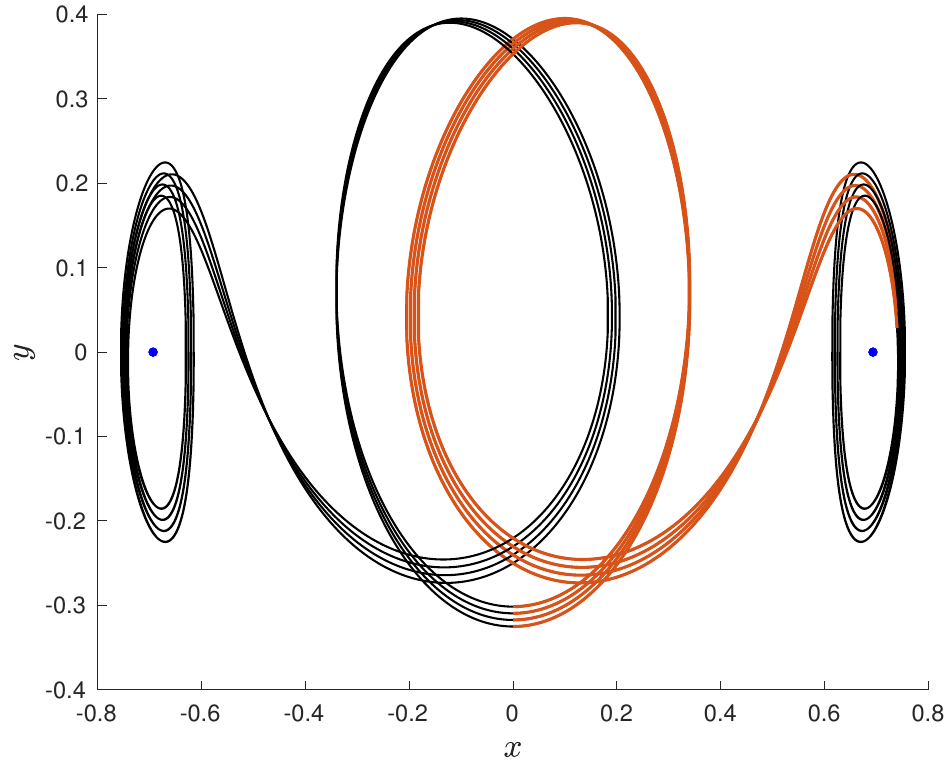}}\hspace{0.001cm}
	\subfigure[{$\Gamma^2_{\textrm{het}}$}]{\includegraphics[scale=0.3]{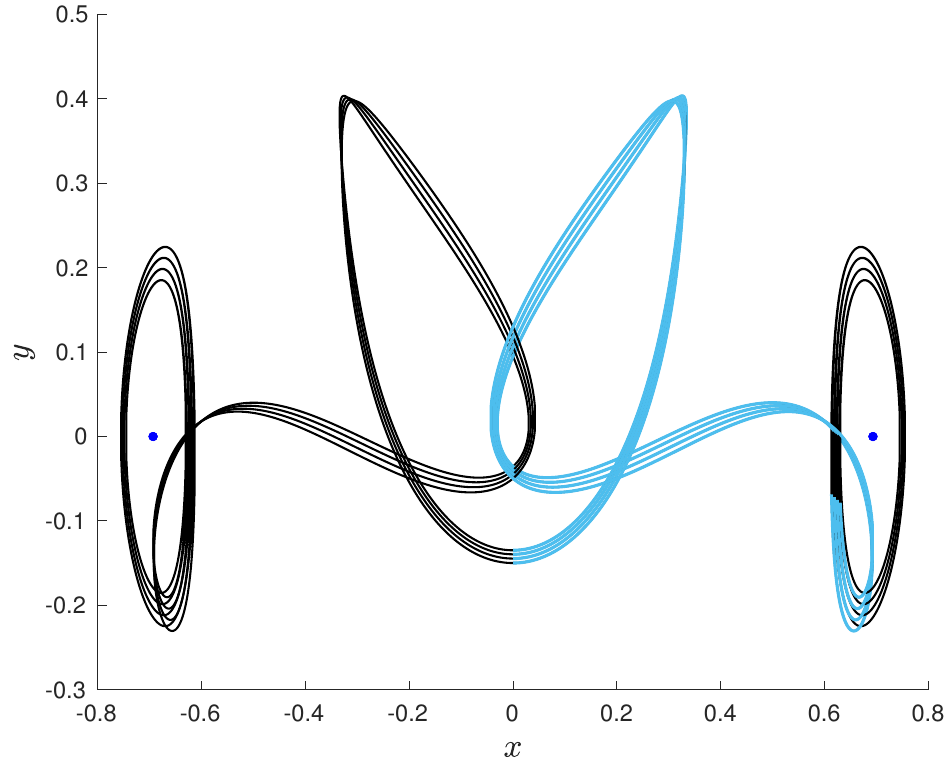}}\hspace{0.001cm}
	\caption{Heteroclinic channels $\Gamma^1_{\textrm{het}}$ and $\Gamma^2_{\textrm{het}}$. \label{Het_channel_inner}}
\end{figure}

We associate  heteroclinic channels of the form
\begin{equation}\label{Het-Channel-inner}
\Gamma^i_{\textrm{het}}=\bigcup_{\substack{ h\in[h_\alpha,h_\beta]\\ \theta\in(\theta_i^-+\theta^*, \theta_i^-+1+\theta^*) }}\Phi_0^{(\theta-\theta_i^-)T} (z_i),
\end{equation}
for some $\theta^*\in \mathbb{R}$.
A sample of such  heteroclinic channels is shown in Fig.~\ref{Het_channel_inner}.

The associated scattering map in $(I,\theta)$-coordinates is a phase-shift  given  by
\[\sigma^i_0(I,\theta)=(I,\theta+\Delta_i)=(I, \theta+(\theta_i^+-\theta_i^-))=(I, \theta - 2\theta_i^-).\]

We stress that in the heteroclinic case the coordinate $\theta_i^-$ is  the  angle coordinate on $\Lambda_0^2$ and
the coordinate $\theta_i^+$ is a the  angle coordinate  on $\Lambda_0^1$, while the value of $I$ is the same in both coordinate systems.


\subsection{Perturbed scattering map associated to homoclinic and heteroclinic connections}
\label{sec:pertubed_scattering_hom_het}
\subsubsection{Perturbed scattering map associated to homoclinic connections}
\label{sec:pertubed_scattering_hom}
Consider an unperturbed  homoclinic channel $\Gamma^i_{\textrm{hom}}$, $i\in\{1,2\}$, defined as in
\eqref{eqn:homoclinic_z2}.
Under the perturbation, $\Gamma^i_{\textrm{hom}}$ persists as the perturbed homoclinic channel $\Gamma^i_{\textrm{hom},\eps}$,
Let $\sigma^i_\eps$ be the associated  scattering map, and $S^i$  the Hamiltonian function that generates it.
We compute $S^i$ in terms of   the (non-symplectic) coordinates $(x^*(h),\theta)$.
using the formula \eqref{eqn:S0-x-star}.  

\begin{figure}
		\centering
		\subfigure[$S^1(x^*,\theta ).$  ]{\includegraphics[scale=0.3]{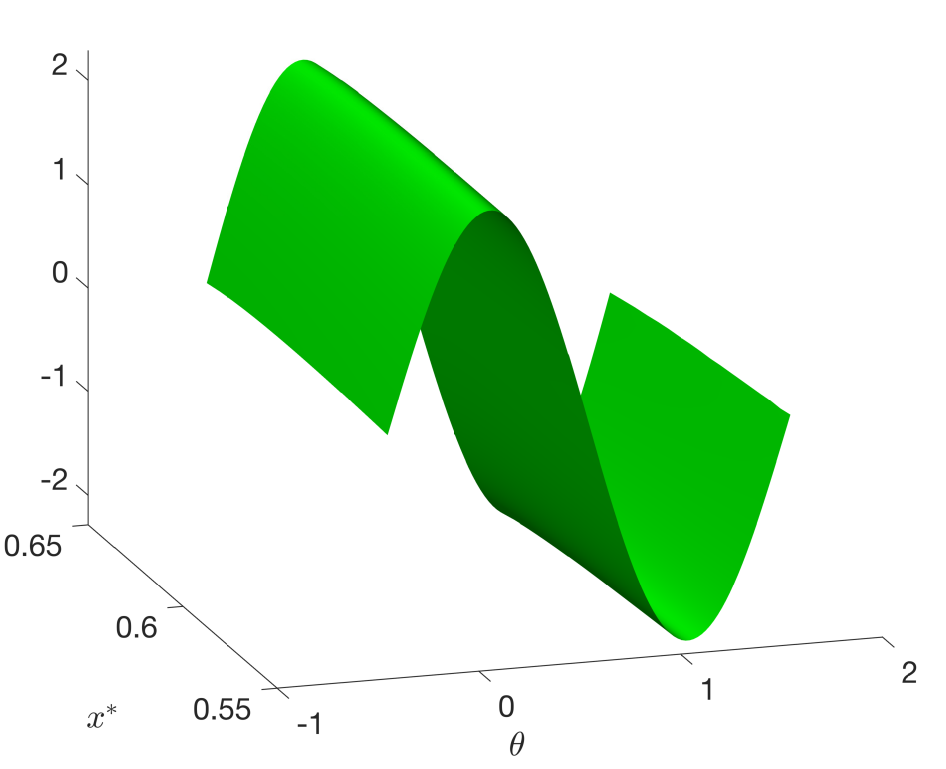}}
		\subfigure[Contour plot of $S^1(x^*,\theta )$. ]{\includegraphics[scale=0.3]{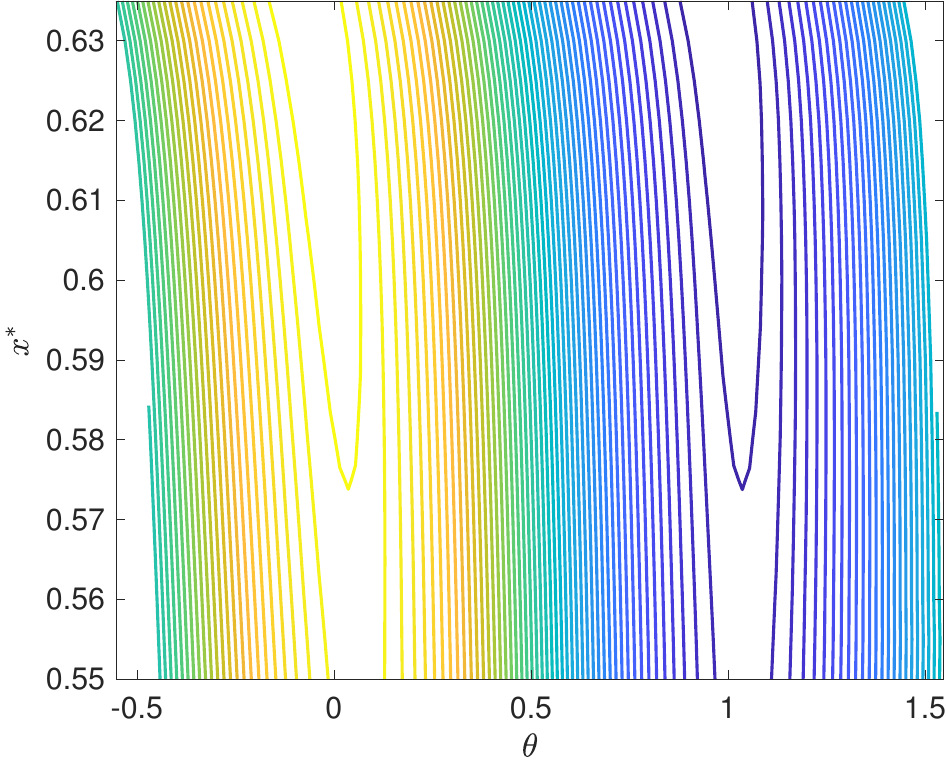}}
		\subfigure[$S^2(x^*,\theta ).$  ]{\includegraphics[scale=0.3]{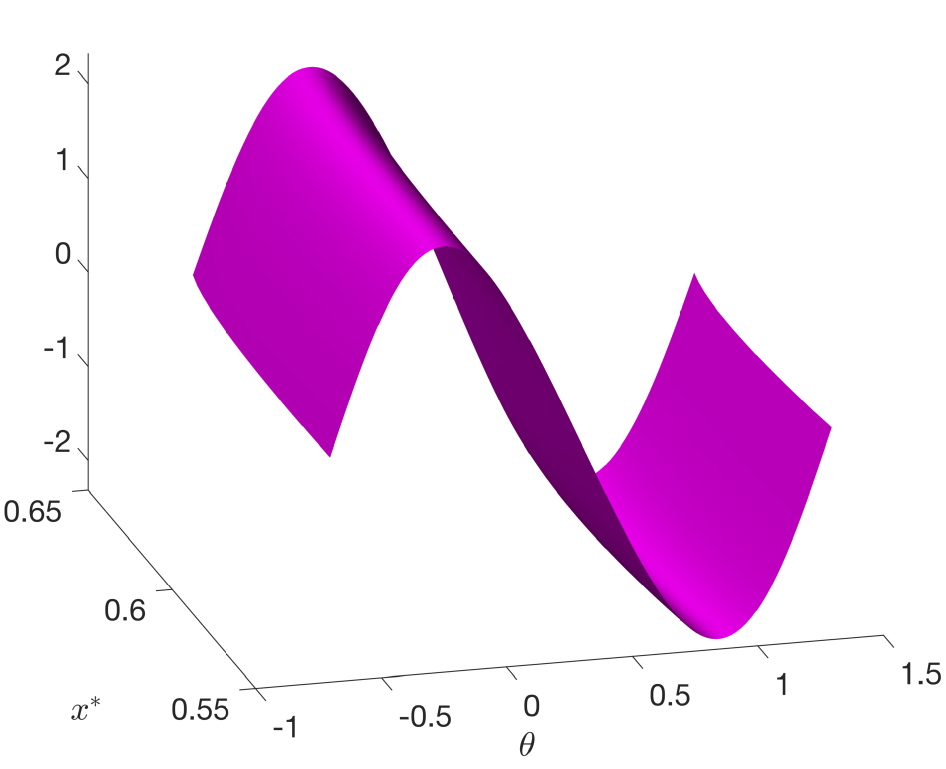}}
		\subfigure[Contour plot of $S^2(x^*,\theta )$. ]{\includegraphics[scale=0.3]{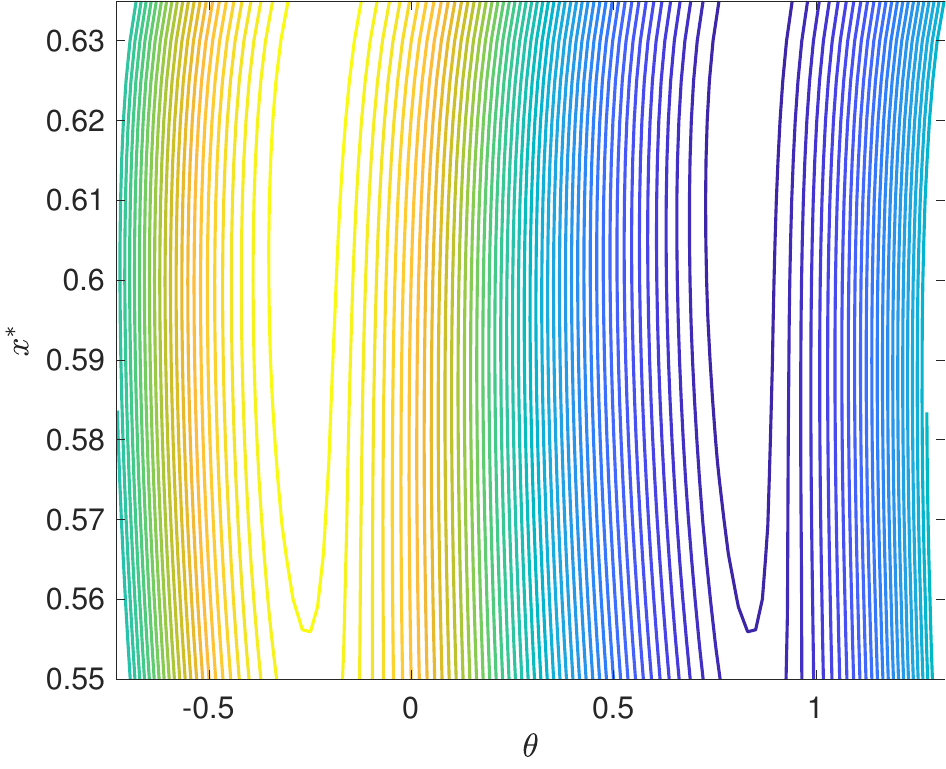}}
		\caption{The plot of $(x^*, \theta)\to S^i(x^*, \theta) $ for $x^*\in[0.55,0.64]$ and $\theta \in[-\theta^-_i(h)-1,-\theta^-_i(h)+1]$ for the homoclinic connections.}\label{3D plot hom}
	\end{figure}

For a range $x^*=x^*(h)\in [0.55, 0.64]$ and $\theta \in [-\theta^-_i(h)-1, -\theta^-_i(h)+1]$ we compute the $3D$-plot and the contour plot of
\[(x^*,\theta)\mapsto S^i(x^*,\theta)\] in Fig.~\ref{3D plot hom}. Here, we take a larger range of $(x^*,\theta)$ in order to obtain a broader view of the scattering map dynamics.

The contour plot represents the level sets of the Hamiltonian $S^i$  in terms of  the (non-symplectic) coordinates $(x^*,\theta)$.
These contour plots are related to the dynamics of the underlying scattering maps $\sigma^i_\eps$ as follows.
An application of $\sigma^i_\eps$ to a point $(I,\theta)$ corresponds, up to an error of $O(\eps^2)$, to an angle shift $(I,\theta+\Delta_i)$, followed by an
$\eps$-step along the contour level at the point $(I,\theta+\Delta_i)$.

\subsubsection{Perturbed scattering map associated to heteroclinic connections}
Consider the heteroclinic connection $W^u(\Lambda_0^2)\cap W^s(\Lambda_0^1)$, and the heteroclinic channels $\Gamma^i_{\textrm{het}}$, $i\in\{1,2\}$, given by
\eqref{Het-Channel-inner}.

Let $\Gamma^i_{\textrm{het},\eps}$ be the perturbed heteroclinic channels.
For each scattering map $\sigma^i_\eps$, the corresponding  Hamiltonian $S^i$ is given by the formula \eqref{eqn:S0-x-star}, where in the case of heteroclinic connection $z_i$ refers to a heteroclinic point in $W^u(\Lambda_0^2)\cap W^s(\Lambda_0^1)$, and   $\theta$ refers to the action-angle coordinate system associated to $\Lambda^{2}_0$.

\begin{figure}
		\centering
		\subfigure[$S^1(x^*,\theta ).$  ]{\includegraphics[scale=0.3]{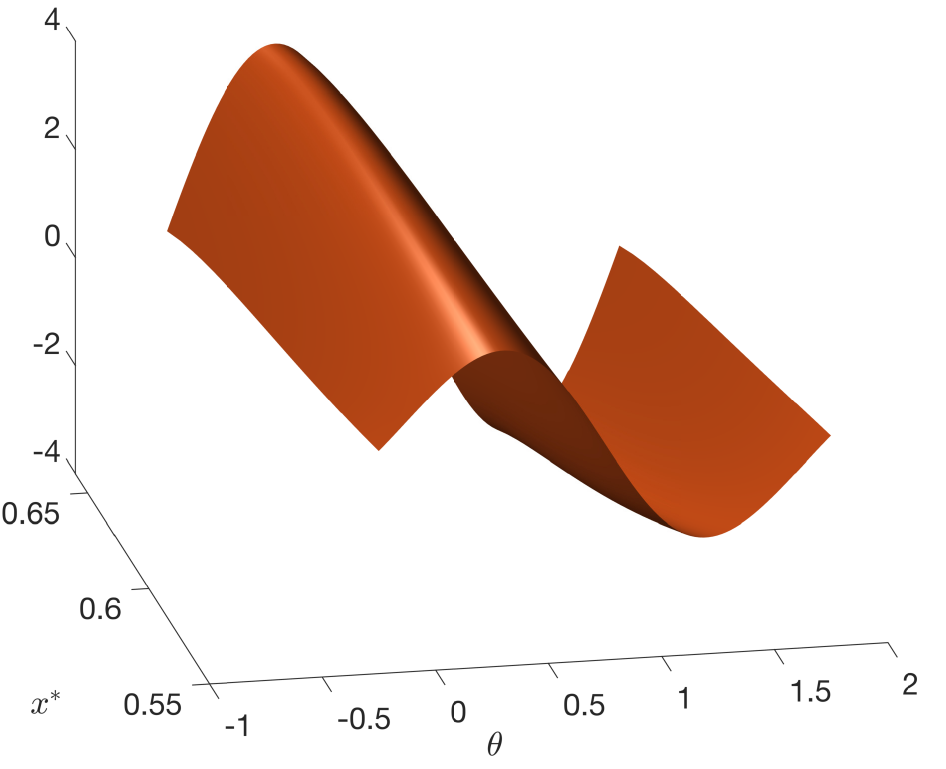}}
		\subfigure[Contour plot of $S^1(x^*,\theta )$. ]{\includegraphics[scale=0.3]{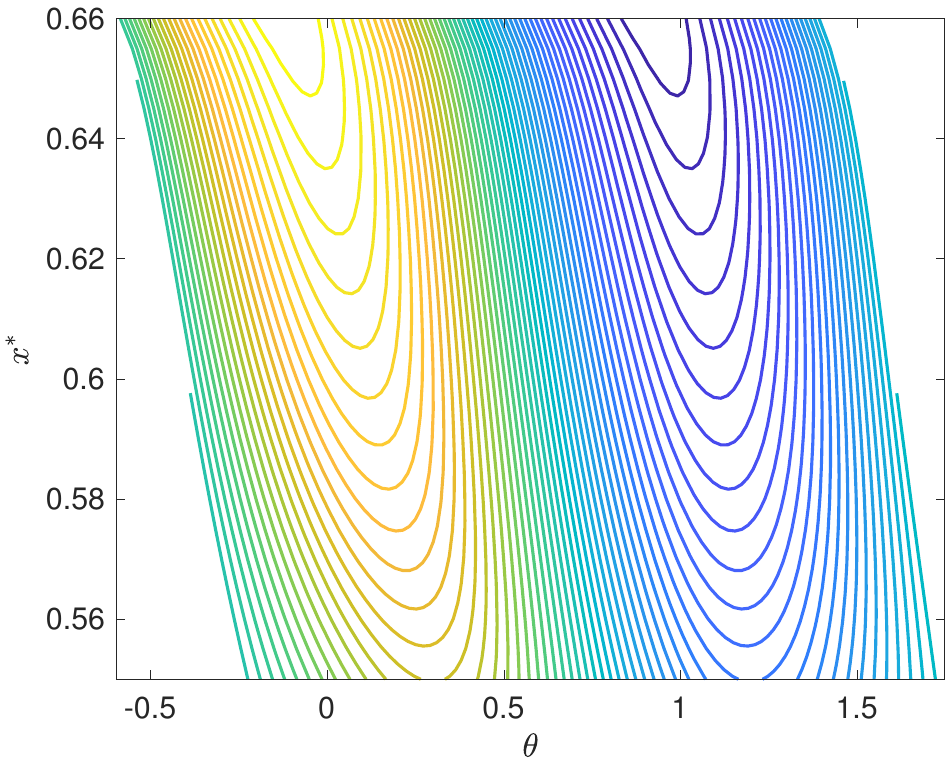}}
		\subfigure[$S^2(x^*,\theta ).$  ]{\includegraphics[scale=0.3]{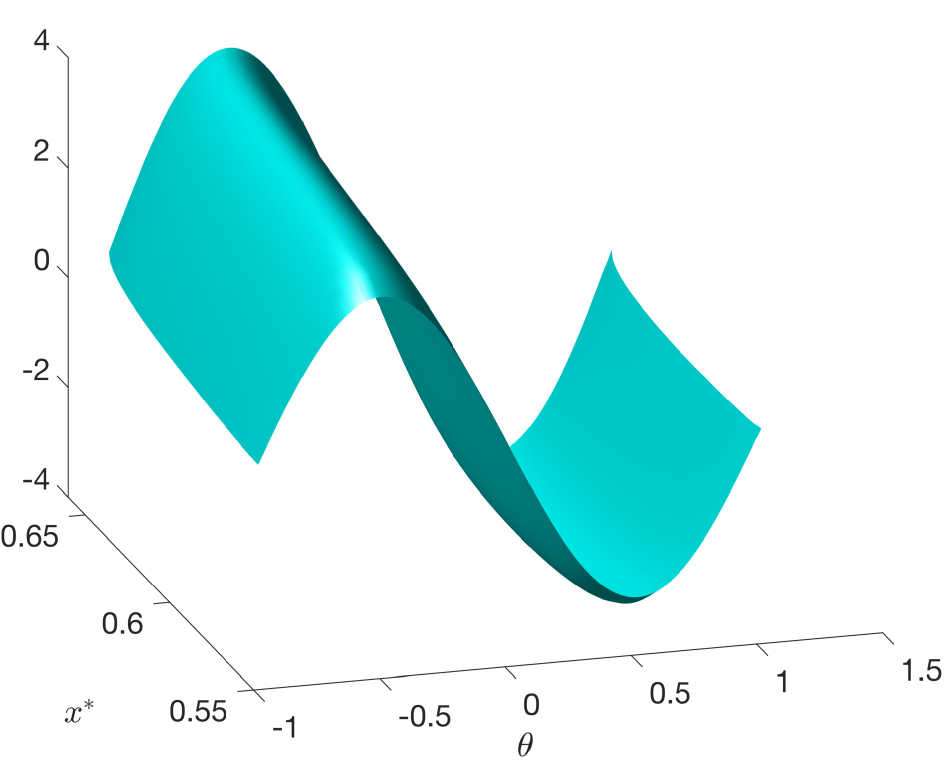}}
		\subfigure[Contour plot of $S^2(x^*,\theta )$. ]{\includegraphics[scale=0.3]{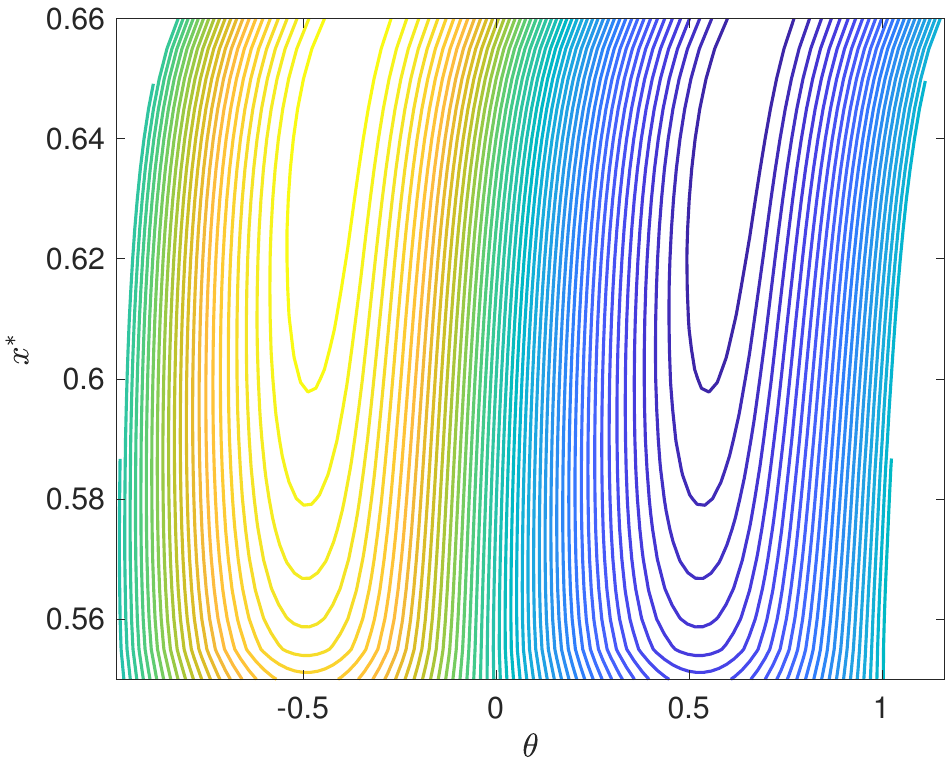}}
		\caption{The plot of $(x^*, \theta)\to S^i(x^*, \theta) $ for $x^*\in[0.55,0.66]$ and $\theta \in[-\theta^-_i(h)-1,-\theta^-_i(h)+1]$ for the heteroclinic connections. }\label{3D plot het}
\end{figure}

As in  the case of homoclinic connections, for  the range $x^*\in [0.55, 0.64]$, we compute the $3D$-plot and the contour plot of the function \[(x^*,\theta)\mapsto S^i(x^*,\theta)\] in Fig.~\ref{3D plot het}. 
The contour plot represents the level sets of the Hamiltonian $S^i$  in terms of  the (non-symplectic) coordinates $(x^*,\theta)$.
The geometric interpretation of these level curves is similar to that for the homoclinic case.

Because of the symmetry $S(x,y,\dot{x},\dot{y},t)=(x,-y,-\dot{x},\dot{y},-t)$,   we also have
a heteroclinic connection $W^u(\Lambda_0^1)\cap W^s(\Lambda_0^2)$, with  heteroclinic channels $\hat\Gamma^i_\textrm{het}$ of the form \eqref{Het-Channel-inner}, for $i=1,2$.
The  scattering map $\hat\sigma^i_\eps$  has the corresponding  Hamiltonian $\hat{S}^i$ given by the  formula \eqref{eqn:S0-x-star},   where  $\theta$ refers to the  angle coordinate  associated to
$\Lambda^{1}_0$. The expression of the scattering map $\hat\sigma^i_\eps$ in the action-angle coordinates on $\Lambda^1_0$  is formally the same as
the expression of the scattering map $\sigma^i_\eps$ in the action-angle coordinates on $\Lambda^2_0$, and the corresponding  Hamiltonian $\hat{S}^i$ and $S^i$ are given by identical formulas.

Thus, one can travel  along the heteroclinic connection $W^u(\Lambda_0^2)\cap W^s(\Lambda_0^1)$,
by using  one of the scattering maps $\sigma^1_\eps$, $\sigma^2_\eps$, and then
then travel along the heteroclinic connection $W^u(\Lambda_0^1)\cap W^s(\Lambda_0^2)$ by using  one of the scattering maps $\hat{\sigma}^1_\eps$, $\hat{\sigma}^2_\eps$. Formally, this is equivalent to applying the formulas for the scattering maps $\sigma^1_\eps$ and $\sigma^2_\eps$ repeatedly, in any order.

\subsection{Numerical verification of the conditions of Theorem \ref{prop:mechanism}}
\label{sec:verification_mechanism}

We verify  the conditions of Theorem~\ref{prop:mechanism} for both homoclinic and heteroclinic orbits.

\subsubsection{Numerical verification of Theorem~\ref{prop:mechanism} in the homoclinic case}
\begin{figure}
		\centering
			\subfigure[$\Gamma_{\text{hom}}$ defined in \eqref{homoclinic channel_1}.]{\includegraphics[scale=0.3]{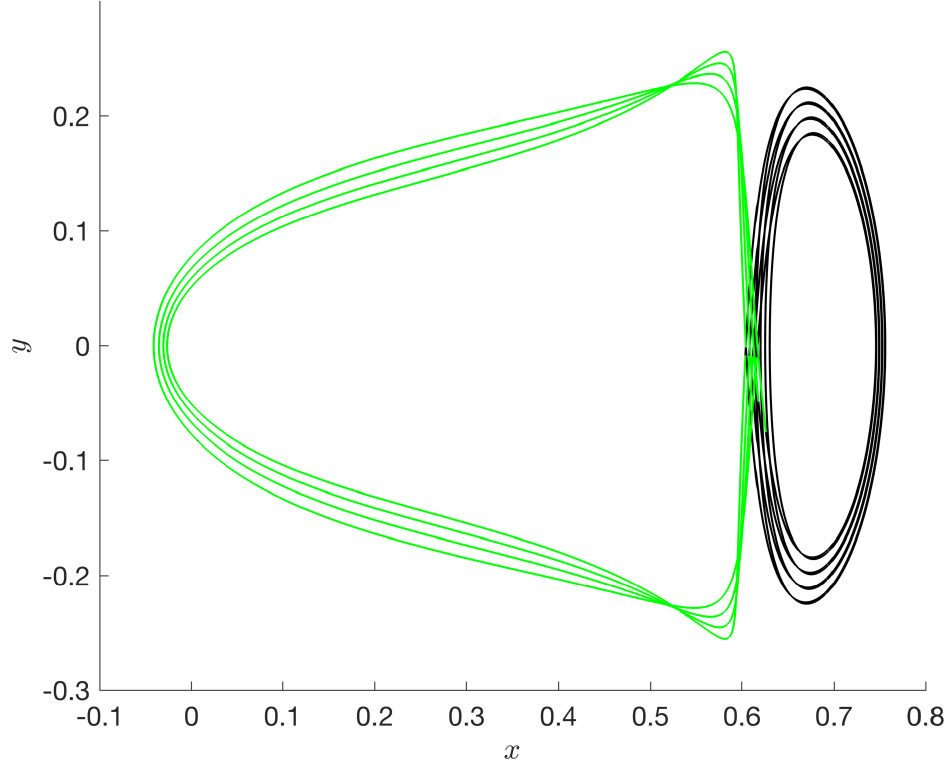}}
             \subfigure[$-\dfrac{\partial S^{1}}{\partial \theta}(x^*,\theta+\Delta_1 (x^*))$.]{\includegraphics[scale=0.3]{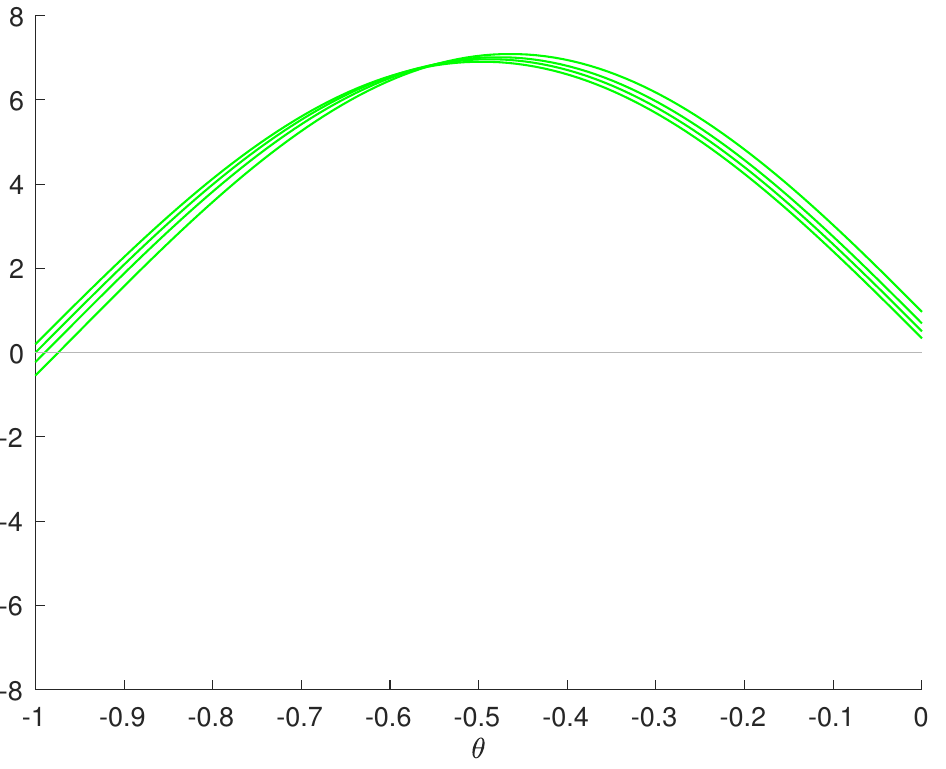}}
		\caption{The plots of $\theta \to -\dfrac{\partial S^{1}}{\partial \theta}(x^*,\theta+\Delta_1 (x^*))$ for the homoclinic channel defined in \eqref{homoclinic channel_1}, i.e, for $\theta\in \left(-1,0\right)$.}\label{dSdtheta_Method1}
	\end{figure}

For $x^*\in [0.615,0.63]$, we select the homoclinic point $z_1(x^*)$, and the
unperturbed homoclinic channel
	\begin{equation}\label{homoclinic channel_1}
		\Gamma^1_{\text{hom}}=\bigcup_{\substack{ x^*\in [0.615,0.63]\\ \theta\in (-1,0) }}\Phi_0^{(\theta-\theta_1^-)T(x^*)} (z_1).
\end{equation}
This is shown in Fig.~\ref{dSdtheta_Method1}.

Let $\Gamma^1_{\text{hom},\eps}$ be the perturbed homoclinic channel, $\sigma^1_\eps$
the associated scattering map, and $S^1$   the Hamiltonian function that generates $\sigma^1_\eps$.
The scattering map  $\sigma^1_\eps$ is defined on $ [0.615,0.63]\times(-1,0)$, that is, on an annuls with a vertical line removed.
This shows that the  condition \eqref{eqn:domains} is satisfied.
The scattering map is of the form as in \eqref{eqn:many_scatterings}.

We  compute $-\dfrac{\partial S^{1}}{\partial \theta}(x^*,\theta+\Delta_1 (x^*))$, whose plot is shown in Fig.~\ref{dSdtheta_Method1}.
The values of $\displaystyle\int_{0}^{1} -\dfrac{\partial S^1}{\partial \theta}(x^*,\theta+\Delta_1(x^*))\, d\theta $ for a sample of values of $x^*\in[0.615,0.63]$ is shown in Table~\ref{tab:B_hom}.
As the integrals are larger than $4.50$,   we have that
\[\int_{\substack{ x^*\in [0.615,0.63]\\ \theta\in (-1,0) }}-\frac{S^1}{d\theta}(x^*,\theta+\Delta_1(x^*))dx^*\wedge d\theta>C_{1,\textrm{hom}}=0.015\cdot 4.50=0.0675.\]

Then the condition \eqref{eqn:sum_integrals} follows. (Since the above integral computed in $(x^*,\theta)$ coordinates is positive, then the  integral in \eqref{eqn:sum_integrals} computed in $(I,\theta)$ coordinates is positive.)

Hence, the conditions of Theorem~\ref{prop:mechanism}  are verified,  implying the existence of pseudo-orbits of the scattering map along which the energy grows.
Theorem \ref{teo:main} yields the existence of diffusing  orbits.
%
%
	
\subsubsection{Numerical verification of Theorem \ref{prop:mechanism} in the heteroclinic case}

\begin{figure}
	\centering
	\subfigure[$\Gamma_{\text{het}}$ defined in \eqref{heteroclinic channel_1}.]{\includegraphics[scale=0.35]{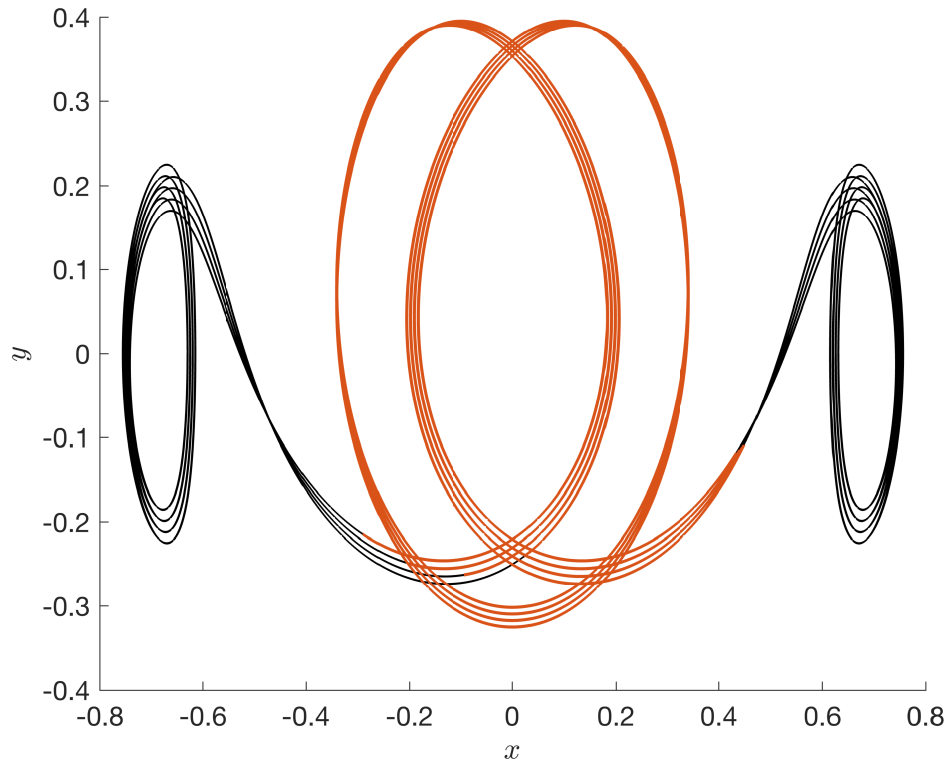}}
	\subfigure[$-\dfrac{\partial S^{1}}{\partial \theta}(x^*,\theta+\Delta_1 (x^*))$.]{\includegraphics[scale=0.35]{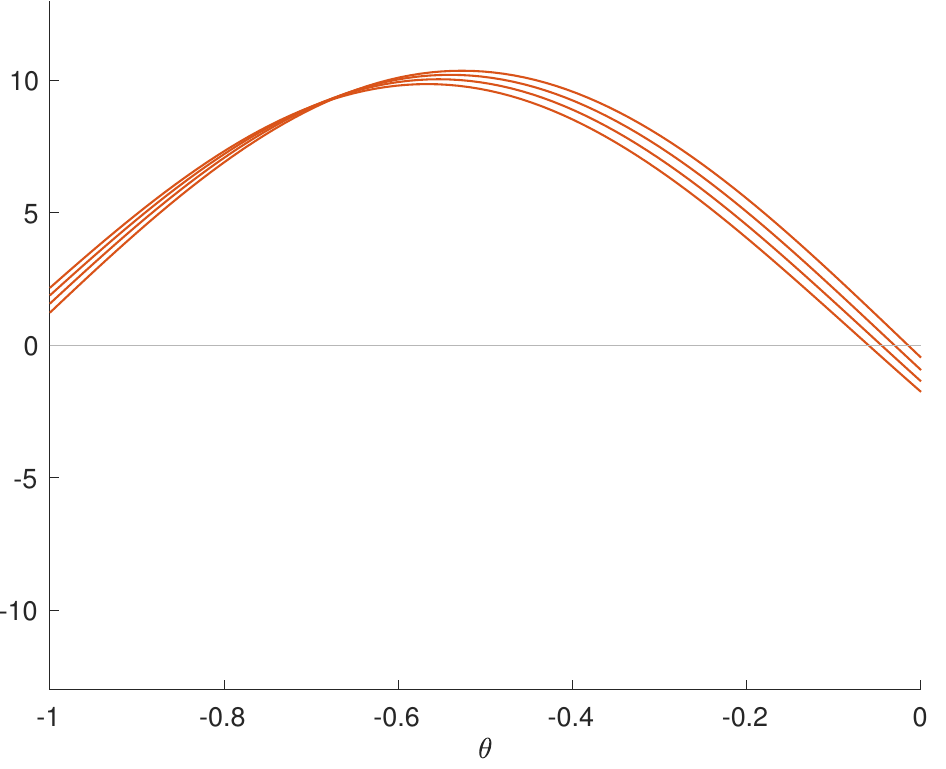}}
	\caption{The plots of $\theta \to -\dfrac{\partial S^{1}}{\partial \theta}(x^*,\theta+\Delta_1 (x^*))$ for the heteroclinic channel defined in \eqref{heteroclinic channel_1}, i.e., for $\theta\in \left(-1,0\right)$.}\label{dSdtheta_Method1_het}
\end{figure}

For $x^*\in [0.615,0.63]$, we consider the heteroclinic point $z_1(x^*)$, and the
heteroclinic  channel
	\begin{equation}\label{heteroclinic channel_1}
		\Gamma^1_{\text{het}}=\bigcup_{\substack{x^*\in [0.615,0.63]\\ \theta\in (-1,0) }}\Phi_0^{(\theta-\theta_1^-)T(x^*)} (z_1).
\end{equation}
This is shown in Fig.~\ref{dSdtheta_Method1_het}.

Let  $\Gamma^1_{\text{het},\eps}$ be the perturbed heteroclinic  channel, $\sigma^1_\eps$ the associated scattering map, and
$S^1$ the Hamiltonian function that generates it.
We  compute $-\dfrac{\partial S^{1}}{\partial \theta}(x^*,\theta+\Delta_1 (x^*))$, whose plot is shown in Fig.~\ref{dSdtheta_Method1_het}.
%

The values of $\displaystyle\int_{0}^{1} -\dfrac{\partial S^1}{\partial \theta}(x^*,\theta+\Delta_1(x^*)) d\theta $ for a sample of values of $x^*\in[0.615,0.63]$ is shown in Table~\ref{tab:B_het}.
As the integrals are larger than $6.18$,
we have that
\[\int_{\substack{ x^*\in [0.615,0.63]\\ \theta\in (-1,0) }}-\frac{S^1}{d\theta}(x^*,\theta+\Delta_1(x^*))dx^*\wedge d\theta>C_{1,\textrm{het}}=0.015\cdot 6.18=0.0927.\]
This implies the condition \eqref{eqn:sum_integrals}.

Hence,  the conditions of Theorem~\ref{prop:mechanism}  are verified, implying the existence of pseudo-orbits  of the scattering map  along which the energy grows.
Theorem \ref{teo:main} yields the existence of diffusing  orbits.

\subsection{Numerical verification of the conditions of Theorem \ref{th:mechanism-main}}
\label{sec:mechanism-main}
We verify  the conditions of Theorem \ref{th:mechanism-main} for both homoclinic and heteroclinic orbits.

\subsubsection{Numerical verification of Theorem \ref{th:mechanism-main} in the homoclinic case}

\begin{figure}
	\subfigure[$\Gamma^1_\text{hom}$ defined in
\eqref{homoclinic channel_12}]{\includegraphics[scale=0.3]{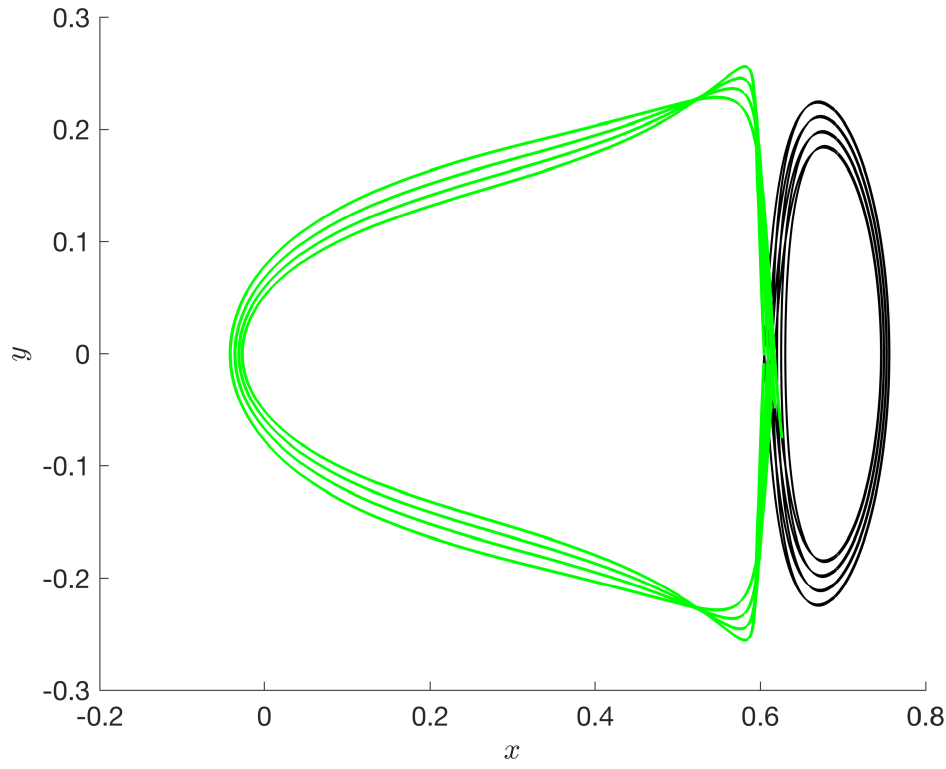}}\hspace{0.5cm}
	\subfigure[$-\dfrac{dS^1}{d\theta}(x^*, \theta+\Delta_1(x^*))$ and line  $\theta=-1$ and $\theta=0$ in blue.]{\includegraphics[scale=0.3]{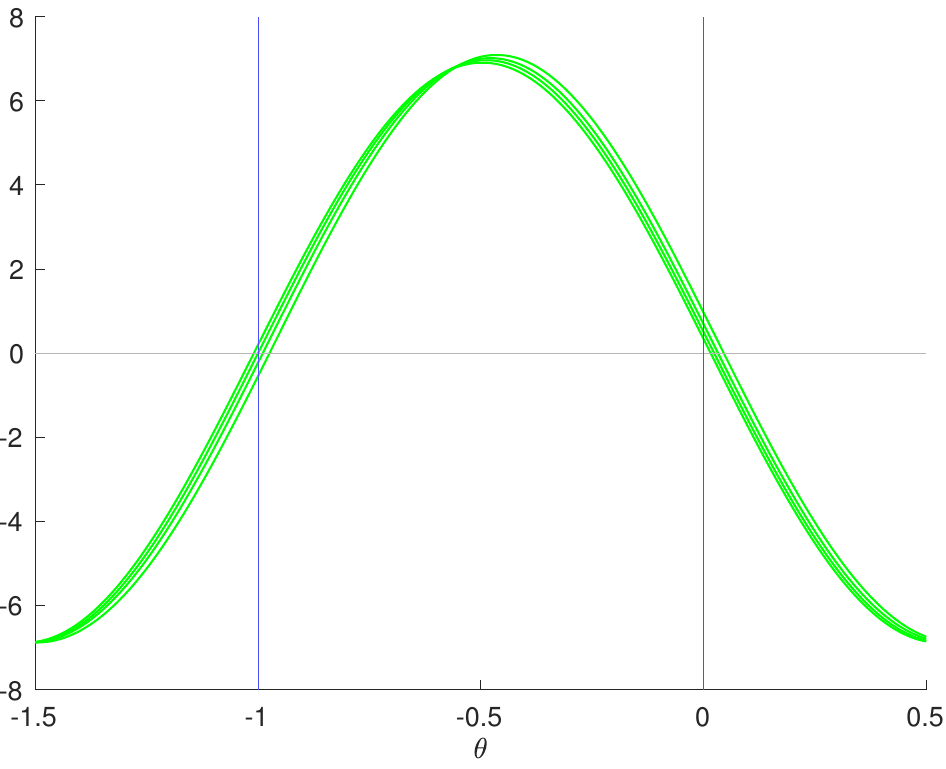}}	
	\subfigure[$\Gamma^2_\text{hom}$  defined in
\eqref{homoclinic channel_12} ]{\includegraphics[scale=0.3]{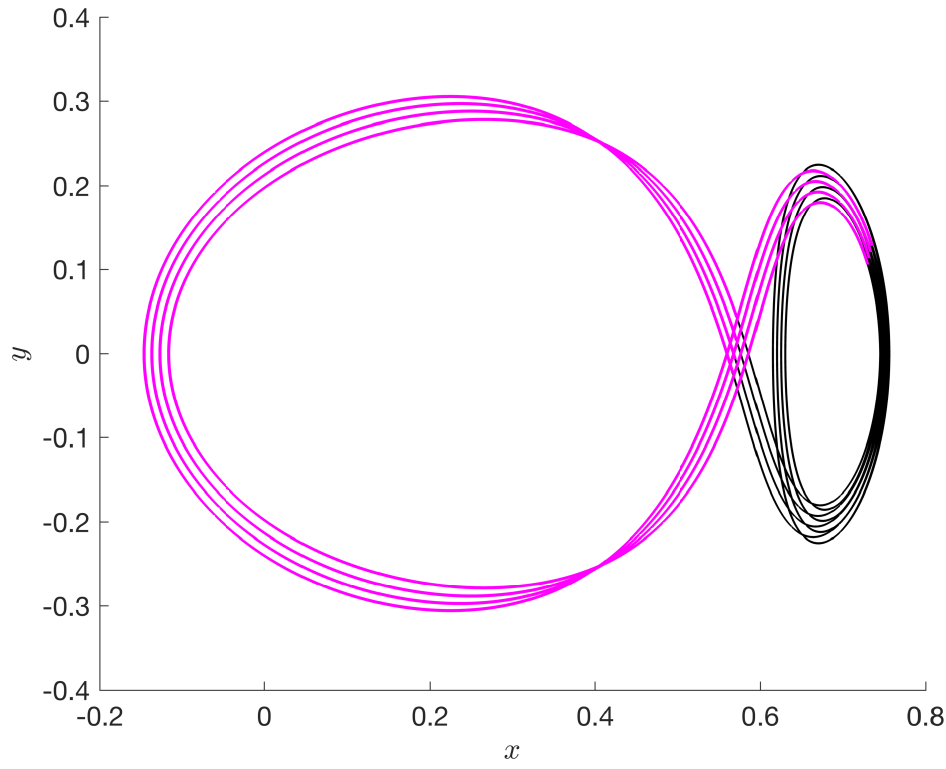}}\hspace{0.5cm}
	\subfigure[$-\dfrac{dS^2}{d\theta}(x^*, \theta+\Delta_2(x^*))$ and line  $\theta=-0.6$ and $\theta=0.4$ in red.]{\includegraphics[scale=0.3]{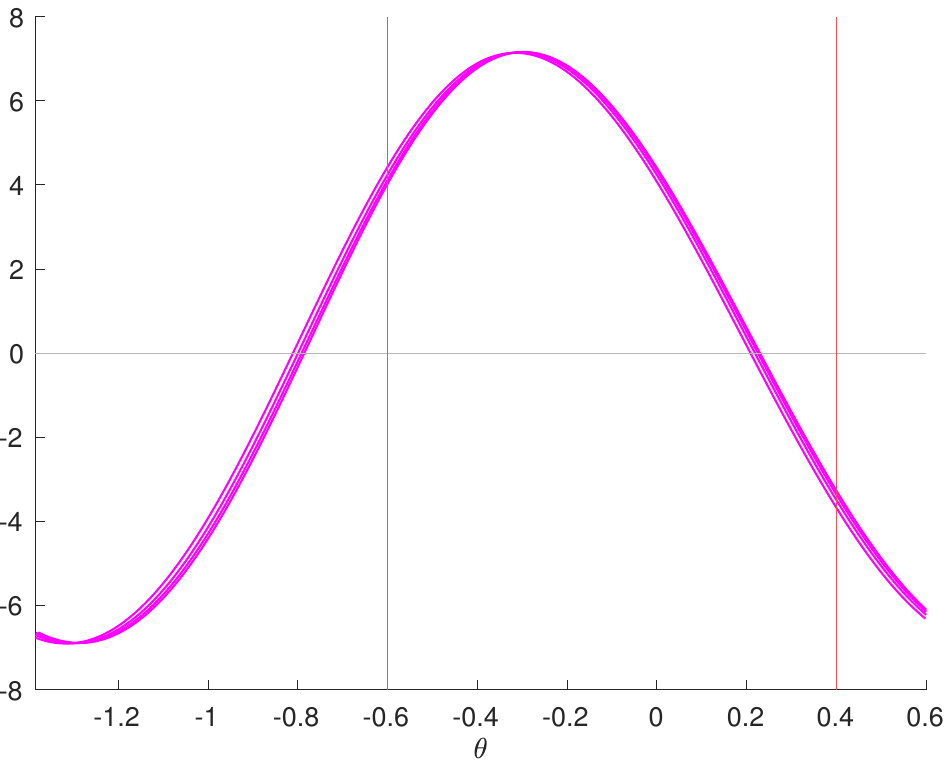}}
\caption{The  homoclinic channels $\Gamma^i_\text{hom}$ defined in
\eqref{homoclinic channel_12}, and the plot of $-\dfrac{dS^i}{d\theta}(x^*, \theta+\Delta_i(x^*))$, for $i=1,2$. The  graph
associated to $\Gamma^1_\text{hom}$ is between the lines in blue and the  graph
associated to $\Gamma^2_\text{hom}$ is between the lines in red.}
\label{fig:Method 2 hom channel}
\end{figure}

For $x^*\in [0.615,0.63]$, we consider the homoclinic points $z_1(x^*)$ and $z_2(x^*)$ and the corresponding homoclinic channels
\begin{equation}\label{homoclinic channel_12}\begin{split}
		\Gamma^{1}_{\text{hom}}=&\bigcup_{\substack{ x^*\in [0.615,0.63]\\ \theta\in (-1,0) }}\Phi_0^{(\theta-\theta_1^-)T(x^*)} (z_1),\\
        \Gamma^{2}_{\text{hom}}=&\bigcup_{\substack{ x^*\in [0.615,0.63]\\ \theta\in (-0.6,0.4) }}\Phi_0^{(\theta-\theta_1^-)T(x^*)} (z_2)
\end{split}
\end{equation}
See Fig. ~\ref{fig:Method 2 hom channel}.

Let $\Gamma^{1}_{\text{hom},\eps}$, $\Gamma^{2}_{\text{hom},\eps}$ be the perturbed homoclinic channels,   $\sigma^1_\eps$, $\sigma^2_\eps$   the associated scattering maps, and $S^1,S^2$   the corresponding  Hamiltonian functions.

\begin{figure}
		\centering
		\includegraphics[scale=0.3]{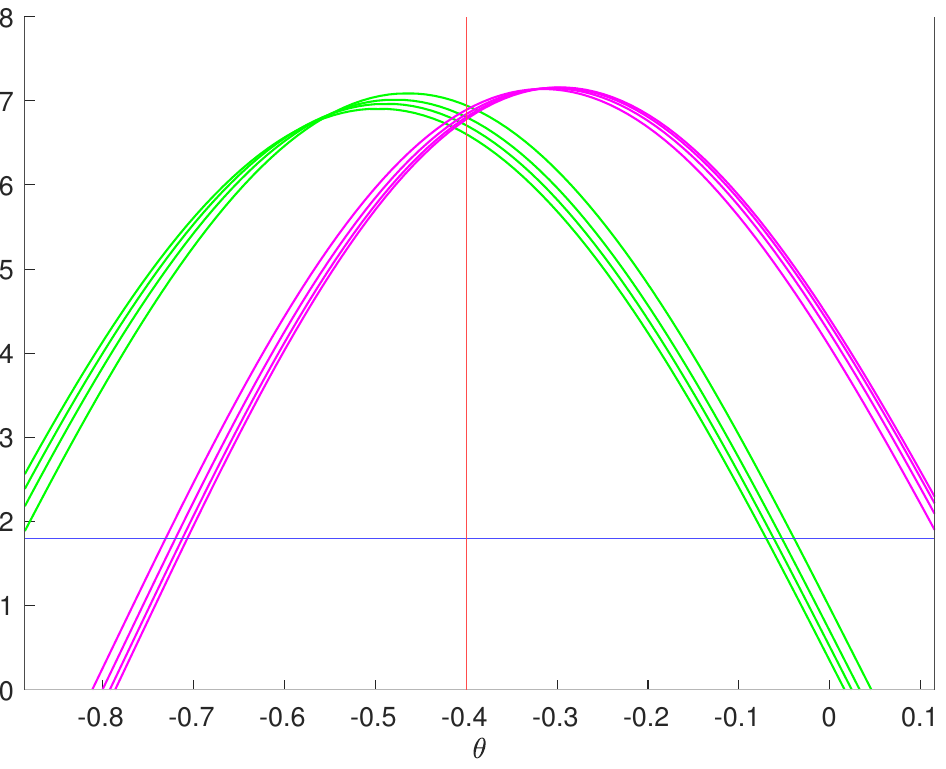}
		\caption{The plots of $-\dfrac{\partial S^{i}}{\partial \theta}(x^*,\theta+\Delta_i(x^*))$, $i=1$ in green and $i=2$ in pink, for $\theta\in [-0.885, 0.115]$. The vertical line $\theta=-0.4$ in orange and the horizontal line $c=1.8$ in blue.} \label{dSdtheta_Method2_zoom}
	\end{figure}

We consider the interval $\theta\in[-0.885,0.115]$ of length $1$. We verify numerically that for all $x^*\in [0.615,0.63]$ we have
\begin{itemize}
\item for $\theta\in [-0.885,-0.4]$ we have $-\dfrac{\partial S^{1}}{\partial \theta}(x^*,\theta+\Delta_1(x^*))>  c_{\textrm{hom}} = 1.8$,
\item for $\theta\in [-0.4, 0.115]$ we have $-\dfrac{\partial S^{2}}{\partial \theta}(x^*,\theta+\Delta_2(x^*))>  c_{\textrm{hom}} = 1.8$.
\end{itemize}
See Figure~\ref{dSdtheta_Method2_zoom}.

This shows that the condition \eqref{eq:scatter-grad-a1} of Theorem~\ref{th:mechanism-main} is verified,  implying the existence of pseudo-orbits of the system of  scattering maps along which the energy grows.
Theorem \ref{teo:main} yields the existence of diffusing  orbits.

\subsubsection{Numerical verification of Theorem \ref{th:mechanism-main} in the heteroclinic case}

\begin{figure}
		\subfigure[$\Gamma^1_\text{het}$  defined in \eqref{heteroclinic channel_12}]{\includegraphics[scale=0.3]{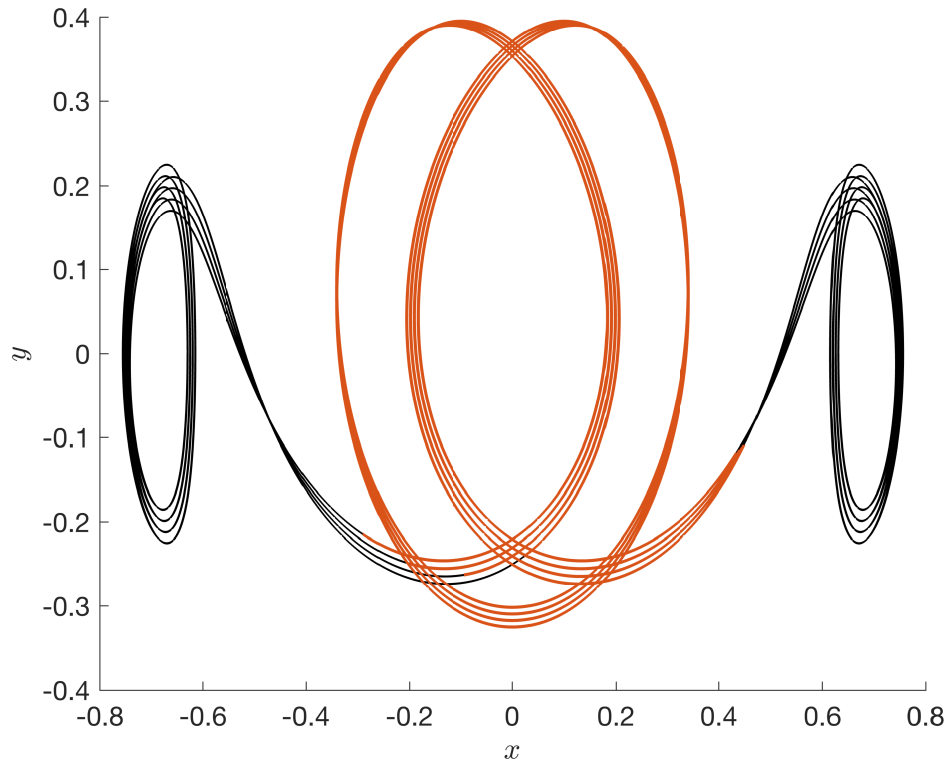}}\hspace{0.5cm}
		\subfigure[$-\dfrac{dS^1}{d\theta}(x^*, \theta+\Delta_1(x^*))$ and line  $\theta=-1$ and $\theta=0$ in blue.]{\includegraphics[scale=0.3]{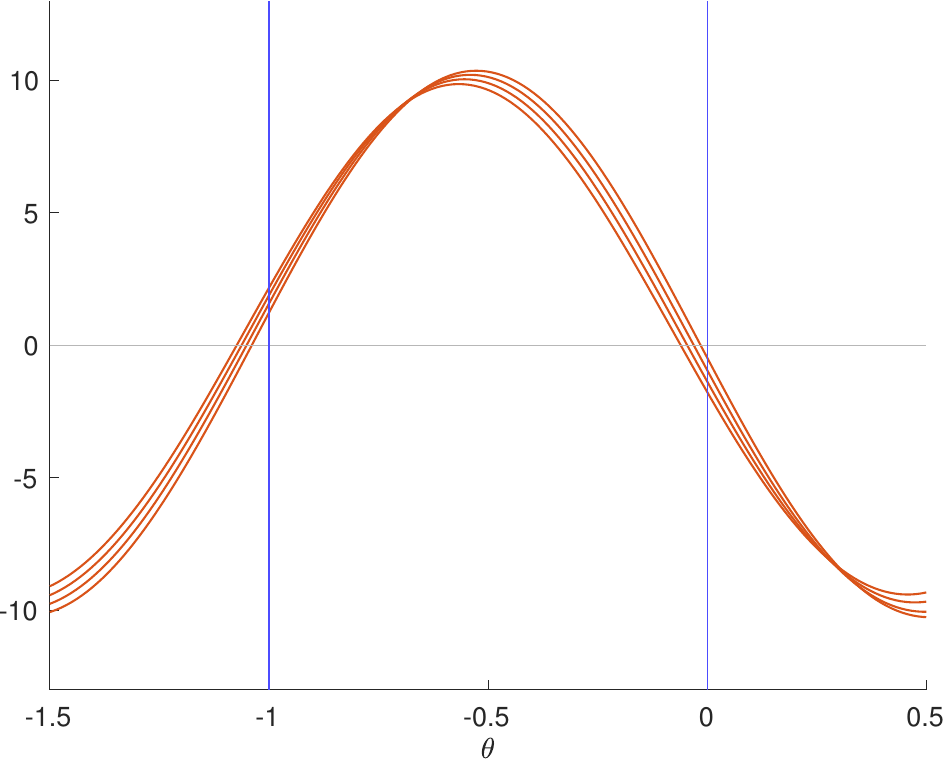}}

			\subfigure[$\Gamma^2_\text{het}$  defined in  \eqref{heteroclinic channel_12} ]{\includegraphics[scale=0.3]{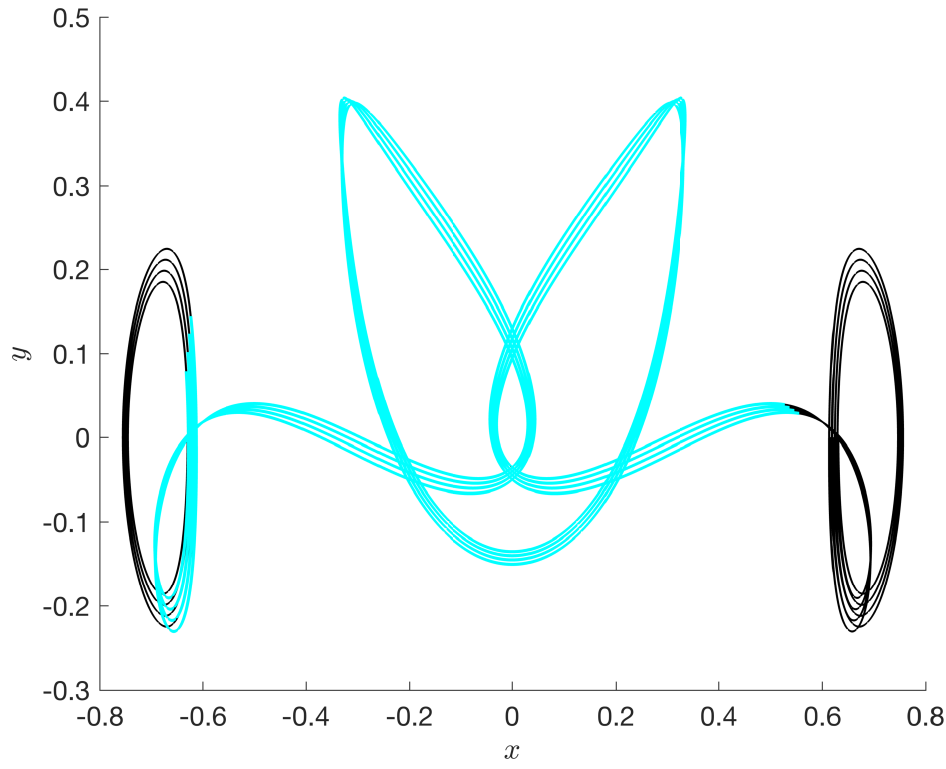}}\hspace{0.5cm}
		\subfigure[$-\dfrac{dS^2}{d\theta}(x^*, \theta+\Delta_2(x^*))$ and line  $\theta=-0.6$ and $\theta=0.4$ in red.]{\includegraphics[scale=0.3]{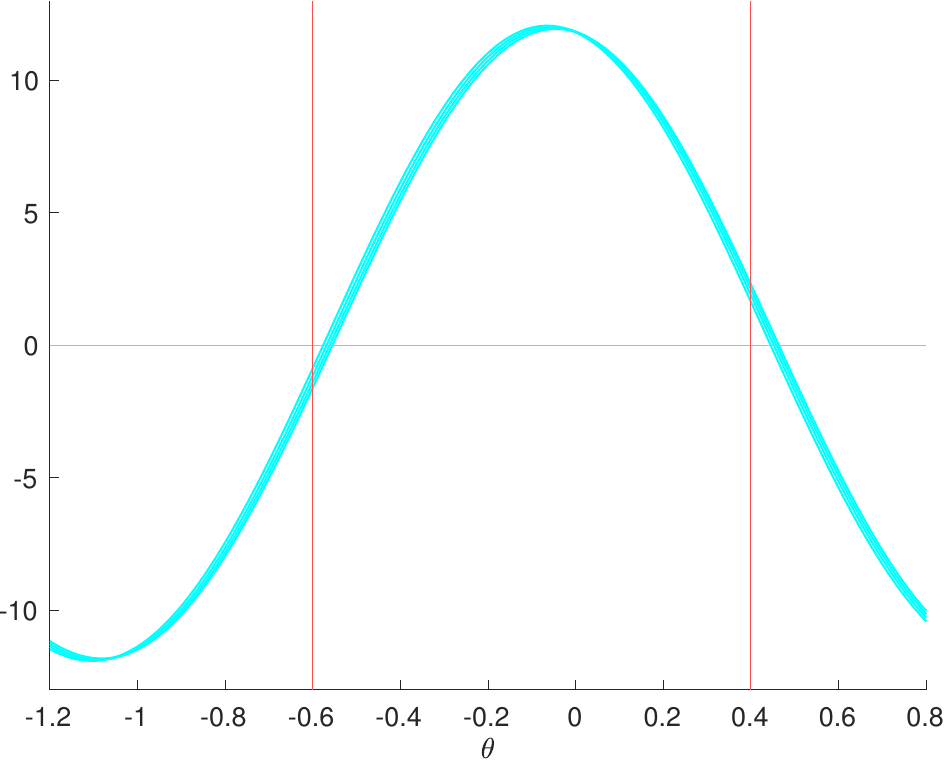}}
\caption{The  heteroclinic channels $\Gamma^i_\text{het}$ defined in  \eqref{heteroclinic channel_12}, and the plot of $-\dfrac{dS^i}{d\theta}(x^*, \theta+\Delta_i(x^*))$, for $i=1,2$.
The graph
associated to $\Gamma^1_\text{het}$ is between the lines in blue and the  graph
associated to $\Gamma^2_\text{het}$ is between the lines in red.}
\label{fig:Method_2_het channel}
	\end{figure}

For $x^*\in [0.615,0.63]$, we consider the heteroclinic points $z_1(x^*)$ and $z_2(x^*)$ and the corresponding heteroclinic channels
\begin{equation}\label{heteroclinic channel_12}
\begin{split}
		\Gamma^{1}_{\text{het}}=&\bigcup_{\substack{ x^*\in [0.615,0.63]\\ \theta\in (-1,0) }}\Phi_0^{(\theta-\theta_1^-)T(x^*)} (z_1),\\
        \Gamma^{2}_{\text{het}}=&\bigcup_{\substack{ x^*\in [0.615,0.63]\\ \theta\in (-0.6,0.4) }}\Phi_0^{(\theta-\theta_1^-)T(x^*)} (z_2)
\end{split}
\end{equation}
See Fig.~\ref{fig:Method_2_het channel}.

Let $\Gamma^{1}_{\text{het},\eps}$, $\Gamma^{2}_{\text{het},\eps}$ be the perturbed homoclinic channels,  $\sigma^1_\eps$, $\sigma^2_\eps$  the associated scattering maps, and $S^1,S^2$  the corresponding  Hamiltonian functions.

\begin{figure}
		\centering
		\includegraphics[scale=0.3]{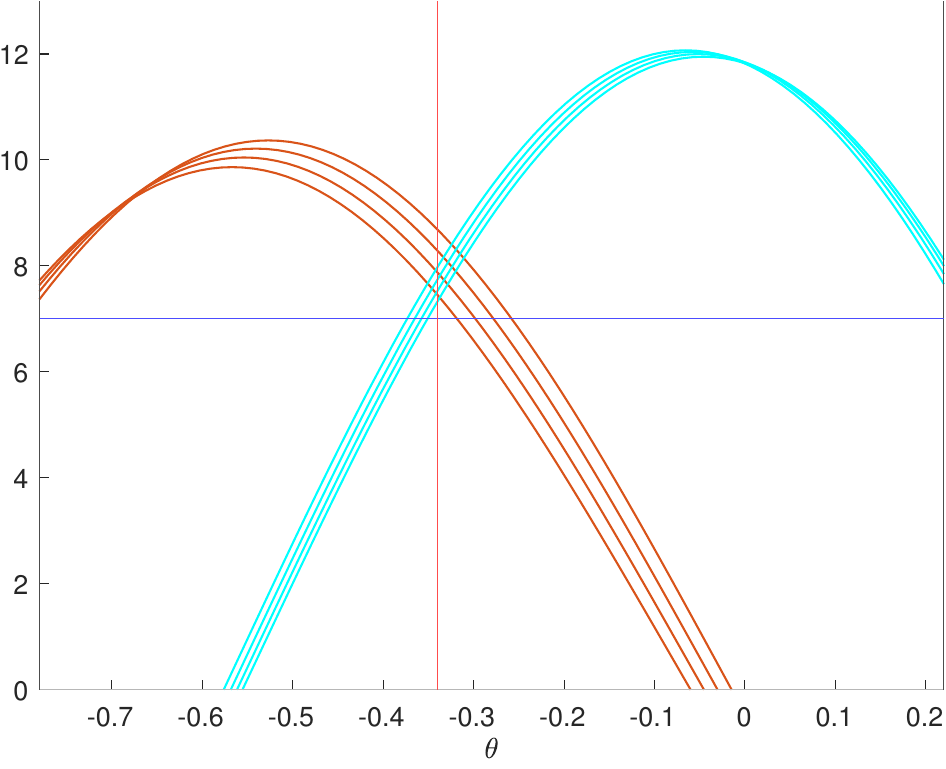}
		\caption{The plots of $-\dfrac{\partial S^{i}}{\partial \theta}(x^*,\theta+\Delta_i(x^*))$, $i=1$ in orange and  $i=2$ in sky blue, for $\theta\in [-0.78, 0.22]$. The vertical line $\theta=-0.34$ in red and the horizontal line $c=7$ in blue.} \label{dSdtheta_Method2_zoom_jhet}
	\end{figure}

We consider the interval $\theta\in[-0.78,0.22]$ of length $1$. We verify numerically that for all $x^*\in [0.615,0.63]$ we have
\begin{itemize}
\item for $\theta\in [-0.78,-0.34]$ we have $-\dfrac{\partial S^{1}}{\partial \theta}(x^*,\theta+\Delta_1(x^*))> c_{\textrm{het}} = 7$,
\item for $\theta\in [-0.34, 0.22]$ we have $-\dfrac{\partial S^{2}}{\partial \theta}(x^*,\theta+\Delta_2(x^*))> c_{\textrm{het}} = 7$.
\end{itemize}
See Figure~\ref{dSdtheta_Method2_zoom_jhet}.
This shows that the condition~\eqref{eq:scatter-grad-a1} of Theorem~\ref{th:mechanism-main} is verified, implying the existence of pseudo-orbits of the system of  scattering maps along which the energy grows.
Theorem \ref{teo:main} yields the existence of diffusing  orbits.
	
\begin{rem}
Since in Section~\ref{sec:verification_mechanism} we have $C_{1,\textrm{hom}}=0.0675 < C_{1,\textrm{het}}=0.0927$,  we see that the mechanism of diffusion based on the Birkhoff Ergodic Theorem is more efficient when using heteroclinic orbits compared to homoclinic orbits.
Similarly, since  $c_{\textrm{hom}}=1.8< c_{\textrm{het}}=7$, we see  that mechanism of diffusion based on multiple scattering maps is more efficient when using heteroclinic orbits compared to homoclinic orbits.
\end{rem}

\clearpage

\appendix
\section{Tables}


\begin{table}[h!]
	\centering	
	{	\begin{eqnarray*}
			\begin{array}{|l| l |  l| l|}
				\hline
			h	&	 \pm	x^*(h)	& \pm p_y(h)& T(h)		  \\ \hline \hline
			-2.07715457 		& \pm 0.615 &\pm	0.48123483 	 & 	  3.05967299	\\  \hline
			-2.08671819	& \pm0.62 & \pm	0.45248801 	 &		3.05678890 		\\  \hline
			-2.09582648 		& \pm	0.625&  \pm 0.42353750 & 3.05406407 		\\  \hline
			-2.10446079	&	\pm	0.63& \pm 0.39437305   & 3.05150060 	\\  \hline
			\end{array}
	\end{eqnarray*}}
	\caption{Initial conditions, periods and energy levels of the Lyapunov periodic orbits shown in Fig.~\ref{fig:families_periodic_orbits}.}\label{Cuadro-initial-conditions-1}
\end{table}


\begin{table}[h]
	\centering	
	{	\begin{eqnarray*}
			\begin{array}{| l| l|
}
				\hline
				x^*	(h)& z_1(h)	  \\ \hline \hline
				0.615    &   (-0.14646739,0,0,  -3.09272039 )                       \\  \hline
				0.62& (-0.13707337 ,0,0,-3.23629663)                               \\  \hline
				0.625&  (-0.12723761 ,0,0,	-3.40227909)                          \\  \hline
				0.63& (-0.11666113 ,0,0,-3.60215892 )                             \\  \hline
			\end{array}
	\end{eqnarray*}}
	\caption{Homoclinic point   $z_1$. \label{Cuadro-homoclinic-points-z1}}
\end{table}


\begin{table}[h]
	\centering	
	{	\begin{eqnarray*}
			\begin{array}{|l|  l|  }
				\hline
				x^*	(h)& z_2(h)	                                 \\ \hline \hline
				0.615 & ( -0.02675921 ,0,0,-8.40169252 )           \\  \hline
				0.62  & (-0.03085457 ,0,0, -7.78778643   )        \\  \hline
				0.625 & (-0.03563430 ,0,0, -7.20679133 )          \\  \hline
				0.63  & (-0.04141244 ,0,0, -6.64009686  )      	\\  \hline
			\end{array}
	\end{eqnarray*}}
	\caption{Homoclinic point   $z_2$. \label{Cuadro-homoclinic-points-z2}}
\end{table}


\begin{table}[h]
	\centering	
	{	\begin{eqnarray*}
			\begin{array}{|  l|  l|    l|}
				\hline
				x^*	(h)& z_1^\pm	(h)	& \theta_1^\pm(h)	 	 \\ \hline \hline
			    0.615 & (0.61500027  , \mp 0.00000005  ,\mp  0.00000100   , 0.48123436 )   & \pm 0.27849625  \\  \hline	
				0.62 &  (0.62000004   ,\mp 0.00000011  ,\mp  0.00000028 ,   0.45248794 )   & \pm 0.28369125  \\  \hline
				0.625 & (0.62500023  ,\mp  0.00000006   ,\mp 0.00000083   , 0.42353712 )	&\pm 0.29106293 \\  \hline
				0.63&   (0.63000040  , \mp 0.00000000   ,\pm 0.00000132  ,  0.39437240 )    &\pm 0.30161988 \\  \hline
			\end{array}
	\end{eqnarray*}}
	\caption{Footpoints $z_1^\pm$  of the homoclinic point $z_1$, and angle coordinates  $\theta_1^\pm$.\label{z1^+_homoc}}
\end{table}


\begin{table}[h]
	\centering	
	{	\begin{eqnarray*}
			\begin{array}{|  l|  l|    l|}
				\hline
				x^*	(h)& z_2^\pm(h)			 & \theta_2^\pm(h) \\ \hline \hline
				0.615 & (0.61499366  , \pm 0.00000076   , \pm0.00001982  , 0.48124537 )  & \pm 0.49677491 \\  \hline
				0.62 &  (0.61999254 ,  \pm 0.00000079   , \pm0.00002317 ,  0.45250024)	 & \pm 0.48835256 \\  \hline
				0.625 & (0.62499943   ,\mp 0.00000216   , \pm0.0000008  ,  0.42353842 )	 & \pm 0.47799047 \\  \hline
				0.63&   (0.63000023   ,\pm 0.00000018   ,\pm 0.0000008   , 0.39437267 )  & \pm 0.46465922 \\  \hline
			\end{array}
	\end{eqnarray*}}
	\caption{Footpoints $z_2^\pm$  of the homoclinic point $z_2$, and  angle coordinates  and $\theta_2^\pm$.\label{z2^+_homoc}}
\end{table}


\begin{table}[h!]
	\centering	
	{	\begin{eqnarray*}
			\begin{array}{|l| l|        }
				\hline
				x^*(h)	& z_1(h),	\hat z_1(h)                                              \\ \hline \hline
				0.615& (0 ,  \mp0.32513154, \pm 1.41324947 , 0)                        \\  \hline
				0.62&  (0  ,\mp0.31736214,\pm 1.45901631 , 0)                          \\  \hline
				0.625& (0 ,\mp0.30951753,\pm 1.50672557 , 0)                            \\  \hline
				0.63& (0 ,\mp0.30157444  ,\pm 1.55664227 , 0)	                             \\  \hline
			\end{array}
	\end{eqnarray*}}
	\caption{Heteroclinic points $z_1, \hat z_1$.\label{Cuadro-heteroclinic-points_z_1}}
\end{table}


\begin{table}[h!]
	\centering	
	{	\begin{eqnarray*}
			\begin{array}{|l| l|  }
				\hline
				x^*(h)& z_2(h),	\hat z_2(h)                               \\ \hline \hline
				0.615& ( 0 , \mp 0.13480207 , \pm  3.26837828 ,   0)        \\  \hline
				0.62&  ( 0,   \mp 0.13957356 ,\pm  3.18684274 , 0)        \\  \hline
				0.625& ( 0 , \mp 0.14460595  ,\pm 3.10468687 ,  0)          \\  \hline
				0.63&  ( 0 ,  \mp 0.14991892 , \pm  3.02186543 , 0)         \\  \hline
			\end{array}
	\end{eqnarray*}}
	\caption{Heteroclinic points $z_2, \hat z_2$.\label{Cuadro-heteroclinic-points}}
\end{table}


\begin{table}[h!]
	\centering	
	{	\begin{eqnarray*}
			\begin{array}{|  l| l|   l|}
				\hline
				x^*	(h)& z_1^\pm(h)  & \theta_1^\pm(h)	            \\ \hline \hline
				0.615& ( \pm0.61500108 , 0.00000001 ,  0.00000360  , \pm0.48123302 )    &  \pm 0.56696772  \\  \hline
				0.62& (  \pm0.62000004 , 0.00000010 ,  0.00000045 ,  \pm0.45248793 )    &  \pm 0.55400246  \\  \hline
				0.625& (\pm0.62500059   , 0.00000015 ,  0.00000215 ,  \pm0.42353653 )   &  \pm 0.54075381   \\  \hline
				0.63& ( \mp0.63000045  ,  0.00000009 ,  0.00000163  , \mp0.39437232 )   &  \pm 0.52702772   \\  \hline
			\end{array}
	\end{eqnarray*}}
	\caption{Footpoints $z_1^\pm$  of the heteroclinic point  $z_1$ and angle coordinates $\theta_1^\pm$.
		\label{z_1^+ hetero}}
\end{table}


\begin{table}[h!]
	\centering	
	{	\begin{eqnarray*}
			\begin{array}{|  l| l|     l|}
				\hline
				x^*	(h)& z_2^\pm(h) 		&	\theta_2^\pm(h)		 \\ \hline \hline
				0.615&  ( \pm0.61499894  , -0.00000004 ,   0.00000329  , \pm0.48123659 )	 &\pm 0.04588521           \\  \hline
				0.62& (   \pm0.61999901  , -0.00000003   , 0.00000316 ,  \pm0.45248962 )     &\pm 0.05160310              \\  \hline
				0.625&(   \pm0.62499909  , -0.00000000   , 0.00000300  , \pm0.42353896 )	 &\pm  0.05811547           \\  \hline
				0.63&  (  \pm0.62999918 ,  -0.00000002   , 0.00000285  , \pm0.39437436 )	 &\pm 0.06560153              \\  \hline
			\end{array}
	\end{eqnarray*}}
	\caption{Footpoints $z_2^\pm$  of  the heteroclinic point $z_2$ and angle coordinates $\theta_2^\pm$.
		\label{z_2^+ hetero}}
\end{table}


	\begin{table}
	{	\begin{eqnarray*}
            \begin{array}{|l|  l| }\hline
			x^*&\displaystyle  \int_{0}^{1} -\dfrac{\partial S^1}{\partial \theta}(x^*,\theta+\Delta_1(x^*)) d\theta \\ \hline\hline
			0.615 &4.504268037535489 \\ \hline
			0.62 & 4.522975396978230\\ \hline
			0.625 & 4.530206829459286 \\\hline
			0.63 &4.546613141589633 \\ \hline
		\end{array}
	\end{eqnarray*}}
		\caption{The integral $ \displaystyle  \int_{0}^{1} -\dfrac{\partial S^1}{\partial \theta}(x^*,\theta+\Delta_1(x^*)) d\theta $ for the homoclinic connections.}
     \label{tab:B_hom}
	\end{table}


\begin{table}
	{	\begin{eqnarray*}
            \begin{array}{|l|  l| }\hline
		x^*&\displaystyle  \int_{0}^{1} -\dfrac{\partial S^1}{\partial \theta}(x^*,\theta+\Delta_1(x^*))\, d\theta \\ \hline\hline
		0.615 & 6.185282908819167 \\ \hline
		0.62 & 6.375516217294260\\ \hline
		0.625 & 6.551536879012402 \\\hline
		0.63 &6.713775974499518 \\ \hline
	\end{array}
	\end{eqnarray*}}
	\caption{The integral $ \displaystyle  \int_{0}^{1} -\dfrac{\partial S^1}{\partial \theta}(x^*,\theta+\Delta_1(x^*))\, d\theta $ for the heteroclinic connections.}
\label{tab:B_het}
\end{table}


\clearpage

\bibliographystyle{alpha}
\bibliography{diffusion,elliptich4bpref}

\end{document}